\documentclass[12pt]{article}
\usepackage{amsmath,amssymb}
\usepackage{amsthm}
\advance\textwidth by +1.0in
\advance\textheight by +1.6in
\advance\oddsidemargin by -0.5in
\advance\evensidemargin by -1.0in
\advance\topmargin by -0.5in
\input psfig.sty

\def\be{\begin{equation}}
\def\ee{\end{equation}}
\def\bq{\begin{eqnarray}}
\def\eq{\end{eqnarray}}
\def\beq{\begin{eqnarray*}}
\def\eeq{\end{eqnarray*}}

\def\ba{\begin{array}}
\def\ea{\end{array}}
\def\bth{\begin{theorem}}
\def\eth{\end{theorem}}
\def\blm{\begin{lemma}}
\def\elm{\end{lemma}}
\def\bdf{\begin{definition}}
\def\edf{\end{definition}}
\def\bpr{\begin{proposition}}
\def\epr{\end{proposition}}
\def\brm{\begin{remark}}
\def\erm{\end{remark}}
\def\bnot{\begin{notation}}
\def\enot{\end{notation}}
\def\bobs{\begin{observation}}
\def\eobs{\end{observation}}
\def\bcrl{\begin{corolary}}
\def\ecrl{\end{corolary}}

\newcommand{\ab}{\mbox{\boldmath $a$}}

\newcommand{\R}{\mathbb{R}}
\newcommand{\C}{\mathbb{C}}

\newcommand{\X}{{\bf X}}

\newcommand{\EProof}{\ \hfill\rule[-1mm]{2mm}{3.2mm}}
\newcommand{\BProof}{\noindent{\it Proof:}\ \,}

\newtheorem {theorem} {Theorem}
\newtheorem {proposition} [theorem]{Proposition}

\newtheorem {lemma}  [theorem]{Lemma}

\newtheorem {remark} [theorem]{Remark}

\newtheorem {notation} [theorem]{Notation}
\newtheorem {observation} [theorem]{Observation}

\begin{document}

\begin{center}
\textbf{ Planar cubic  polynomial differential systems with
     the maximum  number of invariant straight lines }
\end{center}

\begin{center}
{   Jaume Llibre (1) and Nicolae Vulpe (2)}
\end{center}

\begin{center}
( (1) Departament de Matem\`{a}tiques, Universitat Aut\`{o}noma de
Barcelona,  Spain\\ (2) Institute of Mathematics and Computer
Science, Acad. of Sci. of Moldova,
         Chi\c sin\u au)
\end{center}

\begin{abstract}
We classify all cubic systems possessing the maximum number of
invariant straight lines (real or complex) taking into account
their multiplicities. We prove that there are exactly $23$
topological different classes of such systems. For every class we
provide the configuration of its invariant straight lines in the
Poincar\'e disc. Moreover, every class is characterized by a set
of affine invariant conditions.
\end{abstract}

\noindent\textbf{ 1. Introduction and statement of the main results}

We consider here real polynomial differential system
\be\label{PQ}
  \frac {dx}{dt}= P(x,y),\qquad
  \frac {dy}{dt}= Q(x,y),
\ee
where $P,\ Q$ are polynomials in $x,\ y$ with real coefficients,
i.e. $P,\,Q\in \R[x,y]$. We shall say that system (\ref{PQ}) is
{\it cubic} if $\max(\deg(P),\deg(Q))=3$.

A straight line $ux+vy+w=0$ satisfies
$$
u\frac {dx}{dt}+v\frac {dy}{dt}=uP(x,y)+vQ(x,y)=(ux+vy+w)R(x,y)
$$
for some polynomial $R(x,y)$ if and only if it is {\it invariant}
under the flow of the system. If some of the coefficients $u$,
$v$, $w$ of an invariant straight line belongs to $\C\setminus\R$,
then we say that {\it the straight line is complex}; otherwise
{\it the straight line is real.} Note that, since system
(\ref{PQ}) is real, if it has a complex invariant straight line
$ux+vy+w=0,$ then it also has its conjugate complex invariant
straight line $\bar ux+\bar vy+\bar w=0.$

Let
$$
\X=P(x,y)\frac{\partial}{\partial x}
+Q(x,y)\frac{\partial}{\partial y}
$$
be the polynomial vector field corresponding to system (\ref{PQ}).

An invariant straight line $f=0$ for a cubic vector field $\X$ has
{\it geometric multiplicity} $m$ if there exists a sequence of
cubic vector fields $\X_k$ converging to $\X$, such that each
$\X_k$ has $m$ distinct invariant straight lines $f^1_k=0,\ldots,
f^m_k=0$, converging to $f=0$ as $k\to\infty$, and this does not
occur for $m+1$.

An invariant straight line $f=0$ for a cubic  vector field $\X$
has {\it algebraic multiplicity} $m$ if $m$ is the greatest
positive integer such that $f^m$ divides $P\X(Q)-Q\X(P)$. In
\cite{Llib_Per} it is proved that both notions of multiplicity
coincide. The algebraic definition of multiplicity is very useful
for its computation.

We note that this definition of multiplicity can be applied to the
infinite line $Z=0$ in the case when this line is not full of
singular points. So, including the infinite line according with
\cite{Art_Llibre} the maximum number of the invariant straight
lines for cubic systems is $9$.

  In this paper we classify all cubic systems possessing the
maximum number of invariant straight lines taking into account
their multiplicities.

Invariant straight lines for quadratic systems has been studied by
Druzhkova \cite{Druzhkova} and Popa and Sibirskii
\cite{Popa_Sib2,Popa_Sib3}, for cubic systems by Liybimova
\cite{Lyubim1,Lyubim2}, for quartic systems by Sokulski
\cite{Sokulski} and Xiang Zhang \cite{ZX}, for some more general
systems by Popa \cite{Popa2,Popa3} and Popa and Sibirskii
\cite{Popa_Sib1}

The maximum number of invariant straight lines taking into account
their multiplicities for a polynomial differential system of
degree $m$ is $3m$ when we also consider the infinite straight
line, see \cite{Art_Llibre}. This bound is always reached if we
consider the real and the complex invariant straight lines, see
\cite{Llib_Per}.

Using geometric invariants  as well as algebraic ones a
classification of all quadratic systems possessing the maximum
number of  invariant straight lines taking into account their
multiplicities have been made in \cite{Dana_Vlp2}

It is well known that for cubic system (\ref{PQ}) there exist at
most 4 different slopes  for invariant affine straight lines, for
more information about the slopes of invariant straight lines for
polynomial vector fields, see \cite{AL}.

If a cubic system (\ref{PQ}) possesses $9$ distinct invariant
straight lines we say that these lines form a {\it configuration
of  type} (3,3,1,1) if there exist two triples of parallel lines
and two additional lines every set with different slopes.  And we
shall say that these lines form a {\it configuration of  type}
(3,2,2,1) if there exist one triple and two couple of parallel
lines and one additional line every set with different slopes.
Note that in both configurations the straight line which is
omitted is the infinite one.

If a cubic system (\ref{PQ}) possesses $9$ invariant straight
lines taking into account their multiplicities we shall say that
these lines form a {\it potential configuration of type} (3,3,1,1)
(respectively, (3,2,2,1) ) if there exists a sequence of vector
fields $\X_k$ as in the definition of geometric multiplicity
having $9$ distinct line of type (3,3,1,1) (respectively,
(3,2,2,1)).

Consider generic cubic systems of the form:
\be\label{s_1} \ba{ll}
\displaystyle  \frac {dx}{dt}&=p_0+ p_1(x,y)+\,p_2(x,y)+\,p_3(x,y)
\equiv P(x,y), \\[2mm]
\displaystyle  \frac {dy}{dt}&=q_0+
q_1(x,y)+\,q_2(x,y)+\,q_3(x,y)\equiv Q(x,y), \ea \ee with real
homogeneous polynomials $p_i$ and $q_i$ $(i=0,1,2,3)$ of degree
$i$ in $x,y$. We introduce the following polynomials: \beq
   C_i&=&yp_i(x,y)-xq_i(x,y), \\
   D_j&=&\frac{\partial p_j}{\partial x}+
                        \frac{\partial q_j}{\partial y},
\eeq for $i=0,1,2,3$ and $j=1,2,3$ which in fact are $GL$--comitants, see
\cite{Sib1}.

In order to state our Main Theorem we need to construct some
$T$--comitants and $CT$--comitants (see \cite{Dana_Vlp1} for
detailed definitions) which will be responsible for the existence
of the maximum number of invariant  straight lines for  system
(\ref{s_1}). They are constructed by using the polynomials $C_i$
and $D_i$ via the differential operator $(f,g)^{(k)}$ called {\it
transvectant of index $k$} (see, for example, \cite{Gr_Yng}) which
acts on $\R[\ab,x,y]$ as follows:
$$
  (f,g)^{(k)}= \sum_{h=0}^k (-1)^h {k\choose h}
   {\partial^k f\over \partial x^{k-h}\partial y^h}\
   {\partial^k g\over \partial x^h\partial y^{k-h}}.
$$
Here $f(x,y)$ and $g(x,y)$ are polynomials in $x$ and $y$ of the
degree $r$ and $s$, respectively, and $\ab\in\R^{20}$ is the
20--tuple formed by all the coefficients of  system (\ref{s_1}).

 First we  construct the following comitants of  second
degree with respect to the coefficients of the initial system:
$$
\ba{lll}
   T_1=\left(C_0,C_1\right)^{(1)},&\quad
   T_{10}=\left(C_1,C_3\right)^{(1)},&\quad
   T_{19}=\left(C_2,D_3\right)^{(1)},\\
   T_2=\left(C_0,C_2\right)^{(1)},&\quad
   T_{11}=\left(C_1,C_3\right)^{(2)},&\quad
   T_{20}=\left(C_2,D_3\right)^{(2)},\\
   T_3=\left(C_0,D_2\right)^{(1)},&\quad
   T_{12}=\left(C_1,D_3\right)^{(1)},&\quad
   T_{21}=\left(D_2,C_3\right)^{(1)},\\
   T_4=\left(C_0,C_3\right)^{(1)},&\quad
   T_{13}=\left(C_1,D_3\right)^{(2)},&\quad
   T_{22}=\left(D_2,D_3\right)^{(1)},\\
   T_5=\left(C_0,D_3\right)^{(1)},&\quad
   T_{14}=\left(C_2,C_2\right)^{(2)},&\quad
   T_{23}=\left(C_3,C_3\right)^{(2)},\\
   T_6=\left(C_1,C_1\right)^{(2)},&\quad
   T_{15}=\left(C_2,D_2\right)^{(1)},&\quad
   T_{24}=\left(C_3,C_3\right)^{(4)},\\
   T_7=\left(C_1,C_2\right)^{(1)},&\quad
   T_{16}=\left(C_2,C_3\right)^{(1)},&\quad
   T_{25}=\left(C_3,D_3\right)^{(1)},\\
   T_8=\left(C_1,C_2\right)^{(2)},&\quad
   T_{17}=\left(C_2,C_3\right)^{(2)},&\quad
   T_{26}=\left(C_3,D_3\right)^{(2)},\\
   T_9=\left(C_1,D_2\right)^{(1)},&\quad
   T_{18}=\left(C_2,C_3\right)^{(3)},&\quad
   T_{27}=\left(D_3,D_3\right)^{(2)}.\\
\ea
$$
Then we need the following polynomials:
\beq
{\cal D}_1(\ab) &=& 6T_{24}^3-\left[(C_3,T_{23})^{(4)}\right]^2,\\
{\cal D}_2(\ab,x,y) &=& -T_{23}, \\
{\cal D}_3(\ab,x,y) &=& (T_{23},\,T_{23})^{(2)}
-6C_3(C_3,\,T_{23})^{(4)},
\eeq
\beq
{\cal D}_4(\ab)& =& (C_3,\, {\cal D}_2)^{(4)},\\
{\cal V}_1(\ab,x,y)&=& T_{23}+2D_3^2,\\
{\cal V}_2(\ab,x,y)&=&T_{26},\\
{\cal V}_3(\ab,x,y)&=&6T_{25}-3T_{23}-2D_3^2,\\
{\cal V}_4(\ab,x,y)&=&
C_3\left[\left(C_3,T_{23}\right)^{(4)}+36\left(D_3,T_{26}\right)^{(2)}\right],\\
{\cal L}_1(\ab,x,y)&=& 9C_2\left(T_{24}+ 24T_{27}\right)
 -12D_3\left(T_{20}+8T_{22}\right)-12\left(T_{16},D_3\right)^{(2)}\\
&& -3\left(T_{23},C_2\right)^{(2)}
 -16\left(T_{19},C_3\right)^{(2)}
 +12\left(5T_{20}+24T_{22},C_3\right)^{(1)},\\
{\cal L}_2(\ab,x,y)&=&32\left(13T_{19}+33T_{21},D_2\right)^{(1)}+
          84\left(9T_{11}-2T_{14},D_3\right)^{(1)}\\
       &&+8D_2\left(12T_{22}+35T_{18}-73T_{20}\right)
           -448\left(T_{18},C_2\right)^{(1)}\\
       && -56\left(T_{17},C_2\right)^{(2)}
           -63\left(T_{23},C_1\right)^{(2)}
            +756D_3T_{13} -1944D_1T_{26}\\
        &&+112\left(T_{17},D_2\right)^{(1)}
           -378\left(T_{26},C_1\right)^{(1)}
          +9C_1\left(48T_{27}-35T_{24}\right),\\
{\cal L}_3(\ab,x,y)&=&\left(T_{23}, D_3\right)^{(2)}
         \left[\left(D_2,T_{22}\right)^{(1)}-D_1T_{27}\right],\\
{\cal L}_4(\ab,x,y)&=& T_{25},\\
{\cal N}_1(\ab,x,y)&=&4C_2(27D_1D_3-8D_2^2)+
                    2C_2(20T_{15}-4T_{14}+39T_{12})\\
        && +18C_1(3T_{21}-D_2D_3)+54D_3(3T_4-T_7)-288C_3T_9\\
        && +54\left(T_{7},C_3\right)^{(1)}
         -567\left(T_{4},C_3\right)^{(1)}+ 135C_0D_3^2,\\
{\cal N}_2(\ab,x,y)&=&2C_2D_3-3C_3D_2,\\
{\cal N}_3(\ab,x,y)&=&C_2D_3+3T_{16},\\
{\cal N}_4(\ab,x,y)&=&D_2D_3+9T_{21}-2T_{17},\\
{\cal N}_5(\ab,x,y)&=&T_{17}+2T_{19},\\
{\cal N}_6(\ab,x,y)&=&6C_3(T_{12}+6T_{11})-9C_1(T_{23}+T_{25})
                  -8 \left(T_{16},C_2\right)^{(1)}-C_3D_2^2,\\
{\cal N}_7(\ab,x,y)&=&6C_3(12T_{11}-T_{12}-6D_1D_3)-21C_1T_{23}
                  -24\left(T_{16},C_2\right)^{(1)}\\
                 && +3C_1T_{25}+4D_2(T_{16}+2D_2C_3-C_2D_3),\\
{\cal N}_8(\ab,x,y)&=&D_2^2-4D_1D_3,\\
{\cal N}_9(\ab,x,y)&=&C_2^2-3C_1C_3,\\
{\cal N}_{10}(\ab,x,y)&=&2C_2D_1+3T_4.
\eeq

\noindent\textbf{ Main Theorem} {\it Any cubic system having invariant
straight lines with total multiplicity $9$ via affine
transformation and time rescaling can be written as one of the
following $23$ systems. In the figure associated to each system is
presented the configuration of its invariant straight lines in the
Poincar\'e disc. Real invariant straight lines are represented by
continuous lines. Complex  invariant straight lines are
represented by dashed lines. If an invariant straight line has
multiplicity $k>1$, then the number $k$ appears near the
corresponding straight line and this line is more thick. Moreover,
every system has associated a set of affine invariant conditions
which characterize it.}
$$
\ba{llcll}
 (1)\ \ba{l} \dot x=x(x^2-1),\\
                \dot y=y(y^2-1)  \\
         \ea & \Leftrightarrow &
           \left[\ba{l}{\cal D}_1>0,\ {\cal D}_2>0,\ {\cal D}_3>0,\ {\cal L}_3<0,\\
                   {\cal V}_{1}={\cal V}_{2}= {\cal L}_{1} =
                   {\cal L}_{2}= {\cal N}_{1}=0\\
                 \ea\right] &  \Leftrightarrow &  Fig.\, 1;\\[7mm]
 (2)\ \ba{l} \dot x=x(x^2+1),\\
                \dot y=y(y^2+1)  \\
         \ea & \Leftrightarrow &
           \left[\ba{l}{\cal D}_1>0,\ {\cal D}_2>0,\ {\cal D}_3>0,\ {\cal L}_3>0,\\
                   {\cal V}_{1}={\cal V}_{2}= {\cal L}_{1} =
                   {\cal L}_{2}= {\cal N}_{1}=0\\
                 \ea\right] &  \Leftrightarrow &  Fig.\,2;\\[7mm]
  (3)\ \ba{l} \dot x=x^3,\\
                \dot y=y^3  \\
         \ea & \Leftrightarrow &
           \left[\ba{l}{\cal D}_1>0,\ {\cal D}_2>0,\ {\cal D}_3>0,\ {\cal L}_3=0,\\
                   {\cal V}_{1}={\cal V}_{2}= {\cal L}_{1} =
                   {\cal L}_{2}= {\cal N}_{1}=0\\
                 \ea\right] &  \Leftrightarrow &\!\! \ba{l} Fig.\,3;\\
                                                    \ \mbox{ \scriptsize
                                                     (\ref{sys_conf:1-3})}
                                                   \ea\\[7mm]
  (4)\ \ba{l} \dot x\!=\!2x(x^2\!-\!1),\\
                \dot y\!=\!(3x\!-\!y)(y^2\!-\!1)\!\!  \\
         \ea  & \Leftrightarrow &
           \left[\ba{l}{\cal D}_1>0,\ {\cal D}_2>0,\ {\cal D}_3>0,\ {\cal L}_3>0,\\
                   {\cal V}_{3}={\cal V}_{4}= {\cal L}_{1} =
                   {\cal L}_{2}= {\cal N}_{1}=0\\
                 \ea\right] &  \Leftrightarrow &  Fig.\,4;\\[7mm]
 (5)\ \ba{l} \dot x=2x(x^2\!+\!1),\\
                \dot y=(3x\!-\!y)(y^2\!+\!1)\!\!  \\
         \ea  & \Leftrightarrow &
           \left[\ba{l}{\cal D}_1>0,\ {\cal D}_2>0,\ {\cal D}_3>0,\ {\cal L}_3<0,\\
                   {\cal V}_{3}={\cal V}_{4}= {\cal L}_{1} =
                   {\cal L}_{2}= {\cal N}_{1}=0\\
                 \ea\right] &  \Leftrightarrow &  Fig.\,5;\\[7mm]
  (6)\ \ba{l} \dot x=2x^3,\\
                \dot y=y^2(3x-y)  \\
         \ea  & \Leftrightarrow &
           \left[\ba{l}{\cal D}_1>0,\ {\cal D}_2>0,\ {\cal D}_3>0,\ {\cal L}_3=0,\\
                   {\cal V}_{3}={\cal V}_{4}= {\cal L}_{1} =
                   {\cal L}_{2}= {\cal N}_{1}=0\\
                 \ea\right] &  \Leftrightarrow &\!\!  \ba{l} Fig.\,6;\\
                                                    \ \hbox{\footnotesize
                                                      (\ref{sys_conf:16-18})}
                                                   \ea\\[7mm]
 (7)\ \ba{l} \dot x=x(1+x^2),\\
                \dot y=y(1-y^2)  \\
         \ea  & \Leftrightarrow &
           \left[\ba{l}{\cal D}_1<0,\ \ {\cal L}_3\not=0,\ \ {\cal L}_4<0,\\
                   {\cal V}_{1}={\cal V}_{2}= {\cal L}_{1} =
                   {\cal L}_{2}= {\cal N}_{1}=0\\
                 \ea\right] &  \Leftrightarrow &  Fig.\,7;\\[7mm]
  (8)\ \ba{l} \dot x=x^3,\\
                \dot y=-y^3, \\
         \ea  & \Leftrightarrow &
           \left[\ba{l}{\cal D}_1<0,\ \ {\cal L}_3=0,\ \ {\cal L}_4<0,\\
                   {\cal V}_{1}={\cal V}_{2}= {\cal L}_{1} =
                   {\cal L}_{2}= {\cal N}_{1}=0\\
                 \ea\right] &  \Leftrightarrow &\!\!  \ba{l} Fig.\,8;\\
                                                    \ \hbox{\footnotesize
                                                      (\ref{sys_conf:4-5})}
                                                   \ea\\[7mm]
  (9)\ \ba{l} \dot x=x(1\!+\!x^2\!-\!3y^2),\!\!\\
                \dot y=y(1\!+\!3x^2\!-\!y^2)\!\!  \\
         \ea  & \Leftrightarrow &
           \left[\ba{l}{\cal D}_1<0,\ \  {\cal L}_3\not=0,\ \ {\cal L}_4>0,\\
                   {\cal V}_{1}={\cal V}_{2}= {\cal L}_{1} =
                   {\cal L}_{2}= {\cal N}_{1}=0\\
                 \ea\right] &  \Leftrightarrow &  Fig.\,9;\\[7mm]
 (10)\ \ba{l} \dot x=x(x^2-3y^2),\\
                \dot y=y(3x^2-y^2)  \\
         \ea  & \Leftrightarrow &
           \left[\ba{l}{\cal D}_1<0,\ \ {\cal L}_3=0,\ \ {\cal L}_4>0,\\
                   {\cal V}_{1}={\cal V}_{2}= {\cal L}_{1} =
                   {\cal L}_{2}= {\cal N}_{1}=0\\
                 \ea\right] &  \Leftrightarrow &\!\!  \ba{l} Fig.\,10;\\
                                                    \ \hbox{\footnotesize
                                                      (\ref{sys_conf:6-7})}
                                                   \ea\\[7mm]
  (11)\  \ba{l} \dot x=2x(x^2-1),\\
                    \dot y=y(3x^2\!+y^2\!+\!1)\!\!  \\
         \ea  & \Leftrightarrow &
           \left[\ba{l}{\cal D}_1<0,\ \ {\cal L}_3<0,\\
                   {\cal V}_{3}={\cal V}_{4}= {\cal L}_{1} =
                   {\cal L}_{2}= {\cal N}_{1}=0\\
                 \ea\right] &  \Leftrightarrow &  Fig.\,11;\\[7mm]
  (12)\ \ba{l} \dot x=2x(x^2+1),\\
                    \dot y=y(3x^2\!+\!y^2\!-\!1)\!\!  \\
         \ea  &\!\! \Leftrightarrow \!\!&
           \left[\ba{l}{\cal D}_1<0,\ \ {\cal L}_3>0,\\
                   {\cal V}_{3}={\cal V}_{4}= {\cal L}_{1} =
                   {\cal L}_{2}= {\cal N}_{1}=0\\
                 \ea\right] & \!\! \Leftrightarrow \!\! &\  Fig.\,12;\\[6mm]
  (13)\ \ba{l} \dot x=2x^3,\\
                    \dot y=y(3x^2+y^2)  \\
         \ea  &\!\! \Leftrightarrow\!\! &
           \left[\ba{l}{\cal D}_1<0,\ \ {\cal L}_3=0,\\
                   {\cal V}_{3}={\cal V}_{4}= {\cal L}_{1} =
                   {\cal L}_{2}= {\cal N}_{1}=0\\
                 \ea\right] &\!\!  \Leftrightarrow\!\! &  \ba{l} Fig.\,13\\
                                                    \ \hbox{\footnotesize
                                                     (\ref{sys_conf:19-21})};
                                                   \ea\\[6mm]
  (14)\ \ba{l} \dot x=x(x^2-1),\\
                \dot y=2y  \\
         \ea  &\!\! \Leftrightarrow\!\! &
           \left[\ba{l}{\cal D}_1\!=\! {\cal D}_3\!=\!{\cal D}_4=0,\ {\cal D}_2\!\not=\!0,\ {\cal L}_4<0, \\
                   {\cal V}_{1}\!=\!{\cal N}_{1}\!=\! {\cal N}_{2}\!=\!
                   {\cal N}_{3}\!=\! {\cal N}_{7}\!=\!0,\, {\cal N}_8<0 \\
                 \ea\right] &\!\!  \Leftrightarrow \!\! &\!\!  \ba{l}\  Fig.\,14\\
                                                      \hbox{\footnotesize
                                                       (\ref{s:pert1}),(\ref{s:pert2})};\\
                                                   \ea\\[6mm]
  (15)\ \ba{l} \dot x=x(x^2+1),\\
                \dot y=-2y  \\
         \ea  &\!\! \Leftrightarrow\!\! &
           \left[\ba{l}{\cal D}_1\!=\! {\cal D}_3\!=\!{\cal D}_4=0,\ {\cal D}_2\not=0,\ {\cal L}_4<0, \\
                   {\cal V}_{1}\!=\!{\cal N}_{1}\!=\! {\cal N}_{2}\!=\!
                   {\cal N}_{3}\!=\! {\cal N}_{7}\!=\!0,\, {\cal N}_8>0 \\
                 \ea\right] &\!\!  \Leftrightarrow\!\! &\!\!  \ba{l}\  Fig.\,15\\
                                                      \hbox{\footnotesize
                                                     (\ref{s:pert3}),(\ref{s:pert4})};\\
                                                   \ea\\[6mm]
\ea
$$
$$
\ba{llcll}
  (16)\ \ba{l} \dot x=x(x^2-1),\\
                \dot y=-y  \\
         \ea  &\!\! \Leftrightarrow\!\! &
           \left[\ba{l}{\cal D}_1\!=\! {\cal D}_3\!=\!{\cal D}_4=0,\ {\cal D}_2\not=0,\ {\cal L}_4<0, \\
                   {\cal V}_{1}\!=\!{\cal N}_{1}\!=\! {\cal N}_{2}\!=\!
                   {\cal N}_{3}\!=\! {\cal N}_{6}=0,\, {\cal N}_8>0 \\
                 \ea\right] &\!\!  \Leftrightarrow\!\! &\!\!  \ba{l}\ Fig.\,16\\
                                                      \hbox{\footnotesize
                                                     (\ref{s:pert5}),(\ref{s:pert6})};\\
                                                   \ea\\[6mm]
  (17)\ \ba{l} \dot x=x(x^2+1),\\
                \dot y=y \\
         \ea &\!\! \Leftrightarrow\!\! &
           \left[\ba{l}{\cal D}_1\!=\! {\cal D}_3\!=\!{\cal D}_4=0,\ {\cal D}_2\not=0,\ {\cal L}_4<0, \\
                   {\cal V}_{1}\!=\!{\cal N}_{1}\!=\! {\cal N}_{2}\!=\!
                   {\cal N}_{3}\!=\! {\cal N}_{6}=0,\, {\cal N}_8<0 \\
                 \ea\right] & \!\! \Leftrightarrow\!\! &\!\!  \ba{l}\ Fig.\,17\\
                                                      \hbox{\footnotesize
                                                     (\ref{s:pert7}),(\ref{s:pert8})};\\
                                                   \ea\\[6mm]
  (18)\ \ba{l} \dot x=x^3,\\
                \dot y=1 \\
         \ea &\!\! \Leftrightarrow\!\! &
           \left[\ba{l}{\cal D}_1\!=\! {\cal D}_3\!=\!{\cal D}_4=0,\ {\cal D}_2\not=0,\ {\cal L}_4<0, \\
                   {\cal V}_{1}\!=\!{\cal N}_{1}\!=\! {\cal N}_{2}\!=\!
                   {\cal N}_{3}\!=\! {\cal N}_{6}=0,\, {\cal N}_8=0 \\
                 \ea\right] &\!\!  \Leftrightarrow\!\! &\!\!  \ba{l}\ Fig.\,18\\
                                                      \hbox{\footnotesize
                                                     (\ref{s:pert9}),(\ref{s:pert10})};\\
                                                   \ea\\[6mm]
  (19)\ \ba{l} \dot x=x(x^2-1),\\
                \dot y=y(3x^2-1) \\
         \ea  &\!\! \Leftrightarrow\!\! &
           \left[\ba{l}{\cal D}_1\!=\! {\cal D}_3\!=\!{\cal D}_4=0,\ {\cal D}_2\not=0,\ {\cal L}_4>0, \\
                   {\cal V}_{1}\!=\!{\cal N}_{1}\!=\! {\cal N}_{2}\!=\!
                   {\cal N}_{3}\!=\! {\cal N}_{6}=0,\, {\cal N}_8>0 \\
                 \ea\right] &\!\!  \Leftrightarrow\!\! &  \ba{l} Fig.\,19\\
                                                    \ \hbox{\footnotesize
                                                      (\ref{s:pert11})};\\
                                                   \ea\\[6mm]
  (20)\ \ba{l} \dot x=x(x^2+1),\\
                \dot y=y(3x^2+1) \\
         \ea  &\!\! \Leftrightarrow\!\! &
           \left[\ba{l}{\cal D}_1\!=\! {\cal D}_3\!=\!{\cal D}_4=0,\ {\cal D}_2\not=0,\ {\cal L}_4>0, \\
                   {\cal V}_{1}\!=\!{\cal N}_{1}\!=\! {\cal N}_{2}\!=\!
                   {\cal N}_{3}\!=\! {\cal N}_{6}\!=\!0,\, {\cal N}_8<0 \\
                 \ea\right] &\!\!  \Leftrightarrow \!\!&  \ba{l} Fig.\,20\\
                                                    \ \hbox{\footnotesize
                                                     (\ref{s:pert12})};\\
                                                   \ea\\[6mm]
  (21)\ \ba{l} \dot x=2x(x^2-1),\\
                \dot y=y(3x^2+1) \\
         \ea  &\!\! \Leftrightarrow\!\! &
           \left[\ba{l}{\cal D}_1\!=\! {\cal D}_3\!=\!{\cal D}_4=0,\ {\cal D}_2\not=0,\ {\cal L}_4>0, \\
                   {\cal V}_{3}\!=\!{\cal N}_{1}\!=\! {\cal N}_{4}\!=\!
                   {\cal N}_{5}\!=\! {\cal N}_{7}\!=\!0,\, {\cal N}_8>0 \\
                 \ea\right] & \!\! \Leftrightarrow \!\!&   \ba{l} Fig.\,21\\
                                                    \ \hbox{\footnotesize
                                                     (\ref{s:pert13})};\\
                                                   \ea\\[6mm]
  (22)\ \ba{l} \dot x=2x(x^2+1),\\
                \dot y=y(3x^2-1) \\
         \ea  &\!\! \Leftrightarrow\!\! &
           \left[\ba{l}{\cal D}_1\!=\! {\cal D}_3\!=\!{\cal D}_4=0,\ {\cal D}_2\not=0,\ {\cal L}_4>0, \\
                   {\cal V}_{3}\!=\!{\cal N}_{1}\!=\! {\cal N}_{4}\!=\!
                   {\cal N}_{5}\!=\! {\cal N}_{7}\!=\!0,\, {\cal N}_8<0 \\
                 \ea\right] &\!\!  \Leftrightarrow \!\!&   \ba{l} Fig.\,22\\
                                                    \ \hbox{\footnotesize
                                                     (\ref{s:pert14})};\\
                                                   \ea\\[6mm]
  (23)\ \ba{l} \dot x=x,\\
                \dot y=y-x^3 \\
         \ea  &\!\! \Leftrightarrow\!\! &
           \left[\ba{l}{\cal D}_1= {\cal D}_2={\cal D}_3={\cal V}_1=0, \\
                   {\cal N}_{2}= {\cal N}_{3}=
                   {\cal N}_{9}= {\cal N}_{10}=0\\
                 \ea\right] &\!\!  \Leftrightarrow \!\!&\!\!  \ba{l}\ Fig.\,23\\
                                                      \hbox{\footnotesize
                                                     (\ref{s:pert15}),(\ref{s:pert16})}.\\
                                                   \ea\\[6mm]
\ea
$$

Here, a condition of the type $S(\ab,x,y)>0$ (respectively,
$S(\ab,x,y)<0$) means that  the respective homogeneous polynomial
of even degree in $x$ and $y$ is positive (respectively, negative)
defined.  And equality $S(\ab,x,y)=0$  must be understood in
$\R[x,y]$.

Note that, only in the case that some invariant straight lines
have multiplicity $>1$, in the last column of the statement of the
Main Theorem appear some numbers under the figures indicating the
corresponding perturbed systems which show the potential
configurations of the considered cubic system.

\begin{figure}
\centerline {\psfig{figure=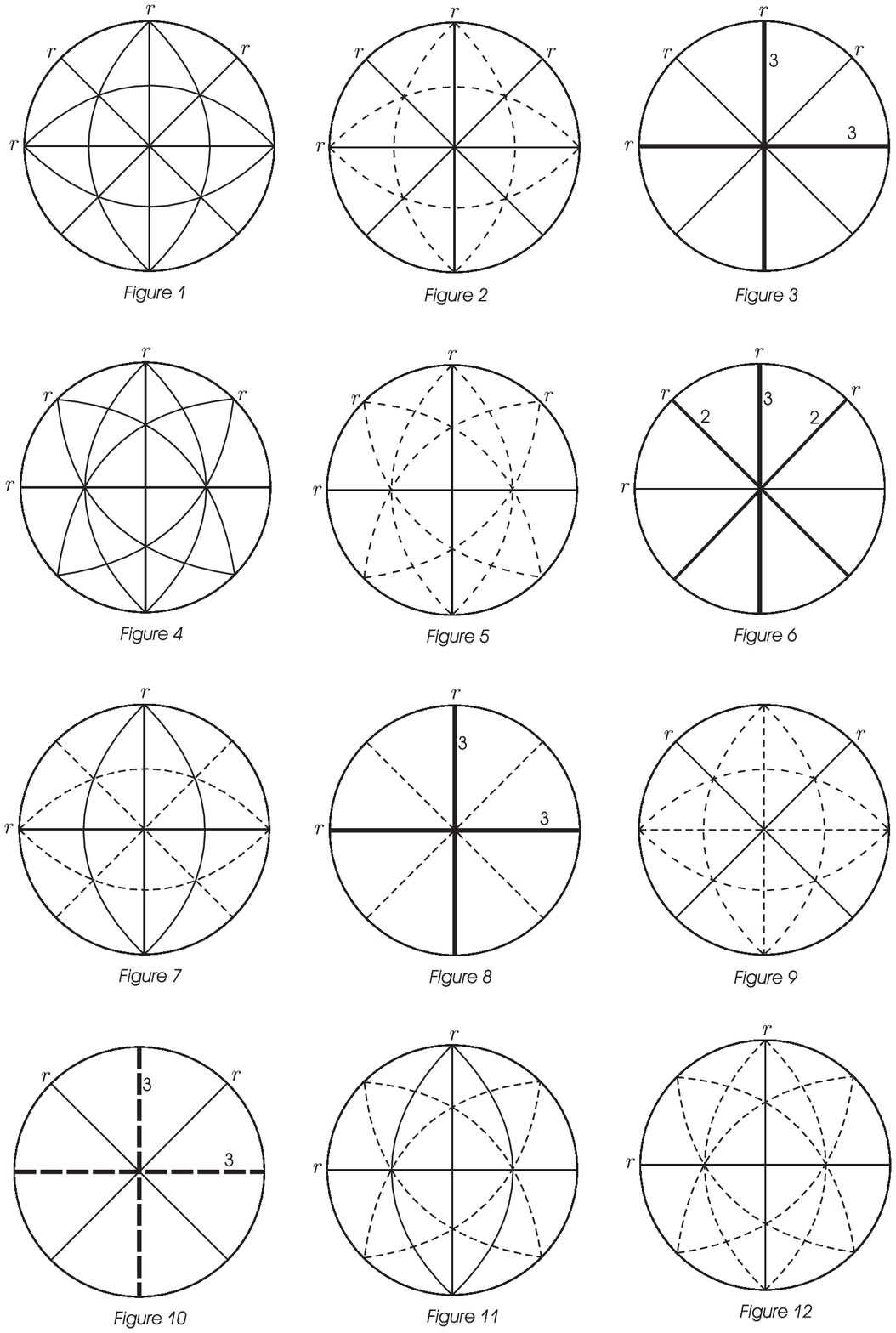,height=20.3cm}}
\end{figure}

\begin{figure}
\centerline {\psfig{figure=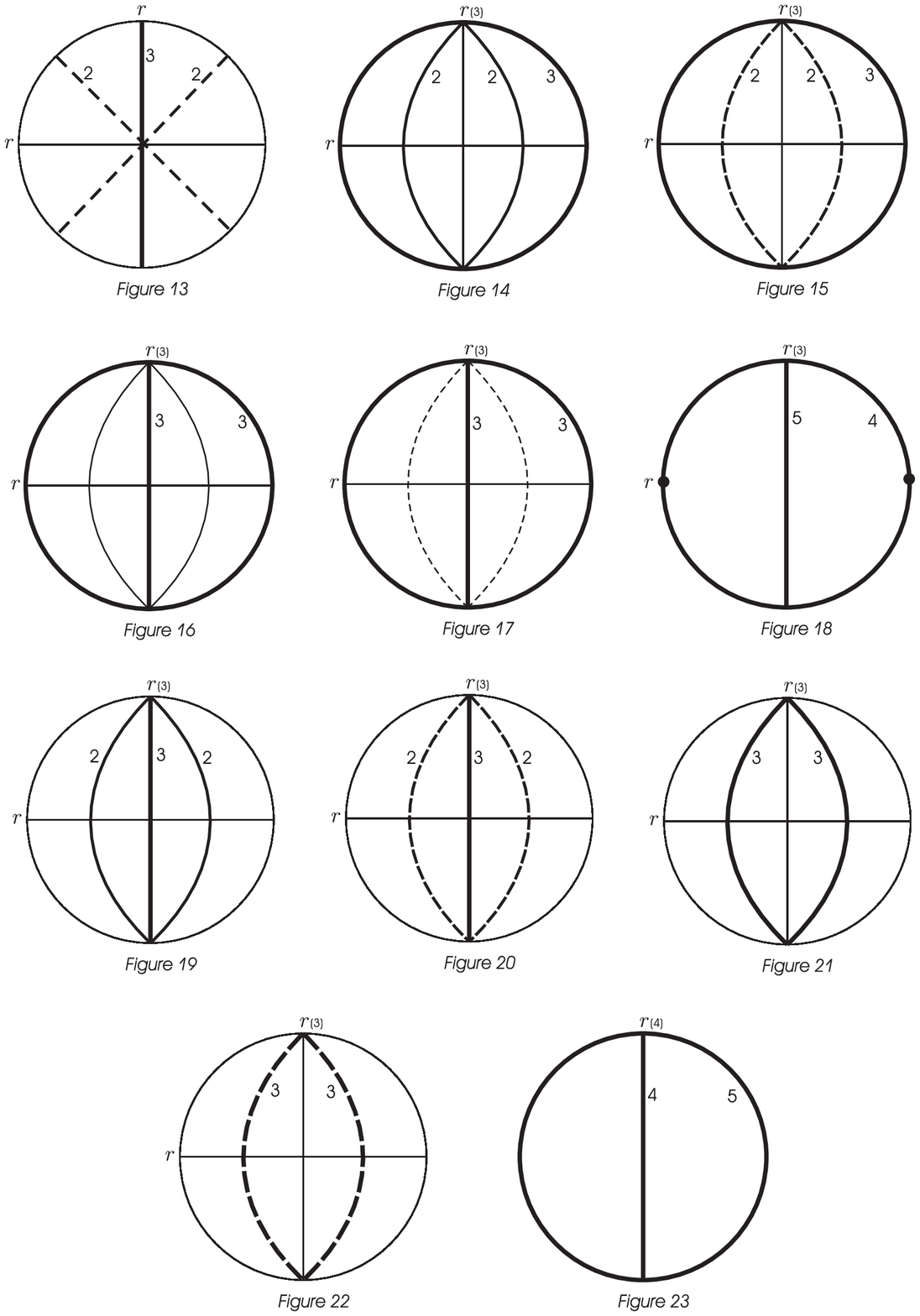,height=20.5cm}}
\end{figure}

\bigskip

\noindent\textbf{ 2. Necessary conditions for the existence of
parallel invariant straight lines}

We define the  auxiliary polynomial ${\cal U}_{\,1}(\ab)=
T_{24}-4T_{27}.$

\blm\label{lm:1}
For cubic systems~(\ref{s_1}) the conditions ${\cal V}_1= {\cal
V}_2= {\cal U}_{\,1}=0$ are necessary  for  the existence of two
triples of parallel invariant straight lines with different slope.
\elm

\BProof Let $ L_i(x,y)= \alpha x+\beta y+\gamma_i=0,\ (i=1,2,3)$
be three parallel invariant straight  lines for a cubic  system
(\ref{s_1}). Then,  we have
$$
      \alpha P(x,y)+\beta Q(x,y)=\xi(\alpha x+\beta y+\gamma_1)
      (\alpha x+\beta y+\gamma_2)(\alpha x+\beta y+\gamma_3),
$$
where the constant $\xi$ can be considered 1 (rescaling  the time,
if necessary). Therefore, from the cubic terms we obtain $\alpha
p_3(x,y)+\beta q_3(x,y) =  (\alpha x+\beta y)^3$. If we denote
$$
p_3(x,y)=px^3+3qx^2y+3rxy^2+sy^3,\quad
q_3(x,y)=tx^3+3ux^2y+3vxy^2+wy^3.
$$
then, for the existence of 3 parallel invariant straight lines it
is necessary the solvability of the following systems of cubic
equations with respect to the parameters $\alpha$ and $\beta$:
\be\label{eq_1} \ba{llllll}
    A_1\equiv &   \alpha p + \beta t - \alpha^3&=0,\quad &
   A_2\equiv  &   \alpha q + \beta u - \alpha^2\beta&=0,\\
   A_3\equiv  &   \alpha r + \beta v - \alpha\beta^2&=0,\quad &
   A_4\equiv  &   \alpha s + \beta w - \beta^3&=0.\\
\ea \ee Without loss of generality we can consider
$\alpha\beta\ne0$, otherwise a rotation of the phase plane can be
done. We have: \bq
B_1\equiv    \alpha A_2- \beta A_1 & = &  q\alpha^2+(u-p)\alpha\beta-t\beta^2=0,\nonumber \\
B_2\equiv    \alpha A_3- \beta A_2 & = &  r\alpha^2+(v-q)\alpha\beta-u\beta^2=0,\label{eq_2} \\
B_3\equiv    \alpha A_4- \beta A_3 & = &
s\alpha^2+(w-r)\alpha\beta-v\beta^2=0.\nonumber \eq Clearly, for
the existence of two directions $(\alpha_1,\beta_1)$ and
$(\alpha_2,\beta_2)$ such that in each of them there are 3
parallel invariant straight lines of a system (\ref{s_1}) it is
necessary that the ${\rm rank}( U)=1$, where
$$
    U=\left(\ba{ccc}  q &\ u\!-\!p\ & -t \\
                     r &\ v\!-\!q\ & -u \\
                     s &  w\!-\!r\ & -v \\
           \ea\right).
$$
We denote by $U^{ij}_{kl}$ the $2\times2$--minor of the matrix $U$
formed only by the columns $i$ and $j$ and by the rows $k$ and $l$
of $U$. We obtain
$$
\ba{lll}
U^{12}_{12}=\left|\ba{cc} q &\, u\!-\!p \\
                          r &\, v\!-\!q \\
                  \ea\right|,\quad &
U^{13}_{12}=\left|\ba{cc} q &\, -t \\
                          r &\, -u \\
                  \ea\right|,\quad&
U^{23}_{12}=\left|\ba{cc}  u\!-\!p &\, -t \\
                           v\!-\!q &\, -u\\
                  \ea\right|,\\
\rule{0pt}{6ex}
U^{12}_{13}=\left|\ba{cc} q &\, u\!-\!p \\
                          s &\, w\!-\!r \\
                  \ea\right|,\quad &
U^{13}_{13}=\left|\ba{cc} q &\, -t \\
                          s &\, -v \\
                  \ea\right|, \quad &
U^{23}_{13}=\left|\ba{cc}  u\!-\!p &\, -t \\
                           w\!-\!r &\, -v\\
                  \ea\right|,\\
\rule{0pt}{6ex}
U^{12}_{23}=\left|\ba{cc} r &\, v\!-\!q \\
                          s &\, w\!-\!r \\
                  \ea\right|,\quad &
U^{13}_{23}=\left|\ba{cc} r &\, -u \\
                          s &\, -v \\
                  \ea\right|,\quad &
U^{23}_{23}=\left|\ba{cc}  v\!-\!q &\, -u \\
                           w\!-\!r &\, -v\\
                  \ea\right|. \\
\ea
$$
Hence, the ${\rm rank}(U)=1$ if and only if $U^{ij}_{kl}=0$ for
all $1\le i<j\le 3$ and $1\le k<l\le 3$.

On the other hand, it is easy to calculate the values of the
$T$--comitants
 \beq {\cal V}_1&=&16
\left[U^{23}_{12}x^4+\left(U^{23}_{13}-2U^{13}_{12}\right)x^3y+
           \left(U^{12}_{12}-2U^{13}_{13}+U^{23}_{23}\right)x^2y^2\right.\\
         &&\quad+ \left.\left(U^{12}_{13}-2U^{13}_{23}\right)xy^3+U^{12}_{23}y^4\right],\\
{\cal V}_2&=&8\left[-\left(2U^{13}_{12}+U^{23}_{13}\right)x^2+
             2\left(U^{12}_{12}-U^{23}_{23}\right)xy+
              \left(U^{12}_{13}+2U^{13}_{23}\right)y^2\right],\\
{\cal U}_{\,1}&=&
2^7\left(U^{12}_{12}+U^{13}_{13}+U^{23}_{23}\right).
\eeq
Thus, it is obvious  that $U^{ij}_{kl}=0$\ ($1\le i<j\le 3$, $1\le
k<l\le 3$) if and only if ${\cal V}_1= {\cal V}_2= {\cal
U}_{\,1}=0.$ This completes the proof of the lemma.\EProof

We assume that  ${\cal V}_1^2+{\cal V}_2^2+{\cal U}_{\,1}^2\ne0$.
Then, by Lemma \ref{lm:1}, there cannot exist two triples of
parallel invariant straight  lines for  system (\ref{s_1}). Now,
we shall examine the case when a system (\ref{s_1})  possesses
only one triple of  parallel invariant straight lines. This means
that  system  (\ref{eq_2}) can have at most one solution
$(\alpha_0,\beta_0)$. By using (\ref{eq_2}) and considering
(\ref{eq_1}), we construct the following linear system with
respect to the parameters $\alpha$ and $\beta$: \be\label{eq_3}
\ba{lll}
\alpha B_1 =& q\alpha^3+(u-p)\alpha^2\beta-t\alpha\beta^2&= B_{11}\alpha+B_{12}\beta=0,\\
\beta B_1 =&  q\alpha^2\beta+(u-p)\alpha\beta^2-t\beta^3 &= B_{21}\alpha+B_{22}\beta=0,\\
\alpha B_2 =& r\alpha^3+(v-q)\alpha^2\beta-u\alpha\beta^2&= B_{31}\alpha+B_{32}\beta=0,\\
\beta B_2 =&  r\alpha^2\beta+(v-q)\alpha\beta^2-u\beta^3 &= B_{41}\alpha+B_{42}\beta=0,\\
\alpha B_3 =& s\alpha^3+(w-r)\alpha^2\beta-v\alpha\beta^2&= B_{51}\alpha+B_{52}\beta=0,\\
\beta B_3 =&  s\alpha^2\beta+(w-r)\alpha\beta^2-v\beta^3 &= B_{61}\alpha+B_{62}\beta=0,\\
\ea
\ee
where
$$
\ba{llll}
B_{11}=-U^{13}_{12}, & B_{12}=-U^{23}_{12},
    & B_{21}=- U^{12}_{12}-U^{13}_{13},
    & B_{22}=- U^{13}_{12}-U^{23}_{13},\\
B_{31}=U^{12}_{12}, & B_{32}=U^{13}_{12},
    & B_{41}=-U^{13}_{23}, & B_{42}=-U^{23}_{23},\\
 B_{51}= U^{12}_{13}+U^{13}_{23},
    & B_{52}= U^{13}_{13}+U^{23}_{23},
    & B_{61}= U^{12}_{23}, & B_{62}= U^{13}_{23}.\\
\ea
$$
We denote by ${\cal B}$ the $6\times2$--matrix of the linear
system (\ref{eq_3}) and by $M_{ij}$ its corresponding $2\times2$
minors:
$$
     {\cal B}=\left( B_{ij}\right)_{\{i=1,\ldots,6;\ j=1,2\}},\quad
     M_{ij}=\left|\ba{ll} B_{i1} & B_{i2}\\
                          B_{j1} & B_{j2}\\
                      \ea\right|,\ \ 1\le i<j\le6.
$$
It is clear that the linear system  (\ref{eq_3}) has a
non--trivial solution if and only if the ${\rm rank}({\cal B})=1$,
i.e. $M_{ij}=0$ for $1\le i<j\le 6$.

On the other hand calculating the polynomial ${\cal V}_3$ as well
as the auxiliary polynomials \beq
{\cal U}_{\,2}(\ab,x,y)&=&6\left(T_{23}-3T_{25},T_{26}\right)^{(1)}-3T_{23}(T_{24}+8T_{27})-24T_{26}^2\\
          &&+2C_3\left(C_3,T_{23}\right)^{(4)} +24
             D_3\left(D_3,T_{26}\right)^{(1)}+24
             D_3^2T_{27},\\
{\cal U}_{\,3}(\ab,x,y)&=&
D_3\left[\left(C_3,T_{23}\right)^{(4)}+36\left(D_3,T_{26}\right)^{(2)}\right],
\eeq
for system (\ref{s_1}) we have
$$
\ba{ll} {\cal V}_4 =&
2^{10}\cdot3^5\left[ \left(M_{12}+M_{13}\right)x^4
   +3\left(M_{16}+M_{34}\right)x^2y^2+\left(M_{46}+M_{56}\right)y^4\right.\\
   &+\left.\left(3M_{14} + M_{15}- M_{23}\right)x^3y
   + \left(M_{26} + 3M_{36}- M_{45}\right)xy^3\right],\\
\rule{0pt}{3ex}
 {\cal U}_{\,2} =& 2^{12}\cdot3^5\cdot5\left[ M_{13}x^4
   +\left(2M_{16}+M_{25}-M_{34}\right)x^2y^2+M_{46}y^4\right.\\
   &+\left.\left(M_{14} + M_{15}+ M_{23}\right)x^3y
   + \left(M_{26} + M_{36}+ M_{45}\right)xy^3\right],\\
\rule{0pt}{3ex} {\cal U}_{\,3} =& 2^{10}\cdot3^6\left[
\left(M_{15}-M_{14}-M_{23}\right)x^2+
   \left(M_{16} -2M_{24}- M_{34}\right)xy\right.\\
   &+\left. \left(M_{36} - M_{26}+ M_{45}\right)y^2\right].\\
\ea
$$
It is not difficult to observe that conditions $ {\cal V}_4={\cal
U}_{\,2}={\cal U}_{\,3}=0$ are equivalent to $M_{ij}=0$ for $1\le
i<j\le 6$. Moreover, taking into account the expressions of the
polynomials ${\cal V}_4$ and ${\cal U}_{\,3}$, we can conclude for
$C_3\ne0$, and that condition ${\cal V}_4=0$ implies ${\cal
U}_{\,3}=0$. Hence, the following lemma was proved:

\blm\label{lm:2} For cubic system (\ref{s_1}) the necessary
conditions for  the existence of one triple of parallel invariant
straight lines are ${\cal V}_4= {\cal U}_{\,2}= 0.$
\elm

The next step is to find  some necessary conditions in order that
system (\ref{s_1}) possesses three couples of parallel invariant straight
lines.

Let $L_i(x,y)= \alpha x+\beta y+\gamma_i=0, (i=1,2)$ be two
parallel invariant straight lines for  a cubic system (\ref{s_1}).
Then, we have
$$
      \alpha P(x,y)+\beta Q(x,y)=(\alpha x+\beta y+\gamma_1)
      (\alpha x+\beta y+\gamma_2)(\mu x+\eta y+\gamma_3).
$$
 Therefore, from  the cubic terms we get
$$
  \alpha p_3(x,y)+\beta q_3(x,y) =
  (\alpha x+\beta y)^2(\mu x+\eta y).
$$
Thus, for the existence of 2 parallel invariant straight lines it
is necessary the solvability of the following system of cubic
equations with respect to parameters $\alpha$, $\beta$, $\mu$ and
$\eta$: \be\label{eq_4} \ba{llll}
   E_1\equiv &   \alpha p  +\beta t - \alpha^2\mu=0,\ &
   E_2\equiv &   3\alpha q + 3\beta u - \alpha^2\eta-2\alpha\beta\mu =0,\\
   E_3\equiv &   3\alpha r + 3\beta v - 2\alpha\beta\eta + \beta^2\mu
   =0,\ &
   E_4\equiv &   \alpha s + \beta w - \beta^2\eta =0.\\
\ea
\ee
Without loss of generality we may consider $\alpha\beta\ne0$,
otherwise a rotation of the phase plane can be done. We have
$$
\ba{lll}
F_1\equiv \ &{\rm Res}_{\eta}\Big({\rm
Res}_{\mu}(E_1,E_2),E_4\Big)/\alpha&=
s\alpha^4+w\alpha^3\beta-3q\alpha^2\beta^2+(2p-3u)\alpha\beta^3+2t\beta^3=0,\\
\rule{0pt}{1.5em}
 F_2\equiv \
&{\rm Res}_{\eta}\Big({\rm Res}_{\mu}(E_1,E_3),E_4\Big)/\beta&=
2s\alpha^4+(2w-3r)\alpha^3\beta-3v\alpha^2\beta^2+p\alpha\beta^3+t\beta^3=0,\\
\ea
$$
where ${\rm Res}_z(f,g)$ denotes the resultant of the polynomials
$f$ and $g$ with respect the variable $z$, for more details on the
resultant see \cite{Olver}. Then \beq G_1(\alpha,\beta)\equiv
\frac{2F_2-F_1}{3\alpha}&=&
s\alpha^3+(w-2r)\alpha^2\beta + (q-2v)\alpha\beta^2 +u\beta^3=0,\\
G_2(\alpha,\beta)\equiv \frac{2F_1-F_2}{3\beta}&=&
r\alpha^3+(v-2q)\alpha^2\beta + (p-2u)\alpha\beta^2 +t\beta^3=0.
\eeq Now, it is clear that for the existence of three distinct
solutions $(\alpha_i,\beta_i)$ $(i=1,2,3)$ of the system
(\ref{eq_4}) it is necessary that the polynomials
$G_1(\alpha,\beta)$ and $G_2(\alpha,\beta)$ be proportional, i.e.
the following identity holds:
$$
 {\cal G}=\left|\ba{cc}\displaystyle \frac{\partial G_1}{\partial \alpha}&
               \displaystyle \frac{\partial G_1}{\partial
               \beta}\\[4mm]
               \displaystyle \frac{\partial G_2}{\partial \alpha}&
                \displaystyle \frac{\partial G_2}{\partial \beta}\\
                     \ea\right|=3\left(K_0\alpha^4+K_1\alpha^3\beta
+K_2\alpha^2\beta^2+K_3\alpha\beta^3+K_4\beta^4\right)=0,\\
$$
where
$$
\ba{l}
K_0 = -rw+2r^2-2qs+vs,\\
K_1 = 2ps-4us-2qr+4vr,\\
K_2 = 2q^2+2v^2+3ts-2pr+ur+pw-2uw-5qv,\\
K_3 = -4tr+2tw+4uq-2uv,\\
K_4 = 2u^2+tq-2tv-pu.
\ea
$$

On the other hand, the comitant ${\cal V}_3$ calculated for system
(\ref{s_1}) gives\\
\hbox{ }\hskip 2.5cm ${\cal V}_3= 2^5\cdot3^2\left(K_4x^4-K_3x^3y
+K_2x^2y^2-K_1xy^3+K_0y^4\right).$\\ Consequently, condition
${\cal G}= 0$ is equivalent to ${\cal V}_3=0$. Hence, we get the
next result.

\blm\label{lm:3} For cubic systems (\ref{s_1}) condition ${\cal
V}_3=0$ is necessary for the existence of three distinct couples
of parallel invariant straight lines.
\elm

Taking into account  Lemmas \ref{lm:2} and \ref{lm:3} we obtain
the next result.

\blm\label{lm:4}
If a cubic system (\ref{s_1}) possesses the configuration or the
potential configuration of
parallel invariant straight lines of the type (3,\,2,\,2), then it
is necessary ${\cal V}_3= {\cal V}_4= {\cal U}_{\,2}= 0$.
\elm

\bigskip

\noindent\textbf{ 3. Infinite singular points  and associated
homogeneous cubic canonical systems}

From \cite{Popa1} (see also \cite{Sib_Dini}) we have the following
result. Here $a\in\C$ is {\it imaginary} if $a\not\in\R$.

\blm\label{lm:5} The number of distinct roots (real and imaginary)
of the polynomial  $C_3=yp_3(x,y)-xq_3(x,y)\not=0$ is determined
by the following conditions:
\begin{description}
\item{[i]} $4$ real if ${\cal D}_1>0,\, {\cal D}_2>0,\,{\cal D}_3>0$;
\item{[ii]} $2$ real and $2$ imaginary if ${\cal D}_1<0$;
\item{[iii]} $4$ imaginary if ${\cal D}_1>0$ and for every $(x,y)$
where ${\cal D}_2{\cal D}_3\ne0$ either  ${\cal D}_2<0$ or ${\cal
D}_3<0$;
\item{[iv]} $3$ real ($1$ double, $2$ simple) if ${\cal D}_1=0,\,
{\cal D}_3>0;$
\item{[v]} $1$ real and $2$ imaginary ($1$ real double) if
${\cal D}_1=0,\, {\cal D}_3<0;$
\item{[vi]} $2$ real ($1$ triple and $1$ simple) if ${\cal D}_1=
{\cal D}_3=0,\, {\cal D}_2\not=0,\, {\cal D}_4=0$;
\item{[vii]} $2$ real ($2$ double) if ${\cal D}_1={\cal D}_3=0,\,
{\cal D}_2>0,\, {\cal D}_4\not=0$;
\item{[viii]} $2$ imaginary ($2$ double) if ${\cal D}_1={\cal D}_3=
0,\, {\cal D}_2<0,\, {\cal D}_4\not=0$;
\item{[ix]} $1$ real (of the multiplicity $4$) if ${\cal D}_1=
{\cal D}_2={\cal D}_3=0$.
\end{description}
where ${\cal D}_i\ $ for  $i=1,2,3,4$  are the T--comitants
defined in the introduction.
\elm

We consider the polynomial ${C_3}({\bf a},x,y)\not=0$ as a quartic
binary form. It is well known that there exists $g\in GL(2,\mathbb
R)$, $g(x,y)=(u,v)$, such that the transformed binary form
$g{C_3}(\ab ,x,y)=C_3(\ab,g^{-1}(u,v))$ is one of the following
$9$ canonical forms:
$$
\ba{clclcl}
 (i)& xy(x-y)(rx+sy),\ rs(r+s)\ne 0; & (iv)& x^2y(x-y); & (vii) &
   x^2y^2;\\
(ii)& x(sx+y)(x^2+y^2);    & (v) & x^2(x^2+y^2); &\ (viii)&
 (x^2+y^2)^2;\\
 (iii)& (px^2+qy^2)(x^2+y^2),\ pq>0; & (vi)& x^3y;
   & (ix)& x^4.\\
\ea
$$
We note that each of such canonical forms corresponds to one of
the cases enumerated in the statement of Lemma \ref{lm:5}.

On the other hand, applying  the same transformation $g$ to the
initial system and calculating for the transformed system its
polynomial $C_3(\ab(g),u,v)$  the following relation hold:
$$
  C_3(\ab(g),\, u,v)= \det(g)\, C_3(\ab,\, x,y) =
  \det(g)\, C_3(\ab ,\, g^{-1}(u,v))= \lambda C_3(\ab ,\,
  g^{-1}(u,v)),
$$
where we may consider $\lambda=1$ (via a time rescaling).

Taking into account  that $C_3(x,y)=yP_3(x,y)-xQ_3(x,y)$, we
construct the canonical forms of the cubic homogeneous systems
having their polynomials $C_3$ the indicated canonical forms
$(i)-(ix)$: \bq\label{pag:canon}
 &&\ba{ll}
    x'=(p+r)x^3+(s+v)x^2y+qxy^2, & C_3=xy(x-y)(rx+sy),\  \\
    y'=px^2y+(r+v)xy^2+(q+s)y^3, &  rs(r+s)\ne0\\
 \ea\label{HSys:1}\\[3mm]
 && \ba{ll}
    x'=(u+1)x^3+(s+v)x^2y+rxy^2, & C_3=x(sx+y)(x^2+y^2),  \\
    y'=-sx^3+ux^2y+vxy^2+(r-1)y^3, & \\
\ea\label{HSys:2}\\[3mm]
&&\ba{ll}
    x'=ux^3+(p+q+v)x^2y+rxy^2+qy^3, & C_3=(px^2+qy^2)(x^2+y^2),  \\
    y'=-px^3+ux^2y+vxy^2+ry^3,  & pq>0\\
\ea\label{HSys:3}\\[3mm]
&&\ba{ll}
    x'=3(u+1)x^3+(v-1)x^2y+rxy^2, & C_3=x^2y(x-y),  \\
    y'=ux^2y+vxy^2+ry^3, & \\
\ea\label{HSys:4}\\[3mm]
&&\ba{ll}
    x'=u x^3+(v+1)x^2y+rxy^2, & C_3=x^2(x^2+y^2),  \\
    y'=-x^3+ ux^2y+vxy^2+ry^3, & \\
\ea\label{HSys:5}\\[3mm]
&&\ba{ll}
    x'=(u+1)x^3+vx^2y+rxy^2, & C_3=x^3y,  \\
    y'=ux^2y+vxy^2+ry^3, & \\
\ea\label{HSys:6}\\[3mm]
&&\ba{ll}
    x'=u x^3+q x^2y+rxy^2, & C_3=(q-v)x^2y^2,  \\
    y'=ux^2y+vxy^2+ry^3,  & q-v\ne0          \\
\ea\label{HSys:7}\\[3mm]
&&\ba{ll}
    x'=u x^3+(v+1) x^2y+rxy^2 +y^3, & C_3=(x^2+y^2)^2,  \\
    y'=-x^3 + ux^2y+3(v-1)xy^2+ry^3,  &         \\
\ea\label{HSys:8}\\[3mm]
&&\ba{ll}
    x'=ux^3+vx^2y+rxy^2, & C_3=x^4,  \\
    y'=-x^3 +ux^2y+vxy^2+ry^3.  &   \\
\ea\label{HSys:9}
\eq

\noindent\textbf{ 4. Criteria for the existence of an invariant straight line
with a given multiplicity}

We consider a cubic system (\ref{s_1}) and the  associated four
polynomials $C_i(x,y)$ for $i=0,1,2,3$.

\bpr\label{pr:ap1} The straight line $\tilde L(x,y)= ux+vy=0$ is
invariant for a cubic system (\ref{s_1}) with $p_0^2+q_0^2\ne0$ if
and only if for $i=0,1,2,3$ the following relations hold:
\be\label{C:i}
   \mbox{either} \quad  C_i(-v,u)= 0,
\ee \be \label{fr:ap1} \mbox{or} \quad   {\rm
Res}_\gamma(C_0,C_i)=0 \quad ( \gamma = \frac{y}{x}
      \quad \mbox{or} \quad \gamma = \frac{x}{y}).
\ee
\epr

\BProof The line $\tilde L(x,y)$=0 is invariant for system
(\ref{s_1}) if and only if
$$
      u(p_0+ p_1+p_2+p_3)+v(q_0+q_1+q_2+q_3) = (ux+vy)(S_0+S_1+S_2),
$$
for some homogeneous polynomials $S_i$ of degree $i$. The last
equality is equivalent to
\beq
     && up_0+vq_0=0,  \\
     && up_1(x,y)+vq_1(x,y)=(ux+vy)S_0,  \\
     && up_2(x,y)+vq_2(x,y)=(ux+vy)S_1(x,y),\\
     && up_3(x,y)+vq_3(x,y)=(ux+vy)S_2(x,y).
\eeq If $x=-v, y=u$,  then the left--hand sides of the previous
equalities become $C_0(-v,u)$, $C_1(-v,u)$, $C_2(-v,u))$ and
$C_3(-v,u))$, respectively. At the same time the right--hand sides
of these identities vanish. Thus, we obtain equations (\ref{C:i})
in which $C_0$ (respectively, $C_1$; $C_2$; $C_3$) is a
homogeneous polynomial of degree 1 (respectively 2; 3; 4) in the
parameters  $u$ and $v$, and $C_0(x,y)\not=0$ because
$p_0^2+q_0^2\ne0$. Hence, the necessary and sufficient conditions
for the existence of a common solution of systems (\ref{C:i}) are
conditions (\ref{fr:ap1}).\EProof

Let $(x_0,y_0)\in \R^2$
 be an arbitrary point on the phase plane of
 systems (\ref{s_1}). Consider a translation $\tau$ bringing the
origin of coordinates to the point $(x_0,y_0)$.  We denote by
($\ref{s_1}{}^\tau$) the system obtained after applying the
transformation $\tau$, and by $\widetilde \ab = \ab(x_0,y_0)\in
\R^{20}$ the $20$--tuple of its coefficients. If $\gamma= y/x$ or
$\gamma= x/y$ then, for $i=1,2,3$ we denote
 \bq
\Omega_i(\ab,x_0,y_0)& = & {\rm Res}_\gamma
    \Big(C_i\big(\widetilde\ab,x,y\big),C_0\big(\widetilde\ab,x,y\big)\Big)
                         \in \R[\ab,x_0,y_0];\nonumber\\
  {\cal G}_i(\ab,x,y)&=&\left.\Omega_i(\ab,x_0,y_0)\right|_{\{x_0=x,\ y_0=y\}}\in \R[\ab,x,y].
  \label{Gam} \eq

\brm\label{obs} For $j=1,2,3$ the polynomials ${\cal G}_j(x,y)=
{\cal G}_j({\bf a},x,y)$ are affine comitants and are homogeneous
in the coefficients of system (\ref{s_1}) and non--homogeneous in
the variables $x$ and $y$. Additionally,
$$
\ba{lll}
  \deg_a {\cal G}_1=3, &\qquad \deg_a {\cal G}_2=4, &\qquad \deg_a {\cal
  G}_3=5,\\
  \deg_{(x,y)} {\cal G}_1=8, &\qquad \deg_{(x,y)} {\cal G}_2=10,
  &\qquad \deg_{(x,y)} {\cal G}_3=12.
\ea
$$
\erm

The geometrical meaning of these affine comitants is given by the
following lemma.

\blm\label{lm2} The straight line $L(x,y)= ux+vy+w=0$ is invariant
for a cubic system (\ref{s_1}) if and only if the polynomial
$L(x,y)$ is a common factor of the polynomials ${\cal G}_1$,
${\cal G}_2$ and ${\cal G}_3$ over $\C$.
\elm

\BProof Let $(x_0,y_0)\in \R^2$ be a non--singular point of system
(\ref{s_1}) (i.e. $P(x_0,y_0)^2+Q(x_0,y_0)^2\ne0$) which lies on
the line $L(x,y)=0$, i.e.
 $ux_0+vy_0+w=0$.
Denote by $\tilde L(x,y)=(L\circ\tau)~(x,y)= ux+vy$ ($\tau$ is a
translation) and consider the line $ux+vy=0$. By Proposition
\ref{pr:ap1}, the straight line $\tilde L(x,y)=0$ will be an
invariant line of systems ($\ref{s_1}{}^\tau$) if and only if
conditions (\ref{fr:ap1}) are satisfied for these systems, i.e.
for $i=1,2,3$, $\Omega_i(\ab,x_0,y_0)=0$, for each point
$(x_0,y_0)$ on the line $L(x,y)= ux+vy+w=0$. Thus, we have
$\Omega_i(\ab,x_0, y_0)= (ux_0+vy_0+w)\tilde \Omega_i(\ab,x_0,
y_0)$. Taking into account relations (\ref{Gam}) the lemma
follows. \EProof

\blm\label{lm3} If $L(x,y)= ux+vy+w=0$ is an invariant straight
line of (geometric) multiplicity $k$ for a cubic system
(\ref{s_1}) then, for $i=1,2,3$, we have that
$$
{\cal G}_i= (ux+vy+w)^k {W}_i(x,y).
$$
\elm

\BProof  By the definition of geometric multiplicity, we denote by
$(\ref{s_1}_\delta)$ the perturbed system from the system
(\ref{s_1}), which has $k$ invariant lines of multiplicity $1$:
$L_{i\delta}(x,y)$ for $i=1,\ldots,k$.

By Lemma \ref{lm2}, for $i=1,2,3$, system $(\ref{s_1}_\delta)$
satisfies ${\cal G}_{i\delta} = L_{1\delta}\cdots L_{k\delta}
\widetilde {W}_i(x,y)$, and when $\delta\to 0$ then
$L_{i\delta}(x,y) \to L(x,y)$. At the same time ${\cal
G}_{i\delta} \to {\cal G}_i = L(x,y)^k W_i$. \EProof

Taking into account Remark \ref{obs} and Lemmas \ref{lm2} and
\ref{lm3} we conclude the following result.

\blm\label{lm4} If a cubic system (\ref{s_1}) possesses the
maximum number of invariant straight lines (counted with their
multiplicities) then ${\cal G}_1(x,y)\mid {\cal G}_2(x,y)$\ and\ $
{\cal G}_1(x,y)\mid {\cal G}_3(x,y).$ \elm

In order to determine the degree of the common factor of the
polynomials ${\cal G}_i(x,y)$ for $i=1,2,3,$ we shall use the
notion of the $k^{th}$ {\sl subresultant}  of two polynomials with
respect to a given indeterminate (see for instance, \cite{Trudi},
\cite{Olver}).

We consider two polynomials \beq f(z)&=& a_0z^n+a_1z^{n-1}+\cdots+
a_n,\qquad g(z)= b_0z^m+b_1z^{m-1}+\cdots+ b_m, \eeq in variable
the $z$ of degree $n$ and $m$, respectively.

We say that the $k$--th {\it subresultant} with respect to
variable $z$ of the two polynomials $f(z)$ and $g(z)$ is the
$(m+n-2k)\times(m+n-2k)$ determinant
 \be\label{Res:k}
    R^{(k)}_z(f,g)= \left|\ba{llllll} a_0 & a_1 & a_2 &\ldots &  \ldots & a_{m+n-2k-1}\\
                    0  & a_0 & a_1 &\ldots & \ldots  & a_{m+n-2k-2}\\
                    0  &  0  &  a_0 &\ldots & \ldots  & a_{m+n-2k-3}\\
                  \ldots & \ldots &\ldots &\ldots &\ldots &\ldots\ \ldots\ \ldots \\
                    0  &  0  &  b_0 &\ldots & \ldots  & b_{m+n-2k-3}\\
                    0  & b_0 & b_1 &\ldots & \ldots  & b_{m+n-2k-2}\\
                  b_0 & b_1 & b_2 &\ldots &  \ldots &
                  b_{m+n-2k-1}\\
 \ea\right|\!\!\!\ba{ll}\left.\rule{0pt}{6ex}\right\}&\!\!\! (m-k)-times\\
                 \left.\rule{0pt}{6ex}\right\}&\!\!\! (n-k)-times\\
            \ea
 \ee
in which there are $m-k$ rows of $a$'s and $n-k$ rows of $b$'s,
and $a_i=0$ for $i>n$, and $b_j=0$ for $j>m$.

For $k=0$ we obtain the standard resultant of two polynomials. In
other words we can say that the $k$--th subresultant with respect
to the variable $z$ of the two polynomials $f(z)$ and $g(z)$ can
be obtained by deleting the first and the last $k$ rows and
columns from its resultant written in the form (\ref{Res:k}) when
$k=0$.

The geometrical meaning of the subresultants is based in the
following lemma.

\blm\label{Trudi:1}{\rm (see \cite{Trudi}, \cite{Olver}).}\
Polynomials $f(z)$ and $g(z)$ have precisely $k$ roots in common
(counting their multiplicities) if and only if the following
conditions hold:
$$
R_z^{(0)}(f,g)=R_z^{(1)}(f,g)=R_z^{(2)}(f,g)= \cdots=
R_z^{(k-1)}(f,g)= 0\ne R_z^{(k)}(f,g).
$$
\elm

For the polynomials in more than one variables it is easy to
deduce from Lemma \ref{Trudi:1} the following result.

\blm\label{Trudi:2} Two polynomials $\tilde f(x_1,x_2,...,x_n)$
and $\tilde g (x_1,x_2,...,x_n)$ have a common factor of degree
$k$ with respect to variable $x_j$ if and only if the following
conditions are satisfied:
$$
R_{x_j}^{(0)}(\tilde f,\tilde g)= R_{x_j}^{(1)}(\tilde f,\tilde
g)= R_{x_j}^{(2)}(\tilde f,\tilde g)= \cdots= R_{x_j}^{(k-1)}
(\tilde f,\tilde g)= 0\not=  R_{x_j}^{(k)}(\tilde f,\tilde  g),
$$
where $R_{x_j}^{(i)}(\tilde f,\tilde g)=0 $ in $\R[x_1,\ldots
x_{j-1},x_{j+1},\ldots, x_n].$ \elm

\bigskip

\noindent\textbf{ 5. Cubic systems with 4 real simple roots of $C_3$}

As it was shown above a cubic homogeneous system having 4 real
distinct infinite singular points via a linear transformation
becomes in the canonical form (\ref{HSys:1}). Therefore, in what
follows we consider the system \be\label{CF4_1} \ba{l}
    x'=p_0 +p_1(x,y)+p_2(x,y) +(p+r)x^3+(s+v)x^2y+qxy^2, \\
    y'=q_0 +q_1(x,y)+q_2(x,y) +px^2y+(r+v)xy^2+(q+s)y^3, \\
\ea \ee where the parameters $r$ and $s$ satisfy the condition
$rs(r+s)\ne0$. For system (\ref{CF4_1}) we obtain $C_3= x y (x-y)
(rx+sy),$ and hence, infinite singular points are situated at the
``ends'' of the following straight lines: $x=0$, $ y=0$, $x-y=0$
and $rx+sy=0$.

The goal of this section is to construct the cubic systems of the
form (\ref{CF4_1}) which have $8$ invariant straight lines with
the configuration $(3,\,3,\,1,\,1)$ or $(3,\,2,\,2,\,1)$.

\bigskip

\noindent\textbf{ 5.1. Systems with the configuration
$(3,\,3,\,1,\,1)$}

In this subsection we construct the cubic system with $4$ real
infinite singular points which possesses $8$ invariant affine
straight lines in the configuration or potential configuration
$(3,\,3,\,1,\,1)$, having total multiplicity $9$, as always the
invariant straight line of the infinity is considered.

According to Lemma \ref{lm:1} if a cubic system possesses $8$
invariant straight lines in the configuration $(3,\,3,\,1,\,1)$,
then necessarily the conditions ${\cal V}_1= {\cal V}_2= {\cal
U}_{\,1}=0$ hold.

A straightforward computation of the values of ${\cal V}_1$ and
${\cal V}_2$  for system (\ref{CF4_1}) yields:
$$
  {\cal V}_1=16\sum_{j=0}^4{\cal V}_{1j}x^{4-j}y^j,\qquad
  {\cal V}_2=8\sum_{j=0}^2{\cal V}_{2j}x^{2-j}y^j,
$$
where
\be\label{val_V1}
\ba{ll}
{\cal V}_{10}&=p(2p+3r),\\
{\cal V}_{11}&=2ps+4pv+2pr+3r^2+3rv,\\
{\cal V}_{12}&=4rs+4pq+3ps+3rq+2sv+2rv-s^2-r^2+2v^2,\\
{\cal V}_{13}&=2sq+4qv+3sv+2rq+3s^2,\\
{\cal V}_{14}&=q(2q+3s),\\
{\cal V}_{20}&=-3rv-3r^2+2ps-2pr,\\
{\cal V}_{21}&=6rq-2s^2-4sv+4rv+2r^2-6ps,\\
{\cal V}_{22}&=2qs+3s^2-2rq+3sv.\\
\ea
\ee
Consequently, relations ${\cal V}_1= {\cal V}_2=0$ provides  the
following equalities:
\be\label{val_U1}
\ba{llll}
{\cal V}_{10}&=p(2p+3r)=0, &  {\cal V}_{14}&=q(2q+3s)=0,\\
{\cal V}_{11}+{\cal V}_{20}&=4p(s+v)=0, &  {\cal V}_{13}-{\cal V}_{22}&=4q(r+v)=0,\\
{\cal V}_{11}-{\cal V}_{20}&=2(r+v)(2p+3r)=0, &  {\cal V}_{13}+{\cal V}_{22}&=2(s+v)(2q+3s)=0.\\
\ea \ee
Thus, we shall consider three cases: $(1)\ pq\ne0;$ $(2)\
pq=0,\, p^2+q^2\ne0;$ $(3)\ p=q=0.$

\noindent\textbf{ Case $p\,q\ne0$.} Then, from (\ref{val_U1}), we
obtain $v=-s$, $r=s$, $p=q=-3s/2\ne0$, and consequently  ${\cal
V}_1= {\cal V}_2= 0$, and ${\cal U}_{\,1}=0$. Therefore, by
changing the time ($t\to -2t/(3s)$) we obtain the following
system: \be\label{CS:1}
    x'=p_0 +p_1+p_2 + x^3+ 3xy^2, \quad
    y'=q_0 +q_1+q_2 +3 x^2y+ y^3,
\ee
for which ${\cal U}_{\,1}(\ab)=0$.

\noindent\textbf{ Case $p\,q=0,\, p^2+q^2\ne0$.} Then, without loss of
generality, we can consider $p=0$ and $q\ne0$ via the
transformation $x\leftrightarrow y$ and the changes
$p\leftrightarrow q$ and $r\leftrightarrow s$. From (\ref{val_U1})
we have $v=-r$, $q=-3s/2\ne0$, and
$$
\ba{l}
{\cal V}_{1i}=0\ (i=0,1,3,4),\ {\cal V}_{2j}=0,\ (j=0,2),\
 {\cal V}_{12}=4{\cal V}_{11}=-{\cal V}_{21}=-4(r+2s)(2r+s).
 \ea
$$
Consequently, we obtain either $s=-2r$, or $r=-2s$. The first case
after a suitable time rescaling writes the   system as
\be\label{CS:2}
    x'=p_0 +p_1+p_2 + x^3- 3x^2y+3xy^2, \quad
    y'=q_0 +q_1+q_2 + y^3, \\
\ee
whereas the second one goes over to the system
\be\label{CS:3}
    x'=p_0 +p_1+p_2+ 4x^3- 6x^2y+3xy^2, \quad
    y'=q_0 +q_1+q_2 + y^3.
\ee
We note that for both  systems we have ${\cal U}_{\,1}(\ab)=0$.

\noindent\textbf{ Case  $p=q=0$.} Then, by (\ref{val_U1}), we have
$r(r+v)=s(s+v)=0$. We claim  that  $rs\ne0$. Indeed, we suppose
$r=0$ (case $s=0$ can be reduced  to this one by changing
$x\leftrightarrow y$). Then, taking into account (\ref{val_V1}),
we obtain
$$
{\cal V}_{13}+6{\cal V}_{21}=12 s^2=0\ \Rightarrow\ {\cal
V}_{12}=16v^2=0.
$$
Thus   $r=s=v=0$ and we obtain $p_3(x,y)= q_3(x,y)=0$. Hence, the
claim  is proved. Considering (\ref{val_U1}), we obtain $t=-r=-s$,
and  after a suitable time rescaling  the system becomes
\be\label{CS:4}
    x'=p_0 +p_1+p_2 + x^3, \quad
    y'=q_0 +q_1+q_2 + y^3, \\
\ee
for which ${\cal U}_{\,1}(\ab)=0$.

\blm\label{lm:eqiv1} Systems (\ref{CS:1}), (\ref{CS:2}),
(\ref{CS:3}) and (\ref{CS:4}) are all affine equivalent. \elm

\BProof It is sufficient to check by straightforward computation
that the transformation $x_1=x-y$, $y_1=y$ writes system
(\ref{CS:4}) into system (\ref{CS:2}), and the transformation
$x_1=x$, $y_1=x-y$ writes system (\ref{CS:1}) into  system
(\ref{CS:3}). It remains to observe that the transformation
$x_1=x$, $y_1=y/2$ and $t_1=4t$ writes system (\ref{CS:3}) into
system (\ref{CS:2}). \EProof

Let $ L(x,y)= U x+V y+W=0 $ be an invariant straight line of
system (\ref{s_1}), which we write explicitly as:
$$
\ba{l}
\dot x= a + c x + d y + gx^2+2hxy+ky^2+px^3+3qx^2y+3rxy^2+sy^3,\\
\dot y= b + e x + f y + lx^2+2mxy+ny^2+tx^3+3ux^2y+3vxy^2+wy^3.\\
 \ea
$$
  Then,  we have
$$
      U P(x,y)+V Q(x,y)=(U x+V y+W)(Ax^2+2Bxy+Cy^2+Dx+Ey+F),
$$
and this identity provides the following 10 relations:
\be\label{eq_g}
\ba{ll}
Eq_1=pU+tV=0, & Eq_6=(2h\!-\!E)U\!+\!(2m\!-\!D)V\!-\!2BW\!=\!0,\\
Eq_2=(3q-2B)U+(3u-A)V=0,& Eq_7=kU+(n-E)V-CW=0,\\
Eq_3=(3r-CU)+(3v-2B)V=0, & Eq_8=(c-F)U+eV-DW=0\\
Eq_4=(s-C)U+VW=0,&  Eq_9= dU+(f-F)V-EW=0,\\
Eq_5=(g-D)U+lV-AW=0,& Eq_{10}=aU+bV-FW=0.
\ea
\ee

We concentrate our attention to the   system with 4 real distinct
infinite  singular points. According to Lemma \ref{lm:eqiv1}  we
can only work with  system (\ref{CS:4}). It is clear that via a
translation of the origin of coordinates at the point
$(-g/3,-n/3)$, we can consider the parameter $g=0$ (respectively,
$n=0$) in the polynomial $p_2$ (respectively, $q_2$). Thus, we
shall work with  the following   system \be\label{sys_c1} \dot x=
a + c x + d y + 2hxy+ky^2+x^3,\quad
\dot y= b + e x + f y + lx^2+2mxy+y^3,\\
\ee for which $C_3(x,y)=xy(x+y)(x-y)$. Therefore, there are the
following 4 directions for the possible invariant straight lines:
$x=0$, $y=0$,  $y=-x$, $y=x$.

We claim that in the direction $y=-x$ as well as in the direction
$y=x$ there can be only one invariant straight line. Indeed, for
the directions $y=-x$ and $y=x$ we have $U=1, V=\pm1$ and then,
from the first 6 equations (\ref{eq_g}), we obtain
$$
A^{\pm}=C^{\pm}=1,\ B^{\pm}=\mp1,\ D^{\pm}=\pm l-W,\ E^{\pm}=\pm
2W+2h-l\pm 2m,
$$
and $Eq_7=-3W \mp 2h+k\pm l-2m=0$. Here, the values with a
superindex $+$ (respectively $-$) correspond to $V=+1$
(respectively $V=-1$). So, from system (\ref{eq_g}), we can obtain
at most one solution $W_0^{\pm}$. Consequently, if system
(\ref{sys_c1}) possesses two couples of triples of parallel
invariant straight lines, then their directions only can be in the
directions $x=0$ and $y=0$. So, the claim is proved. Now we shall
investigate the conditions in order to have two couples of triples
of parallel invariant straight lines.

\noindent\textbf{ Direction $x=0$.} Then, $U=1,\ V=0$ and, from
(\ref{eq_g}), we obtain \beq
&&A=1,\ B=C=0,\ D=-W,\ E=2h,\ F=W^2+c,\\
&& Eq_7=k,\ Eq_9=2hW+d,\ Eq_{10}=-W^3-cW+a.
\eeq Thus, for the
existence of three solutions $W_i$ counted with their
multiplicity, it is necessary and sufficient that $k=h=d=0$.

\noindent\textbf{ Direction $y=0$.} In this case $U=0,\ V=1$ and, from
(\ref{eq_g}), we obtain \beq
&&A=B=0,\ C=1,\ D=2m,\ E=-W,\ F=W^2+f,\\
&& Eq_5=l,\ Eq_8=-2mW+e,\ Eq_{10}=-W^3-fW+b.
\eeq Hence,  for the
existence of three solutions $W_i$ counted with their multiplicity
it is necessary and sufficient that $l=m=e=0$.

Taking into account the conditions obtained  we have the following
system
\be\label{sys_c2} \ba{l} \dot x= a + c x + x^3,\quad
\dot y= b + f y + y^3.\\
 \ea
\ee
Now, it remains  to find out the conditions for the existence of
one invariant  straight line in each of the directions $y=-x$ and
$y=x$.

\noindent\textbf{ Direction $y\pm x=0$}. Considering equations
(\ref{eq_g}) for system (\ref{sys_c2}) in the directions $y\pm
x=0$ (i.e., $U=1,V=\pm 1$) we obtain: \beq
&&A^{\pm}=C^{\pm}=1,\ B^{\pm}=\mp1,\ D^{\pm}=-W,\ E^{\pm}=\pm 2W,\ F^{\pm}=W^2+c,\\
&& Eq_7^{\pm}=-3W,\ Eq_9^{\pm}=\mp 3W^2 \pm(f-c),\
Eq_{10}^{\pm}=-W^3+a\pm b.
\eeq
Thus, for both directions  the unique solution can be $W^{\pm}=0$,
and in order to  have in each direction an invariant straight line
it is necessary and sufficient that $f-c=a=b=0$. Thus, we have
obtained the    system
\be\label{sys_conf:1-3}
\ba{l}
\dot x=  c x + x^3,\quad
\dot y=  c y + y^3,\\
 \ea
\ee which possesses the invariant straight lines $ x=0,\
x=\pm\sqrt{-c},\ y=0,\ y=\pm\sqrt{-c},\ y=\pm x.$ It is clear that
the lines $x=\pm\sqrt{-c}$ (respectively, $y=\pm\sqrt{-c}$) are
real for $c<0$, imaginary for $c>0$, and coincide with the axes
for $c=0$. Hence we obtain  Figure 1 (respectively, 2; 3) for
$c<0$ (respectively, $c>0$; $c=0$).

\vspace{-2mm}
\brm\label{rm:transf} Assume $\alpha\in \R$. Then the transformation
$x=|\alpha|^{1/2}x_1$, $y=|\alpha|^{1/2}y_1$ and $t=|\alpha|^{-1}t_1$
does not change the coefficients of the  cubic part of the generic
cubic system. Whereas each coefficient of the quadratic (respectively,
linear; constant ) part will be multiplied by  $|\alpha|^{-1/2}$ (respectively,
by $|\alpha|^{-1}$; $|\alpha|^{-3/2}$).
\erm
\vspace{-2mm}
By Remark \ref{rm:transf} for system (\ref{sys_conf:1-3}) we can consider  $c\in \{-1,0,1\}$.

In order to obtain equivalent invariant conditions we shall use
the constructed T-comitants ${\cal L}_i$ $(i=1,2,3)$ and ${\cal
N}_1$. We note that the T-comitants ${\cal L}_i$ $(i=1,2,3)$ were
constructed by Calin \cite{Calin}.

For system (\ref{sys_c1}) we have ${\cal L}_1= -2^8\, 3^4
\left(lx^3+2mx^2y-2hxy^2-ky^3\right).$ Thus, conditions $k= h= l=
m= 0$ are equivalent to ${\cal L}_1=0$. Moreover, if for system
(\ref{sys_c1}) condition ${\cal L}_1=0$ holds, then we obtain $
{\cal L}_2=2^7\, 3^5x^2y^2 \left[-ex^2+6(f-c)xy+dy^2\right],
$\linebreak $ {\cal L}_3=2^9\,3^5 (c+f)\left(x^2+y^2\right). $
This means that condition ${\cal L}_2=0$ is equivalent to $d= e=
c-f= 0$ and, hence, since $f=c$ we obtain ${\rm sign} ({\cal
L}_3)= {\rm sign}(c)$. Therefore,  for ${\cal V}_1= {\cal V}_2=
  {\cal L}_1= {\cal L}_2=0$, we obtain the system
\be\label{sys_cc2}
\ba{l}
\dot x= a + c x + x^3,\quad
\dot y= b + c y + y^3.\\
 \ea
\ee
for which ${\cal N}_1=-2^33^5 xy(x^2+y^2)(ax-by).$

\brm\label{rm:N1} It is necessary  to underline that the
$GL$-comitants  ${\cal L}_i, \ i=1,2,3$ in fact are $T$-comitants
for the initial system (i.e. their coefficients are absolute
invariants under translations). But this is not the case for the
$GL$-comitant ${\cal N}_1$. However, for the   system
(\ref{sys_cc2}) ${\cal N}_1$ is a $CT$-comitant (see
\cite{Dana_Vlp1} for detailed definitions). \erm

\noindent{\it Proof}: Indeed, we  consider the system:
\be\label{sys_cc2t} \ba{ll}
\dot x_1=& a+c\gamma +\gamma ^3+ (c+3\gamma ^2)x_1+ 3\gamma x_1^2+x_1^3,\\
\dot y_1=& b+c\delta  +\delta  ^3+ (c+3\delta  ^2)y_1+ 3\delta  y_1^2+y_1^3 \\
 \ea
\ee which is obtained from system (\ref{sys_cc2}) via the
translation $x=x_1+\gamma $, $y=y_1+\delta  $, where $(\gamma
,\delta  )$ is an arbitrary point of the phase plane. For   system
(\ref{sys_cc2t}) we  calculate the value of the $GL$-comitant
${\cal N}_1= -2^33^5 x_1y_1(x_1^2+y_1^2)(ax_1-by_1).$ As we can
observe, the value of this polynomial does not depend on the
coordinates of the arbitrary point $(\gamma ,\delta  )$ and,
consequently for system (\ref{sys_cc2}) condition ${\cal N}_1=0$
is equivalent to $a= b= 0$ and this is an affine invariant
condition. \EProof

In short, we have the  following result.

\bpr\label{prop:1-3} A cubic system (\ref{sys_c1}) possesses
invariant straight lines with total multiplicity 9  if and only if
${\cal L}_1= {\cal L}_2= {\cal N}_1=0.$ Moreover, the
configuration or the potential configuration of the lines corresponds
 with $(3,3,1,1)$ given in Figure $1$ (respectively, $2$; $3$) for ${\cal L}_3$
negative (respectively, positive; zero). \epr

\bigskip

\noindent\textbf{ 5.2 Systems with configuration $(3,\,2,\,2,\,1)$}

In this subsection we construct the cubic systems with $4$ real
infinite singular points which possess $8$ distinct
invariant affine  straight lines with configuration or potential
configuration $(3,\,2,\,2,\,1)$, having total multiplicity $9$, as
always the invariant straight line of the infinity is considered.

For having the configuration $(3,\,2,\,2,\,1)$ a cubic system has
to possess three couples of parallel invariant lines and,
moreover, one couple must increase up to a triplet. Thus,
according to Lemma \ref{lm:4}, if a cubic system
possesses $8$ invariant straight lines in the configuration
$(3,\,2,\,2,\,1)$, then necessarily the conditions ${\cal V}_3=
{\cal V}_4= {\cal U}_2=0$ hold.

A straightforward computation of the value of ${\cal V}_3$ for
system (\ref{CF4_1}) yields: $\displaystyle {\cal V}_3= 32
\sum_{j=0}^4 {\cal V}_{3j}x^{4-j}y^j, $\ where
\be\label{val_V4}
\ba{l} {\cal V}_{30}=-p(p+3r),\quad
{\cal V}_{31}=2p(r-2s-v),\\
{\cal V}_{32}=4rs+3rq-sv-vr+3ps-2pq+2s^2+2r^2-v^2,\\
{\cal V}_{33}=-2q(2r-s+v),\quad {\cal V}_{34}=-q(q+3s).
\ea\ee

If $pq\ne0$, by (\ref{val_V4}), conditions ${\cal V}_{3i}= 0$,
$i=0,1,3,4$ yield $p= q= -3r= -3s= 3v\ne 0$, and then, condition
${\cal V}_{32}= -27v^2=0$ implies $v=0$, a contradiction. Thus,
condition $p q=0$ occurs, and we can suppose $q=0$; otherwise we
interchange $x\leftrightarrow y$, $p \leftrightarrow q$ and
$r\leftrightarrow s$. We consider two cases: $p=0$ and $p\ne0$.

\noindent\textbf{ Case $p=0$.} Then, we have ${\cal V}_{30}={\cal
V}_{31}=0$, and from (\ref{val_V4}) we obtain\linebreak $ {\cal
V}_{32}=(2s+2r+v)(s+r-v)=0. $

\noindent\textbf{ Subcase $v=s+r$.} For system (\ref{CF4_1}) we obtain
\be\label{val_W:1} {\cal V}_4 = 2^{10}\ 3^2
(s-r)(2r+s)(2s+r)xy(x-y)(rx+ys)=-\frac{1}{20}{\cal U}_{\,2}. \ee
Via the transformation $x\leftrightarrow y$, $r\leftrightarrow s$
and $v\leftrightarrow v$, we can consider only two subcases: $
r=s$ and $r=-2s$.

\noindent\textbf{ (a)} We assume $r=s$. Then, $v=2s$, and after a
suitable time rescaling we obtain the following system:
\be\label{CSS:1}
    x'=p_0 +p_1+p_2 + x^3 +3x^2y, \quad
    y'=q_0 +q_1+q_2 + 3 xy^2+ y^3. \\
\ee

\noindent\textbf{ (b)} If $r=-2s$, then $v=-s$, and  as above we
obtain the   system \be\label{CSS:2}
    x'=p_0 +p_1+p_2+ 2x^3, \quad
    y'=q_0 +q_1+q_2 + 3 xy^2- y^3.
\ee

\noindent\textbf{ Subcase $v=-2(s+r)$.} Then, for system (\ref{CF4_1})
we again obtain the values of the comitants ${\cal V}_4$ and
${\cal U}_{\,2}$  indicated in (\ref{val_W:1}). Thus, we consider
only two subcases: $ r=s$ and $r=-2s$.

\noindent\textbf{ (a)} If $r=s$ then $v=-4s$, and after a suitable
time rescaling we obtain the  system: \be\label{CSS:3}
    x'=p_0 +p_1+p_2 + x^3-3x^2y, \quad
    y'=q_0 +q_1+q_2 - 3 xy^2+ y^3. \\
\ee

\noindent\textbf{ (b)} Assume that $r=-2s$. Then, $v=2s$ and this
leads to the system
\be\label{CSS:4}
    x'=p_0 +p_1+p_2 -2x^3+3x^2y, \quad
    y'=q_0 +q_1+q_2 +  y^3. \\
\ee

\noindent\textbf{ Case $p\ne0$.} From (\ref{val_V4}) we obtain
$p=-3r$, $v=r-2s$. Then, for  system (\ref{CF4_1}) we again obtain
the values of the comitants ${\cal V}_4$ and ${\cal U}_{\,2}$
indicated in (\ref{val_W:1}). However, as $q=0$ and $p\ne0$  we
consider three subcases : $r=s$, $r=-2s$ and $s=-2r$.

\noindent\textbf{ Subcase $r=s$.} Then, we have $v=-s$, $p=-3s\ne0$,
and this provides the  system:
\be\label{CSS:5}
    x'=p_0 +p_1+p_2 +2x^3, \quad
    y'=q_0 +q_1+q_2 + 3x^2y- y^3.
\ee

\noindent\textbf{ Subcase $r=-2s$.}  Then $v=-4s$, $p=6s\ne0$, and via
a suitable change of the time we obtain  the  system:
\be\label{CSS:6}
    x'=p_0 +p_1+p_2 +4x^3-3x^2y, \quad
    y'=q_0 +q_1+q_2 + 3x^2y- y^3.
\ee

\noindent\textbf{ Subcase $s=-2r$.}  Then  $v=5r$, $p=-3r\ne0$, and we
get  the  system: \be\label{CSS:7}
    x'=p_0 +p_1+p_2 +2x^3-3x^2y, \quad
    y'=q_0 +q_1+q_2 + 3x^2y-6xy^2+ 2 y^3.
\ee

\blm\label{lm:eqiv2} Canonical systems
(\ref{CSS:1})--(\ref{CSS:7}) are all affine equivalent. \elm

\BProof To prove this assertion it is sufficient to verify that
the following changes go over to the systems as it is indicated:
$$
\ba{lll}
{}[x=x_1-y_1,\ y=-y_1]: & (\ref{CSS:5}) & \Rightarrow\quad (\ref{CSS:7});\\
 {}[x=y_1,\ y=2x_1,\ t=t_1/2]:\ & (\ref{CSS:4}) & \Rightarrow \quad
(\ref{CSS:2});\\
 {}[x=x_1-2y_1,\ y=-2y_1,\ t=t_1/2]:\ & (\ref{CSS:1}) & \Rightarrow \quad
(\ref{CSS:7});\\
{}[x=y_1-x_1,\ y=-x_1]: & (\ref{CSS:2}) & \Rightarrow\quad (\ref{CSS:7});\\
{}[x= x_1,\ y=-y_1]: & (\ref{CSS:3}) & \Rightarrow\quad (\ref{CSS:1});\\
 {}[x=x_1,\ y=2y_1,\ t=t_1/2]:\ & (\ref{CSS:6}) & \Rightarrow \quad
(\ref{CSS:7}).\\
\ea
$$
{}\EProof

\brm We note that for    system
(\ref{CSS:1}) the comitant  ${\cal V}_1=xy(x+y)^2\ne0$, and
according to Lemma \ref{lm:1} this system cannot possesses two
triples of parallel invariant lines.
\erm

\brm\label{W:2:3}
As it follows from (\ref{val_V4}) and (\ref{val_W:1}) for the
system with 4 real distinct infinite singular points the
conditions ${\cal V}_3= {\cal V}_4=0$ imply ${\cal U}_{\,2}=0$, as
well as ${\cal U}_{\,3}=0$, because for systems (\ref{CSS:1}) --
(\ref{CSS:7}) we have ${\cal U}_{\,3}=0$.
\erm

According to Lemma \ref{lm:eqiv2} we can  work for studying the
configuration $(3,\,2,\,2,\,1)$ with the system (\ref{CSS:2}). It
is clear that, via translation of the origin of coordinates at the
point $(-g/6,-m/3)$, we can consider $g=m=0$. Thus, we only need
to consider the system
\be\label{sys_c1_} \dot x= a + c x + d y +
2hxy+ky^2+2x^3,\quad
\dot y= b + e x + f y + lx^2+ny^2+3xy^2-y^3,\\
\ee
for which $C_3(x,y)=xy(x-y)(2x-y)$. Therefore, there are the
following $4$ directions for the possible invariant straight
lines: $x=0$, $y=0$, $y=x$, $y=2x.$

\noindent\textbf{ Direction $y=x$}. We show that this direction it can
be  only one invariant straight line. Indeed, for $U=1$ and $V=-1$
from the first $6$ equations (\ref{eq_g}) we obtain
$$
A=B=2,\ C=-1,\ D=-2W-l,\ E=-4W+2h-l,
$$
and $Eq_7=-3W+2h+k-l-n$. So, from system (\ref{eq_g}) we can
obtain at most one solution $W_0$.

\noindent\textbf{ Direction $x=0$.} Then $U=1,\ V=0$ and, from
(\ref{eq_g}), we obtain \beq
&&A=2,\ B=C=0,\ D=-W,\ E=2h,\ F=2W^2+c,\\
&& Eq_7=k,\ Eq_9=-2hW+d,\ Eq_{10}=-2W^3-cW+a.
\eeq
Thus, for the existence of at least two solutions $W_i$ it is
necessary $k=h=d=0$. Then  there exist 3 invariant straight lines
(which can coincide) in the direction $x=0$. So, in what follows
we shall suppose $k=h=d=0$.

\noindent\textbf{ Direction $y=0$.} In this case $U=0,\ V=1$ and, from
(\ref{eq_g}), we obtain
\beq
&&A=0,\ B=3,\ C=-1,\ D=-3W,\ E=W+n,\ F=-W^2-nW+f,\\
&& Eq_5=l,\ Eq_9=3W^2+e,\ Eq_{10}=W^3+nW^2-fW+b.
\eeq
Thus, we obtain $l=0$, and
$$
   Eq_9=3W^2+e=0,\quad Eq_{10}=W^3+nW^2-fW+b=0.
$$
For having two common solutions $W_i$, according to Lemma
\ref{Trudi:1}, we obtain the  relations:
$$
   R^{(0)}_{W}(Eq_9,Eq_{10})= -3(3b-en)^2-e(e+3f)^2=0,\quad
   R^{(1)}_{W}(Eq_9,Eq_{10})= 3(e+3f)=0.
$$
Hence, we obtain $e=-3f$ and $b=en/3$, and the  system (we set
$n_1=n/3$)
\be\label{sys_cc1_}
\dot x= a + c x + 2x^3,\quad \dot y= en_1 -3f x + f y +
3n_1y^2+3xy^2-y^3.
\ee

\noindent\textbf{ Direction $2x-y=0$.} Then  $U=2,\ V=-1$ and from
(\ref{eq_g}) we obtain
$$
\ba{l}
 A=2,\ B=1,\ C=-1,\
 D=-W+6n_1,\ E=-W+3n_1,\   F=-W^2+3n_1W+f,\\
  Eq_5=-12n_1,\
  Eq_9=3W^2-12n_1W+2c+f,\
  Eq_{10}=W^3-3n_1W^2-fW+2a-en_1.
\ea
$$
So, $n_1=0$ and we obtain the following system of equations
$$
   Eq_9=3W^2+2c+f=0,\quad Eq_{10}=W^3-fW+2a=0.
$$
For having two common solutions $W_i$, by Lemma \ref{Trudi:1}, we
obtain the relations:
$$
  R^{(0)}_{W}(Eq_9,Eq_{10})= -4\left[27a^2
 +(2c+f)(c+2f)^2\right]=0,\quad
   R^{(1)}_{W}(Eq_9,Eq_{10})= 6(c+2f)=0.
$$
Hence, we obtain $c+2f=a=0$ and, from (\ref{sys_cc1_}), the
following system \be\label{sys_conf:16-18} \dot x= -2f x +
2x^3,\quad
\dot y= -3f x + f y + 3xy^2-y^3,\\
\ee
which possesses the invariant straight lines $ x=0,\
x=\pm\sqrt{f},\ y=\pm\sqrt{f},\ y=x,\ 2x-y=\pm\sqrt{f}.$

We note, that by Remark \ref{rm:transf} we can consider $f=\in
\{-1,0,1\}$. Thus, we obtain that system (\ref{sys_c1_}) possesses
the indicated invariant straight lines if and only if $h= k= l= n= 0,\quad
d=c+2f=e+3f=0,\quad a=b=0.$ We shall prove that these conditions
are equivalent to $ {\cal L}_1= {\cal L}_2= {\cal N}_1= 0. $
Indeed, for system (\ref{sys_c1_}) we have $ {\cal
L}_1=2^9\,3^4\left[lx^3+2(h+n)x^2y-3hxy^2-ky^3\right]. $ Thus,
conditions $k=h=l=n=0$ are equivalent to ${\cal L}_1=0$. Next, for
system (\ref{sys_c1_}) with ${\cal L}_1=0$ we obtain ${\cal
L}_2=2^8\,3^5\left[(5e-8c-f)x^2+\right.$ $\left.(3c+7d+6f)xy-5dy^2\right]$ and
${\cal L}_3= -2^{10}\, 3^6(x-y)^2(c+f).$ Therefore, condition
${\cal L}_2=0$ is equivalent to $d= c+2f= e+3f= 0$, and for
$c=-2f$ we obtain that ${\rm sign}({\cal L}_3)= {\rm sign}(f)$.

Finally, for system (\ref{sys_c1_}) with ${\cal L}_1= {\cal
L}_2=0$ we calculate
$$
{\cal N}_1=2^23^5
xy(2x-y)\left[(2a+5b)x^2-(5a+2b)+(b-a)y^2\right],
$$
and hence, conditions $a=b=0$ are equivalent to ${\cal N}_1=0$. In
the same manner as  above (see Remark \ref{rm:N1}) it can be
verified that the value of the $CT$--comitant ${\cal N}_1$ will
not be changed after a translation of the origin of coordinates of
the phase plane of the system
$$
\dot x= a -2f x + 2x^3,\quad \dot y= b -3f x + f y + 3xy^2-y^3,
$$
at an arbitrary point $(\gamma ,\delta)$. Thus, considering Remark
\ref{W:2:3} we obtain the next result.

\bpr\label{prop:4-6} Cubic system (\ref{sys_c1_}) possesses
invariant straight lines with total multiplicity $9$ if and only
if ${\cal L}_1= {\cal L}_2= {\cal N}_1=0.$ Moreover, the
configuration or the potential configuration of the lines
corresponds with $(3,2,2,1)$ given in Figure $4$ (respectively,
$5$; $6$) for ${\cal L}_3$ positive (respectively, negative;
zero). \epr

From Lemma \ref{lm:5} and Propositions \ref{prop:1-3} and
\ref{prop:4-6} the next theorem follows.

\bth\label{th:1-6} For cubic system (\ref{s_1}) we assume that
the conditions ${\cal D}_1>0,$ ${\cal D}_2>0,$ ${\cal D}_3>0$ hold
(i.e. system has  4 real infinite singular points). Then, this
system will possess the maximum number of invariant straight lines
(with total multiplicity 9) if and only if one of the following
sequences of conditions holds:
$$
\ba{ll}
 {\cal V}_{1}={\cal V}_{2} = {\cal L}_{1} = {\cal L}_{2}= {\cal
 N}_{1}=0; &\hskip2cm ({\cal A}_1)\\
 {\cal V}_{3}={\cal V}_{4} = {\cal L}_{1} = {\cal L}_{2}= {\cal
 N}_{1}=0. &\hskip2cm ({\cal A}_2)\\
\ea
$$
Moreover, the configuration of the straight lines corresponds to
Figure 1, 2 or 3 for $({\cal A}_1)$ if ${\cal L}_3$ negative,
positive or zero, respectively; and to Figure 4, 5 or 6 for
$({\cal A}_2)$ if ${\cal L}_3$ positive, negative  or zero,
respectively.
\eth

\bigskip

\noindent\textbf{ 6. Cubic systems with 2 simple real and 2 simple
ima\-gi\-nary roots of $C_3$}

In this case, from (\ref{HSys:2}), the cubic system  after a
linear transformation becomes: \be\label{CF5_1} \ba{l}
    x'=p_0 +p_1(x,y)+p_2(x,y) +(u+1)x^3+(s+v)x^2y+rxy^2, \\
    y'=q_0 +q_1(x,y)+q_2(x,y) -sx^3+ux^2y+vxy^2+(r-1)y^3. \\
\ea \ee
For system (\ref{CF5_1}) we obtain $ C_3= x(sx+y)
(x^2+y^2)$, and hence, infinite singular points are situated at
the ``ends'' of the straight  lines: $x=0$, $ y=-sx$ and $y=\pm
ix$.

In this section we shall construct the cubic systems of the form
(\ref{CF5_1}) which can possess 8 invariant straight lines with
the configuration $(3,\,3,\,1,\,1)$ or $(3,\,2,\,2,\,1)$.

\bigskip

\noindent\textbf{ 6.1 Systems  with  the configuration
$(3,\,3,\,1,\,1)$}

In this subsection we construct the cubic system with $2$ real and
$2$ imaginary infinite singular points which possess $8$
invariant affine  straight lines with configuration or potential
configuration $(3,\,3,\,1,\,1)$, having total multiplicity $9$, as
always the invariant straight line of the infinity is considered.

By Lemma \ref{lm:1}, if a cubic system possesses $8$ invariant
straight lines in the configuration $(3,\,3,\,1,\,1)$, then  the
conditions ${\cal V}_1= {\cal V}_2= {\cal U}_{\,1}=0$ hold. A
straightforward computation of the values of ${\cal V}_1$ and
${\cal V}_2$ for system (\ref{CF5_1}) yields:\quad $\displaystyle
{\cal V}_1=16\sum_{j=0}^4{\cal V}_{1j}x^{4-j}y^j,$\quad $
\displaystyle {\cal V}_2=8\sum_{j=0}^2{\cal V}_{2j}x^{2-j}y^j, $\
where
\be\label{val_VV1} \ba{llll}
{\cal V}_{10}&=3s^2+2u^2+3u,& {\cal V}_{14}&=r(2r-3),\\
{\cal V}_{11}&=2su+4uv+9s+3v,& {\cal V}_{20}&=12sr-9s-3v+2su,\\
{\cal V}_{12}&=2sv+4ur-s^2+3r-3u+2v^2,&{\cal V}_{21}&=6r-2s^2-4sv+6u,\\
{\cal V}_{13}&=2r(s+2v)-3(s+v), & {\cal V}_{22}&=2sr-3s-3v.\\
\ea\ee
Consequently, relations ${\cal V}_1= {\cal V}_2=0$ imply
\be\label{val_VV2}
\ba{l}
{\cal V}_{14}=r(2r-3)=0,\quad   {\cal V}_{13}-{\cal V}_{22}=4rv=0,\quad {\cal V}_{11}+{\cal V}_{20}=2(2r-3)(s+v)=0. \\
\ea\ee
 Thus, we shall consider two cases: $r=0$ and $ r\ne0.$

\noindent\textbf{ Case $r=0$.} From (\ref{val_VV2}) we obtain $v=-s$
and, then the equation $ {\cal V}_{10}+{\cal V}_{12}= 2(u^2+s^2)=0
$ implies $u=s=0$, and after the time rescaling $t\to t/3$, we get
the system:
\be\label{CSS_3}
    x'=p_0 +p_1+p_2 +x^3, \quad
    y'=q_0 +q_1+q_2 -y^3, \\
\ee
for which ${\cal U}_{\,1}(\ab)=0$.

\noindent\textbf{ Case $r\ne0$.} Then, from (\ref{val_VV2}) we obtain
$v=0$, $r=3/2$ and then we have $ 2{\cal V}_{10}+ {\cal V}_{21}=
(2u+3)^2+s^2= 0.$ Hence, $s=0$, $u=-3/2$, and after a time
rescaling $t\to -2t/3$, we obtain the system
\be\label{CSS_4}
    x'=p_0 +p_1+p_2 +x^3-3xy^2, \quad
    y'=q_0 +q_1+q_2 +3x^2y -y^3, \\
\ee
for which ${\cal U}_{\,1}(\ab)=0$. It is easy to check  that
systems (\ref{CSS_3}) and (\ref{CSS_4}) are not linearly
equivalent.  This immediately follows from the following remark.

\brm\label{rm:L5} For  systems (\ref{CSS_3}) and (\ref{CSS_4}) we
have $ {\cal L}_4= \alpha(x^4+6x^2y^2+y^4), $ where $\alpha=-6$
for system (\ref{CSS_3}) and $\alpha=24$ for  system
(\ref{CSS_4}). Since the T--comitant ${\cal L}_4$ is of degree 2
in the coefficients of the system and is of zero weight (i.e. is
an  absolute T--comitant), when conditions ${\cal V}_1= {\cal
V}_2=0$ are satisfied we obtain system (\ref{CSS_3}) for ${\cal
L}_4<0$ and  system (\ref{CSS_4}) for ${\cal L}_4>0$.
\erm

As it has been  proved above we must examine
systems (\ref{CSS_3}) and (\ref{CSS_4}).

\noindent\textbf{ System (\ref{CSS_3})} Evidently via translation of
the origin of coordinates at the point $(-g/3,-n/3)$ we can
consider $g=n=0$. Thus, we examine the system
\be\label{sys_c1:3}
\dot x= a + c x + d y + 2hxy+ky^2 +x^3,\quad
\dot y= b + e x + f y + lx^2+2mxy -y^3,\\
\ee for which $C_3(x,y)=xy(x^2+y^2)$. Therefore, there are the
following $4$ directions for the possible invariant straight
lines: $x=0$, $y=0$, $y=\pm ix.$ We show that in each imaginary
direction $y=\pm ix $ there is only one invariant  straight line.
Indeed, to show this it is sufficient to examine only  the
direction $x+iy=0$, i.e. $U=1$ and $V=i$. From the first 6
equations (\ref{eq_g}) we obtain
$$
A=1,\ B=-i,\ C=-1,\ D=-W-1,\ E= 2h+i(2W+2m+1),
$$
and $Eq_7=3W+2m+k+1-2ih$. So, from (\ref{eq_g}) we obtain at most
one solution $W_0$. We examine other directions in each of which
it has to be one triple of parallel lines.

\noindent\textbf{ Direction $x=0$.}  Then, $U=1,\ V=0$ and, from
(\ref{eq_g}), we obtain
\beq
&&A=1,\ B=C=0,\ D=-W,\ E=2h,\ F=W^2+c,\\
&& Eq_7=k,\ Eq_9=-2hW+d,\ Eq_{10}=-W^3-cW+a.
\eeq
Thus, for the
existence of three solutions $W_i$ we must have $k= h= d= 0$,
these three solutions can coincide in the direction $x=0$. So, in
what follows we shall suppose $k= h= d= 0$.

\noindent\textbf{ Direction $y=0$.} In this case $U=0,\ V=1$ and,
again  from  (\ref{eq_g}), we obtain
\beq
&&A=B=0,\ C=-1,\ D=2m,\ E=W,\ F=-W^2+f,\\
&& Eq_5=l,\ Eq_8=-2mW+e,\ Eq_{10}=W^3-fW+b.
\eeq
Hence, for the existence of one triple of parallel lines in this
direction it is necessary $l=m=e=0$. Thus, we obtain the system
\be\label{sys_c1:3_}
\dot x= a + c x + x^3,\quad
\dot y= b +  f y  -y^3,\\
\ee
for which is necessary to examine the imaginary directions.

\noindent\textbf{ Direction $y+ix=0$.} We have $U=1,\ V=i$ and,
 from  (\ref{eq_g}), we obtain
\beq
&&A=1,\ B=-i,\ C=-1,\ D=-W,\ E=2iW,\ F=W^2+W +c,\\
&& Eq_7=3W,\ Eq_8=3W^2+c-f,\ Eq_{10}=2W^3-fW+a+ib.
\eeq
So, there can exist a unique solution $W_0=0$. Hence, since the
 system is real we get up the conditions $c-f=a=b=0$.
Thus, we obtain the  system
\be\label{sys_conf:4-5}
\dot x=  c x + x^3,\quad
\dot y=  c y  -y^3,\\
\ee
which possesses the invariant lines: $ x=0,\ y=0,\
x=\pm\sqrt{-c},\ y=\pm\sqrt{c},\ y=\pm ix. $

We can consider $c\in \{0,1\}$. Indeed, firstly we can suppose
$c\ge0$, otherwise the transformation $x\leftrightarrow y$,
$t\leftrightarrow -t$ can be used. Then, for $c>0$, we can apply
Remark~\ref{rm:transf}. Thus, we obtain that  system (\ref{sys_c1:3})
possesses the indicated invariant straight lines if and only if
the following conditions hold: $h= k= l= m= 0$, $d= e= c-f= 0$,
$a= b= 0.$ We prove that these conditions are equivalent to $
{\cal L}_1= {\cal L}_2= {\cal N}_1= 0. $ Indeed, for system
(\ref{sys_c1:3}) we have $ {\cal L}_1= 2^8\, 3^4\left(lx^3+
2mx^2y- 2hxy^2- ky^3\right).$ Hence, conditions $k=h=l=m=0$ are
equivalent to ${\cal L}_1=0$. Next, for system (\ref{sys_c1:3})
with ${\cal L}_1=0$ we obtain $ {\cal L}_2= 2^7\,3^5 \left[ex^2+
6(c-f)xy- dy^2\right]$, $ {\cal L}_3= 2^{9}\,3^5(c+f) \left(x^2-
y^2\right)$. Therefore, condition ${\cal L}_2= 0$ is equivalent to
$d= e= c-f= 0$, and for $f=c$ we obtain that condition $c=0$ is
equivalent to ${\cal L}_3=0$.

Finally, if for  system (\ref{sys_c1:3}) the conditions ${\cal
L}_1= {\cal L}_2=0$ hold, then we obtain the system
\be\label{sys:N1}
\dot x= a+ c x + x^3,\quad
\dot y= b+ c y  -y^3,\\
\ee
for which  we calculate
$
{\cal N}_1= 2^3\,3^5xy\left(x^2-y^2\right)\left(ax-by\right).
$
Moreover, for the system
$$
\ba{l}
\dot x= a+ c\gamma +\gamma ^3+(c+3\gamma ^2)x_1+3\gamma x_1^2 + x_1^3,\\
\dot y= b+c\delta  -\delta  ^3+ (c+3\delta  ^2)y_1-3\delta  y_1^2  -y_1^3,\\
 \ea
$$
which is obtained from system (\ref{sys:N1}) via translation
$x=x_1+\gamma $, $y=y_1+\delta  $, where $(\gamma ,\delta  )$ is
an arbitrary point of the phase plane of system (\ref{sys:N1}), we
have
$$
{\cal N}_1\left(\ab(\gamma ,\delta  ),x_1,y_1)\right)=
2^3\,3^5x_1y_1\left(x_1^2-y_1^2\right)\left(ax_1-b_1y\right).
$$
As we can observe that the value of this polynomial does not
depend on the coordinates of the arbitrary point $(\gamma ,\delta
)$ and, consequently for system (\ref{sys:N1}) condition ${\cal
N}_1=0$ is equivalent to $a=b=0$, and this condition  is affine
invariant. Thus, we have the next result.

\bpr\label{prop:7-8} A cubic  system (\ref{CSS_3}) possesses
invariant straight lines with total multiplicity 9  if and only if
${\cal L}_1= {\cal L}_2= {\cal N}_1=0.$ Moreover, the
configuration or the potential configuration of the lines
corresponds with $(3,3,1,1)$ given in Figure $7$ for ${\cal
L}_3\not=0$, and in Figure $8$ for ${\cal L}_3=0$.
\epr

\noindent\textbf{ System (\ref{CSS_4})} Doing a translation of the
origin of coordinates at the point $(-g/3,-n/3)$, we can consider
$g=n=0$. Thus, we have the  system
\be\label{sys_CSS_4_1} \ba{l}
\dot x= a + c x + d y + 2hxy+ky^2 +x^3-3 xy^2,\\
\dot y= b + e x + f y + lx^2+2mxy + 3x^2y-y^3,\\
 \ea
\ee for which $C_3(x,y)=-2xy(x^2+y^2)$. Therefore, there are the
following 4 directions for possible invariant straight lines:
$x=0$, $y=0$, $y=\pm ix$. We shall show that in the direction
$x=0$ as well as in the direction $y=0$ there can exist at most
one invariant straight line. Indeed, for the direction $x=0$ (i.e.
$U=1,\ V=0$) from the first 7 equations (\ref{eq_g}) we obtain
$A=1,\ B=0,\ C=-3,\ D=-W,\ E=2h,\ Eq_7=3W+k.$ Whereas for the
direction $y=0$ (i.e. $U=0,\ V=1$) we have $ A=3,\ B=0,\ C=-1,\
D=2m,\ E=W,\ Eq_5=-3W+l.$ Thus, in both cases there can exist only
one line in each considered direction.

We examine the imaginary directions, and since the system is real
it is sufficient  to consider only one direction.

\noindent\textbf{ Direction $x+iy=0$.}  Then $U=1,\ V=i$ and, from
(\ref{eq_g}) we obtain
\beq
&&A=1,\ B=2i,\ C=-1,\ D=-W+il,\ E=-i(W+k),\ F=W^2-ilW +c+ie,\\
&& Eq_6=l+2h+i(k+2m),\ Eq_9=(ik-l)W + i(f-c)+d+e,\\
&& Eq_{10}=-W^3+ilW^2 -(c+ie)W+a+ib.
\eeq
Thus, for the existence of three solutions $W_i$ it is necessary
$l+2h=k+2m=k=l=c-f=d+e=0$, i.e. $l=k=h=m=0$, $f=c$ and $d=-e$. In
this  case  we get the  system
$$
\dot x= a + c x - e y + x^3-3 xy^2,\quad
\dot y= b + e x + c y  + 3x^2y-y^3,\\
$$
and we have to examine the real directions.

\noindent\textbf{ Direction $x=0$.} In this case $U=1,\ V=0$, and from
(\ref{eq_g}), we obtain
\beq
&&A=1,\ B=0,\ C=-3,\ D=-W,\ E=0,\  F=W^2+c,\\
&&Eq_7=3W,\ Eq_9=-e,\ Eq_{10}=-W^3 -cW+a.
\eeq So, there can exist
only one solution $W_0=0$, and for the existence of one invariant
straight line with the direction $x=0$ it is necessary and
sufficient $e=a=0$.

\noindent\textbf{ Direction $y=0$.}  Then $U=0,\ V=1$ and, from
(\ref{eq_g}), we obtain \beq
&&A=3,\ B=0,\ C=-1,\ D=2m,\ E=W,\ F=-W^2+c, \\
 &&Eq_5=-3W,\quad  Eq_{10}=W^3 -cW+b.
\eeq So, again there is only one solution $W_0=0$, and for the
existence of the invariant straight line with the direction $y=0$
it is necessary and sufficient $b=0$. Thus, we obtain the system
\be\label{sys_conf:6-7} \dot x=  c x + x^3-3xy^2,\quad
\dot y=  c y +3x^2y -y^3,\\
\ee
which possesses the invariant straight lines: $ x=0,\ y=0,\
x+iy=\pm\sqrt{-c},\ x-iy=\pm\sqrt{-c},\ x\pm iy=0. $ We  can
consider $c\in \{0,1\}$. Indeed, firstly we can suppose $c\ge0$,
otherwise the transformation $x\leftrightarrow y$,
$t\leftrightarrow -t$ can be used. Then, for $c>0$ we can apply
Remark \ref{rm:transf}. Thus, we obtain that system (\ref{sys_CSS_4_1})
possesses the indicated invariant straight lines if and only if $
h=k=l=m=0,\quad d=e=c-f=0,\quad a=b=0.$ We shall prove that these
conditions are equivalent to $ {\cal L}_1= {\cal L}_2= {\cal N}_1=
0. $ Indeed, for  system (\ref{sys_CSS_4_1}) we have $ {\cal
L}_1=2^8\,3^4\left[(2h+3l)x^3+(k+6m)x^2y-(6h+l)xy^2-(3k+2m)y^3\right].
$ Hence, conditions $k=h=l=m=0$ are equivalent to ${\cal L}_1=0$.
Next, for  system (\ref{sys_CSS_4_1}) with ${\cal L}_1=0$ we
obtain
$$
{\cal
L}_2=2^8\,3^5\left[(7d-5e)x^2+2(c-f)xy+(5d-7e)y^2\right],\quad
{\cal L}_3= 2^{14}\,3^5(c+f)\left(x^2-y^2\right).
$$
Therefore, condition ${\cal L}_2=0$ is equivalent to $d=e=c-f=0$,
and for $f=c$ condition $c=0$ is equivalent to ${\cal L}_3=0$.

Finally, system (\ref{sys_CSS_4_1}) for ${\cal L}_1= {\cal L}_2=0$
becomes the system
\be\label{sys:N1b}
\dot x=  a+ c x + x^3-3xy^2,\quad
\dot y=  b+ c y +3x^2y -y^3,\\
\ee
for which  we calculate
$
{\cal N}_1= 2^4\,3^5\left(x^4-y^4\right)\left(ax+by\right).
$
Moreover, for the system
$$
\ba{l}
\dot x\!=\! (a+ c\gamma +\gamma ^3\!-\!3\gamma \delta
^2)+(c+3\gamma ^2\!-\!3\delta  ^2)x_1\!-\!6\gamma \delta  +3\gamma
x_1^2\!-\!6\delta
x_1y_1\!-\!3\gamma y_1^2 +x_1^3\!-\!3x_1y_1^2,\\
\dot y\!=\! (b+c\delta  +3\gamma ^2\delta  -\delta  ^3)+
(c+3\gamma ^2-3\delta  ^2)y_1
  +3\delta  x_1^2+6\gamma x_1y_1-3\delta  y_1^2+3x_1^2y_1-y_1^3,\\
 \ea
$$
which is obtained from system (\ref{sys:N1b}) via the translation
$x=x_1+\gamma $, $y=y_1+\delta  $, where $(\gamma ,\delta  )$ is
an arbitrary point of the phase plane of system (\ref{sys:N1b}),
we have $ {\cal N}_1(\ab(\gamma ,\delta  ),x_1,y_1)= 2^4\,3^5
\left(x_1^4-y_1^4\right)\left(ax_1+by_1\right). $ So, the value of
this polynomial does not depend on the coordinates of the
arbitrary point $(\gamma ,\delta  )$ and, consequently for system
(\ref{sys:N1b}) condition ${\cal N}_1=0$ is equivalent to $a=b=0$,
and this condition  is affine invariant. Thus, we get the
following result.

\bpr\label{prop:9-10} A cubic system (\ref{CSS_4}) possesses
invariant straight lines with total multiplicity 9  if and only if
$ {\cal L}_1= {\cal L}_2= {\cal N}_1=0.$ Moreover, the
configuration or the potential configuration of the lines
corresponds with $(3,3,1,1)$ given in Figure $9$ for ${\cal
L}_3\not=0$ or in Figure $10$ for ${\cal L}_3=0$.
\epr

\bigskip

\noindent\textbf{ 6.2 Systems with configuration $(3,\,2,\,2,\,1)$}

In this subsection we construct the cubic systems with $2$ real
and $2$ imaginary infinite singular points which possess $8$
 invariant affine  straight lines with configuration or potential
configuration $(3,\,2,\,2,\,1)$, having total multiplicity $9$.

By Lemma \ref{lm:4} if a cubic system possesses
$8$ invariant straight lines with configuration $(3,\,2,\,2,\,1)$,
then the conditions ${\cal V}_3= {\cal V}_4= {\cal U}_{\,2}=0$
hold.

We consider again system (\ref{CF5_1}). A straightforward
computation of the value of ${\cal V}_3$ for this system yields:
$\displaystyle {\cal V}_3=32\ \sum_{j=0}^4{\cal
V}_{3j}x^{4-j}y^j,$\ \ where \be\label{val_V4a} \ba{l} {\cal
V}_{30}=-3s^2-3u+3sv-u^2,\quad
{\cal V}_{31}=6sr-4su-2uv-18s,\\
{\cal V}_{32}=-sv-2ur-9+2s^2+3r-3u-v^2,\\
{\cal V}_{33}=2r(s-v),\qquad {\cal V}_{34}=r(3-r).
\ea\ee
So, we need to consider the cases: $r\ne0$ and $r=0$.

If $r\ne0$ then, from (\ref{val_V4a}), we have $v=s,$ $r=3$ and,
then
$$
 {\cal V}_{30}=-u(u+3)=0,\ {\cal V}_{31}=-6su=0,\  {\cal V}_{32}=-9u=0.
$$
Hence, $u=0$ and for system (\ref{CF5_1}) we have:
\be\label{val_W1a} {\cal V}_4 = 2^{11}\ 3^2
s(s^2+9)x(x^2+y^2)(sx+y)=-\frac{1}{20}{\cal U}_{\,2}. \ee
Consequently condition ${\cal V}_4= 0$ implies $s=0$, and this
provides, after the rescaling of the time $t\to t/3$, the system
\be\label{sys_SS4a}
    x'=p_0 +p_1+p_2 +x^3+3xy^2, \quad
    y'=q_0 +q_1+q_2 +2y^3. \\
\ee

For $r=0$, from (\ref{val_V4a}), we obtain ${\cal V}_{32}=
-sv-9+2s^2-3u-v^2=0.$ Hence $u=(2s^2-sv-v^2-9)/3$. Then, by
(\ref{val_V4a}), we obtain
$$
{\cal V}_{30}= -\frac{1}{9}(s-v)^2\left[(2s+v)^2+9\right]=0,\quad
{\cal V}_{31}= -\frac{2}{3}(s-v)\left[(2s+v)^2+9\right]=0,\
$$
that implies $v=s$, and then $u=-3$. Now, calculating the
polynomials ${\cal V}_4= 2^{11}\ 3^2
s(s^2+9)x(x^2+y^2)(sx+y)=-{\cal U}_{\,2}/20=0$ we obtain  $s=0$.
Thus, via the time rescaling $t\to -t/3$ we obtain the system
\be\label{sys_SS4b}
    x'=p_0 +p_1+p_2 +2x^3, \quad
    y'=q_0 +q_1+q_2 +3x^2y+y^3. \\
\ee
Evidently,  systems (\ref{sys_SS4a}) and (\ref{sys_SS4b}) are
affine equivalent via the transformation $x\leftrightarrow y$.

We consider  system (\ref{sys_SS4b}). It is clear that via
translation of the origin of coordinates at the point
$(-g/6,-n/3)$ we can consider $g=n=0$. Thus, we shall examine the
system \be\label{sys_cs4b} \ba{l}
\dot x= a + c x + d y + 2hxy+ky^2+2x^3,\\
\dot y= b + e x + f y + lx^2 + 2mxy + 3x^2y+y^3,\\
 \ea
\ee for which $C_3(x,y)=-xy(x^2+y^2)$. Therefore, there are the
following $4$ directions for the possible invariant straight
lines: $x=0$, $y=0$, $y=\pm ix.$ Since we are looking for the
configuration $(3,2,2,1)$ in system (\ref{sys_cs4b}), it follows
that the two couples of parallel invariant lines must be  in the
imaginary directions.

First we show that in the direction $y=0$ there cannot be a triple
of parallel lines. Indeed, for $U=0$ and $V=1$, from the first
five equations of (\ref{eq_g}), we obtain
$$
A=3,\ B=0,\ C=1,\ Eq_5=-3W+l=0
$$
and, hence, there can exist at most one solution of system
(\ref{eq_g}).

\noindent\textbf{ Direction $x=0$.} Then $U=1,\ V=0$ and, from
(\ref{eq_g}), we obtain \beq
&&A=2,\ B=C=0,\ D=-2W,\ E=2h,\ F=2W^2+c,\\
&& Eq_7=k,\ Eq_9=-2hW+d,\ Eq_{10}=-2W^3-cW+a. \eeq Thus, for the
existence of three solutions $W_i$ (which can coincide) it is
necessary and sufficient $k=h=d=0$, and in what follows we shall
assume  that these conditions hold.

\noindent\textbf{ Direction $y+ix=0$.} We have $U=1,\ V=i$ and, from
(\ref{eq_g}), we obtain \beq
&&A=2,\ B=i,\ C=1,\ D=-2W+il,\ E=iW,\ F=-W^2 +f,\\
&& Eq_6=l+2im,\ Eq_8=3W^2-ilW+ c-f+ie,\ Eq_{10}=W^3-fW+a+ib. \eeq
Thus, since the considered cubic system is real, for the existence
of two solutions $W_i$ it is necessary that $l=m=0$. Then,
according to Lemma \ref{Trudi:1}, for having two common solutions
$W_i$ it is necessary $R^{(1)}_{W}(Eq_8,Eq_{10})= 3(c+2f+ie)=0.$
Hence, since the system is real, the last condition yields $e=0$,
$c=-2f$. Therefore, $R^{(0)}_{W}= 27(a+ib)^2=0$, that gives
$a=b=0$. It remains to note that the obtained conditions are
sufficient for the existence of an invariant straight line in the
real direction $y=0$, as it is observed below. Thus, we obtain the
system
\be\label{sys_conf:19-21} \dot x=  -2f x +2x^3,\quad
\dot y=  f y + 3x^2y+y^3,\\
\ee which possesses the  invariant straight lines: $x=0$, $y=0$,
$x=\pm\sqrt{f}$, $x-iy=\pm\sqrt{f}$, $x+iy=\pm\sqrt{f}$. It is
clear that the lines $x=\pm\sqrt{f}$ (respectively, $y= \pm
\sqrt{f}$) are real for $f>0$, imaginary for $f<0$, and coincide
for $f=0$. Hence, we obtain  Figure 11 (respectively, 12 and 13)
for $f>0$ (respectively for $f<0$ and $f=0$). We note, that by Remark
\ref{rm:transf} we can consider $f\in \{-1,0,1\}$.

In short, we obtain that system (\ref{sys_cs4b}) possesses the
indicated invariant straight lines if and only if $h=k=l=m=0$,
$d=e=c+2f=0$, $a=b=0.$ We shall prove that these conditions are
equivalent to $ {\cal L}_1= {\cal L}_2= {\cal N}_1= 0. $ Indeed,
for system (\ref{sys_cs4b}) we have\quad $ {\cal L}_1= 2^9\,3^4
\left[(h-l)x^3-2mx^2y+3hxy^2+ky^3\right]. $ Thus, conditions
$k=h=l=m=0$ are equivalent to ${\cal L}_1=0$. Next, for system
(\ref{sys_cs4b}) with ${\cal L}_1=0$, we obtain
$$
{\cal L}_2=2^8\,3^5\left[(7d-5e)x^2-3(c+2f)xy+5dy^2\right],\quad
{\cal L}_3=2^{10}\,3^6(c+f)y^2.
$$
Hence, condition ${\cal L}_2=0$ is equivalent to $d=c+2f=e+3f=0$
and, for $c=-2f$, we obtain ${\rm sign}({\cal L}_3)=-{\rm sign}(
f)$.

Finally, for  system (\ref{sys_cs4b}) with ${\cal L}_1= {\cal
L}_2=0$ we calculate $ {\cal N}_1= -2^33^5 x(x^2+y^2) (4bx^2-
7axy- by^2)$, and hence, conditions $a=b=0$ are equivalent to
${\cal N}_1=0$. Moreover, the $GL$--comitant ${\cal N}_1$  is a
$CT$--comitant for  system (\ref{sys_cs4b}) when conditions ${\cal
L}_1= {\cal L}_2=0$ are satisfied. Thus, the following result
holds.

\bpr\label{prop:11-13} A cubic system (\ref{sys_cs4b}) possesses
invariant straight lines with total multiplicity $9$ if and only
if ${\cal L}_1= {\cal L}_2= {\cal N}_1=0.$ Moreover, the
configuration or the potential configuration of the lines
corresponds with $(3,2,2,1)$ given in Figure $11$ (respectively,
$12$; $13$) for ${\cal L}_3$ negative (respectively, positive;
zero). \epr

From Lemma \ref{lm:5} and Propositions \ref{prop:7-8},
\ref{prop:9-10} and \ref{prop:11-13} the next theorem follows.

\bth\label{th:7-13} We assume that for a cubic system (\ref{s_1})
the condition ${\cal D}_1<0 $  holds (i.e. there are 2  real and 2
imaginary infinite singular points). Then, this system  possesses
the maximum number of invariant straight lines (with total
multiplicity 9) if and only if one of the following sequences of
conditions holds:
$$
\ba{ll}
 {\cal V}_{1}={\cal V}_{2} = {\cal L}_{1} = {\cal L}_{2}= {\cal
 N}_{1}=0; &\hskip2cm ({\cal B}_1)\\
 {\cal V}_{3}={\cal V}_{4} = {\cal L}_{1} = {\cal L}_{2}= {\cal
 N}_{1}=0. &\hskip2cm ({\cal B}_2)\\
\ea
$$
Moreover, the configuration of the lines corresponds to Figures 7,
8, 9 or 10 for $({\cal B}_1)$ if $|{\cal L}_3|{\cal L}_4<0$,\
${\cal L}_3=0$ and ${\cal L}_4<0$,\ $|{\cal L}_3|{\cal L}_4>0$,\
or ${\cal L}_3=0$ and ${\cal L}_4>0$, respectively;\ to Figures
11, 12, or 13 for $({\cal B}_2)$ if ${\cal L}_3$ negative,
positive or zero, respectively.
\eth

From Theorems  \ref{th:1-6} and \ref{th:7-13}, and the fact that
the $GL$--comitants ${\cal L}_1$ and $ {\cal L}_2$ are
$T$--comitants for the initial cubic system (\ref{s_1}) we obtain
the next result.

\blm\label{lm:L12} Conditions $ {\cal L}_1= {\cal L}_2=0 $ are
necessary in order that a cubic system possesses the maximum
number of the invariant straight lines counted with their
multiplicities. \elm

\BProof The lemma is obvious for cubic systems having a {\it
generic behavior at infinity}, i.e. the multiplicity of all
singular points at infinity as roots of the polynomial $C_3$ is
one. Then, by continuity, it follows for the cubic systems having
non--generic behavior at infinity. \EProof

\bigskip

\noindent\textbf{ 7. Cubic systems with 1 triple and 1 simple real
roots of $C_3$}

By (\ref{HSys:6}) a cubic system having one triple and one real
distinct infinite singular points via a linear transformation
becomes the system:
\be\label{S6:1} \ba{l}
    x'=p_0 +p_1(x,y)+p_2(x,y) +3(u+1)x^3+3vx^2y+3rxy^2, \\
    y'=q_0 +q_1(x,y)+q_2(x,y) +3ux^2y+3vxy^2+3ry^3. \\
\ea \ee For  system (\ref{S6:1}) we have $ C_3=3x^3y$. Hence, the
infinite singular points are situated at the ``ends'' of the
following straight lines: $x=0$ and $ y=0$.

The aim of this section is to construct cubic systems of the form
(\ref{S6:1}) which possess invariant straight lines with total
multiplicity $9$, having $8$ affine lines with
potential confi\-gu\-ration $(3,\,3,\,1,\,1)$ or $(3,\,2,\,2,\,1)$.

\bigskip

\noindent\textbf{ 7.1 Systems with the potential configuration
$(3,\,3,\,1,\,1)$}

By Lemma \ref{lm:1}, if a cubic system possesses 8 invariant
straight lines with configuration or potential configuration
$(3,\,3,\,1,\,1)$, then the conditions ${\cal V}_1= {\cal V}_2=
{\cal U}_{\,1}=0$ hold.

A straightforward computation of the value of ${\cal V}_1$ for
system (\ref{S6:1}) provides that $\displaystyle {\cal V}_1$ is
equal to $\displaystyle 16\sum_{j=0}^4 {\cal V}_{1j}x^{4-j}y^j$,
where
$$
\ba{ll}
{\cal V}_{10}=u(2u+3),\ & {\cal V}_{12}=4ru +3r+2v^2,\\
{\cal V}_{11}=v(4u+3),\ &  {\cal V}_{13}=4vr,\ \ {\cal
V}_{14}=2r^2.
\ea
$$
Consequently, relations ${\cal V}_1=0$ imply $v=r=0$ and
$u(2u+3)=0$, and we have to consider two subcases $u=0$ and
$u=-3/2$.

For $u=0$, by the time rescaling $t\to t/3$, we obtain the  system
\be\label{SC6}
    x'=p_0 +p_1+p_2 +x^3, \quad
    y'=q_0 +q_1+q_2, \\
\ee whereas for $u=-3/2$, after the time rescaling $t\to -2t/3$,
we have the system \be\label{SC7}
    x'=p_0 +p_1+p_2 +x^3, \quad
    y'=q_0 +q_1+q_2 +3x^2y. \\
\ee It has to be underlined that for systems (\ref{SC6}) and
(\ref{SC7}) the relation ${\cal V}_2(\textbf{ a})={\cal U}_{\,1}(\textbf{
a})=0$ holds. On the other hand, by calculating the value of the
$T$--comitant ${\cal L}_4$ for system (\ref{SC6}) (respectively,
(\ref{SC7})), we obtain $ {\cal L}_4=-6x^4$ (respectively, ${\cal
L}_4=24 x^4)$. Hence, we get the next result.

\bpr\label{prop:L5} If for a cubic system the conditions $ {\cal
D}_1={\cal D}_3={\cal D}_4=0,\, {\cal D}_2\not=0$ and ${\cal
V}_1=0$ hold, then via a linear transformation and time rescaling  this  system
becomes into the form  (\ref{SC6}) for ${\cal L}_4<0$,  and into
the system (\ref{SC7}) for ${\cal L}_4>0$. \epr

\brm\label{rem:SC6} We note that for  system (\ref{SC6}) the
relations ${\cal V}_3= {\cal V}_4= {\cal U}_{\,2}=0$ are also
satisfied. By Lemma \ref{lm:4} this system can
present the potential configuration $(3,\,2,\,2,\,1).$ \erm

As it was proved above, we must examine systems (\ref{SC6}) and
(\ref{SC7}).

\noindent\textbf{ System (\ref{SC6})} Via the translation of the
origin of coordinates at the point $(-g/3,0)$, we can consider
$g=0$ and, hence, we get the system
\be\label{SC6:1} \dot x= a + c
x + d y + 2hxy+ky^2+x^3,\quad
\dot y= b + e x + f y + lx^2+2mxy+ ny^2.\\
\ee
For this system we have $C_3(x,y)=x^3y$ and therefore, there exist
two directions for the possible invariant straight lines: $x=0$
and $y=0$.

\noindent\textbf{ Direction $x=0$.} In this case  $U=1,\ V=0$ and,
from  (\ref{eq_g}), we obtain \beq
&&A=1,\ B=C=0,\ D=-W,\ E=2h,\ F=W^2+c,\\
&& Eq_7=k,\ Eq_9=-2hW+d,\ Eq_{10}=-W^3-cW+a.
\eeq
So, in order to have the maximum number of invariant straight
lines we obtain the conditions: $k=h=d=0$.

\noindent\textbf{ Direction $y=0$.} In this case  $U=0,\ V=1$ and,
from  (\ref{eq_g}), we obtain \beq
&&A=B=C=0,\ D=2m,\ E=n,\ F=-nW+f,\\
&& Eq_5=l,\ Eq_8=-2mW+e,\ Eq_{10}=nW^2-fW+b. \eeq So, we get the
conditions : $l=m=e=0$. Thus,  system (\ref{SC6:1}) becomes
$$
\dot x= a + c x  +x^3,\quad
\dot y= b +  f y + ny^2,\\
$$
for which we calculate ${\cal L}_1=0,\quad {\cal L}_2= -2^8 3^5
n^2 xy.$ Hence, since in order to reach the total multiplicity
$9$, by Lemma \ref{lm:L12}, it is necessary ${\cal L}_2=0$, we
obtain the additional condition: $n=0$. Therefore, system
(\ref{SC6}) goes over to \be\label{SC6:2} \dot x= a + c x
+x^3,\quad
\dot y= b +  f y. \\
\ee By Lemmas \ref{lm2} and \ref{lm3}, in order to determine the
possible invariant straight lines we shall use the affine
comitants ${\cal G}_i$ $(i=1,2,3)$. We consider the homogenized
system
\be\label{SC6:2H} \dot X= aZ^3 + c XZ^2 + X^3,\quad
\dot Y= bZ^3 +  f YZ^2, \\
\ee associated to system (\ref{SC6:2}) and calculate the following
polynomial:
$$
  H(\ab,X,Y,Z)= \gcd({\cal G}_1,{\cal G}_2,{\cal G}_3)=
  2Z^2(fY+Zb)(aZ^3+cXZ^2+X^3).
$$
Therefore, by Lemmas \ref{lm2} and \ref{lm3}, we obtain that
system (\ref{SC6:2H}) has $7$  invariant affine straight lines
(counted with their multiplicities) $Z=0$, $fY+Zb=0$, $aZ^3+
cXZ^2+ X^3= 0.$ We observe that $Z=0$ has multiplicity $3$. So,
for having the total multiplicity $9$, by Lemmas \ref{lm2} and
\ref{lm3}, the  polynomial $H(\ab,X,Y,Z)$ must have degree $8$. In
order to find out the conditions to reach this situation we shall
calculate for system (\ref{SC6:2H}) the polynomials:
\be\label{S6}
\ba{ll}
{\cal G}_1/H=&3X^2+(c-f)Z^2=   T(X,Z),\\
{\cal G}_3/H=& 12(aZ^3+cXZ^2+X^3)^2= 12 S^2(X,Z).
\ea
\ee
From Lemma \ref{lm4}, in order to have the maximum number of
invariant  straight lines  it is necessary that $T(X,Z)\mid
S^2(X,Z)$. We  consider two subcases: $c-f\ne0$ and $c-f=0$.

\noindent\textbf{ Case $c-f\ne0$.} In this case polynomial $T(X,Z)$
has two distinct factors and hence, it is necessary $T(X,Z)\mid
S(X,Z)$. By Lemma \ref{Trudi:2}, we  have
$$
 R^{(1)}_{X}\left(T,S\right)= 3(2c+f)Z^2=0,\quad
 R^{(0)}_{X}\left(T,\,S\right)=-\left[27a^2+((2c+f)^2(c-f)\right]=0,
$$
and these conditions yield $2c+f=a=0$. Thus, $f=-2c\ne0$ (by
condition $c-f\ne0$), and then we can consider  $b=0$ via a
translation of the origin of coordinates at the point
$(0,b/(2c))$. As a result we obtain the  system
\be\label{sys_conf:11-12}
\dot x=  c x  +x^3,\quad
\dot y=  -2c y, \\
\ee
which possesses the invariant  straight lines: $ x=0,\ y=0,\
x=\pm\sqrt{-c}. $
 It is clear, that  $x=\pm\sqrt{-c}$
 are real for $c<0$, and imaginary for $c>0$.
Hence, we obtain Figure $14$ for $c<0$ and Figure $15$ for $c>0$.

By Remark \ref{rm:transf} we can consider $c\in \{-1,1\}$. Since for system
(\ref{sys_conf:11-12}) the conditions ${\cal V}_1= {\cal V}_2=0$
hold as well as the conditions ${\cal V}_3= {\cal V}_4=0$, we
conclude that this system has invariant straight lines with total
multiplicity $9$ and presents both potential configurations:
$(3,3,1,1)$ and $(3,2,2,1)$. This is proved by the four perturbed
systems constructed below.

\noindent\textbf{ Subcase $c=-1$}. The system \be\label{s:pert1} \dot
x=  -x +x^3,\quad \dot y= 2y+ 3\varepsilon y^2+ \varepsilon^2 y^3,
\ee with the invariant straight lines  $x=0$, $y=0$, $x=\pm 1$,
$\varepsilon y+1=0$, $\varepsilon y+2=0$, $\varepsilon y \pm
x+1=0$, has the configuration $(3,3,1,1)$; and the system
\be\label{s:pert2} \dot x= -x +x^3-3\varepsilon^2xy^2,\quad
\dot y= 2y - 2\varepsilon^2 y^3,\\
\ee with the invariant straight lines $x=0$, $y=0$, $\varepsilon
y=\pm1$, $x+\varepsilon y=\pm1$, $x-\varepsilon y =\pm1$, has the
configuration $(3,2,2,1)$.

\noindent\textbf{ Subcase $c=1$}. The system \be\label{s:pert3} \dot
x= x +x^3,\quad
\dot y= -2y+ 3\varepsilon y^2- \varepsilon^2 y^3,\\
\ee
with the invariant straight lines $x=0$, $y=0$, $\varepsilon
y=1$, $\varepsilon y=2$, $\varepsilon y \pm ix=1$, $x=\pm i$, has the
configuration $(3,3,1,1)$; and the system \be\label{s:pert4} \dot
x= x +x^3-3\varepsilon^2xy^2,\quad
\dot y= -2y - 2\varepsilon^2 y^3,\\
\ee with the invariant straight lines $x=0$, $y=0$, $\varepsilon
y=\pm1$, $\varepsilon y+ix=\pm i$, $\varepsilon y-ix =\pm i$ has the
configuration $(3,2,2,1)$.

\noindent\textbf{ Case $c-f=0$.} From (\ref{S6}) it follows that
polynomial $R(X,Z)$ has the double root $X=0$ and in virtue of
condition $R(X,Z)\mid S^2(X,Z)$ it is necessary $a=0$. Thus, we
obtain the  system \be\label{sys_conf:8-10} \dot x=  c x
+x^3,\quad
\dot y=  b+c y, \\
\ee and by Remark \ref{rm:transf} we can consider $c\in \{-1,0,1\}$. Moreover,
if $c\ne0$ then the translation of the origin of coordinates at
the point $(0,-b/c)$  implies $b=0$;\ and for $c=0$, by replacing
$y\to by$, we can set $b=1$. So, we obtain Figure $16$
(respectively, Figure $17$; $18$) for $c=-1$, $b=0$ (respectively
for $c=1$, $b=0$; $c=0$, $b=1$).

Since for system (\ref{sys_conf:8-10}) the conditions ${\cal V}_1=
{\cal V}_2=0$ hold as well as the conditions ${\cal V}_3={\cal
V}_4=0$, we conclude that this system has invariant straight lines
with total multiplicity 9 and presents both potential
configurations: $(3,3,1,1)$ and $(3,2,2,1)$. It remains to
construct the six perturbed systems which will prove the
realization of all possibilities.

\noindent\textbf{ Subcase $c=-1$}. The system
\be\label{s:pert5} \dot
x=  -x +x^3,\quad \dot y= -y+ \varepsilon^2 y^3, \ee with the
invariant straight lines $x=0$, $y=0$, $x=\pm1$, $\varepsilon
y=\pm1$, $x-\varepsilon y=0$, $x+\varepsilon y=0$, has the
configuration $(3,3,1,1)$; and the system \be\label{s:pert6} \dot
x= (4\varepsilon^2-1)x-6\varepsilon x
y+x^3-12\varepsilon^2xy^2,\quad \dot y=
-2\varepsilon-(8\varepsilon^2+1)y-6\varepsilon
y^2-8\varepsilon^2y^3, \ee with the invariant straight lines
$x=0$, $4\varepsilon y=-1$, $2x+4\varepsilon y+1= \pm \sqrt{1-
16\varepsilon^2}$, $2\varepsilon y^2+y=-2\varepsilon$, $2x-
4\varepsilon y -1= \pm \sqrt{1-16\varepsilon^2}$, has the
configuration $(3,2,2,1)$.

\noindent\textbf{ Subcase $c=1$}. The system
\be\label{s:pert7}
\dot x= x +x^3,\quad \dot y= y+ \varepsilon^2 y^3,
\ee
with the invariant straight lines $x=0$, $y=0$, $x=\pm i$, $\varepsilon
y=\pm i$, $x-\varepsilon y=0$, $x+\varepsilon y=0$, has the
configuration $(3,3,1,1)$; and the system
\be\label{s:pert8} \dot
x=x+6\varepsilon xy+x^3+12\varepsilon^2xy^2,\quad
\dot y= y+6\varepsilon y^2+8\varepsilon^2 y^3,\\
\ee with the invariant straight lines $x=0$, $y=0$, $4\varepsilon
y=-1$, $2\varepsilon y= -1$, $2\varepsilon y\pm ix= 0$, $2
\varepsilon y\pm ix= -1$, has the configuration $(3,2,2,1)$.

\noindent\textbf{ Subcase $c=0$}. The system \be\label{s:pert9} \dot
x= 18\varepsilon^2 x  + 9\varepsilon x^2+x^3,\quad \dot y=
1-6\varepsilon^2 y-24\varepsilon^4y^2 +64\varepsilon^6y^3, \ee
with the invariant straight lines $x=0$, $x=-3\varepsilon$, $x= -6
\varepsilon$, $x+ 8\varepsilon^3 y= -2\varepsilon$, $8
\varepsilon^2 y=1$, $4 \varepsilon^2 y= -1$, $2 \varepsilon^2 y=
1$, $x- 8 \varepsilon^3 y = -4 \varepsilon$, has the configuration
$(3,3,1,1)$; and the system \be\label{s:pert10} \dot
x=\varepsilon(1-4\varepsilon) x  +
12\varepsilon^3xy+x^3-12\varepsilon^4 xy^2,\quad
\dot y= 1-2\varepsilon(1+2\varepsilon)y+12\varepsilon^3y^2 -8\varepsilon^4y^3,\\
\ee with the invariant straight lines  $x=0$, $2 \varepsilon y=
1$, $x- 2 \varepsilon^2 y + \varepsilon = \pm \sqrt{
\varepsilon^2- \varepsilon}$, $\varepsilon(2 \varepsilon y-1)^2=
\varepsilon- 1$, $x+ 2\varepsilon^2 y -\varepsilon = \pm \sqrt{
\varepsilon^2- \varepsilon}$ has in configuration $(3,2,2,1)$.

We construct necessary and sufficient affine invariant conditions
for the realization of each possible  configuration given in
Figures $14$--$18$ for  system (\ref{SC6:1}).

As it was proved above system (\ref{SC6:1})  possesses 9 invariant
straight lines (counted with their multiplicities) if and only if
$h= k= l= m= n= 0$, $d= e= (c-f)(2c+f)= 0$, $a=0.$ First, we shall
show that conditions $h= k= l= m= n= 0$ are equivalent to
conditions ${\cal N}_2= {\cal N}_3=0$. To get this goal we
consider the system:
\be\label{SC6:1t1} \dot x_1= p^\tau_0 +
p^\tau_1 +3\gamma x_1^2+2hx_1y_1+ky_1^2+x_1^3,\quad
\dot y_1= q^\tau_0 + q^\tau_1 +lx_1^2+ 2mx_1y_1+ny_1^2,\\
\ee which is obtained from system (\ref{SC6:1}) via the
translation $x=x_1+\gamma $, $y=y_1+\delta  $, where $(\gamma
,\delta  )$ is an arbitrary point of the phase plane of system
(\ref{SC6:1}) and $p^\tau_i(a(\gamma ,\delta  ),x_1,y_1)$,
$q^\tau_i(a(\gamma ,\delta),x_1,y_1)$ $(i=0,1)$ are the
corresponding homogeneous polynomials of degree $i$ in $x_1$ and
$y_1$.

For system (\ref{SC6:1t1}) we  calculate the values of the
$GL$--comitants  ${\cal N}_2$ and ${\cal N}_3$ which depend only
on the coefficients of the quadratic and cubic parts of the
system: \beq
 && {\cal N}_2 = -6x_1^2\left[lx_1^3+3mx_1^2y_1+(2n-h)x_1y_1^2-ky_1^3\right], \\
 && {\cal N}_3 =-12 x_1^2\left[lx_1^3+(2h-n)x_1y_1^2+2ky_1^3\right].
\eeq
As we can observe, the values of these polynomials do not depend
on the coordinates of the arbitrary point $(\gamma ,\delta  )$
and, consequently conditions ${\cal N}_2=0$ and ${\cal N}_3=0$ are
affine invariant ones. It is obvious to find out that these
conditions yield $h=k=l=m=n=0$, i.e. all quadratic coefficients
vanish. Thus, system (\ref{SC6:1t1}) becomes the system
$$
\ba{ll}
\dot x_1=& (a+c\gamma +d\delta  +\gamma ^3) + (c+3\gamma ^2)x_1 + dy_1 +3\gamma x_1^2+x_1^3,\\
\dot y_1=& (b +e\gamma +f\delta  ) +ex_1 + f y_1,\\
 \ea
$$
for which we have:
\beq
 {\cal N}_6 &=& -36x_1^4\left[6ex_1^2+(c-f)x_1y_1-16dy_1^2\right],\\
 {\cal N}_7 &=& -72x_1^4\left[5ex_1^2+(2c+f)x_1y_1-18dy_1^2\right],\\
 {\cal N}_8 &=& -12(c+f)x_1^2.\\
 {\cal N}_1 &=& -2^33^5x_1^3y_1\left[ax_1+d(\delta  x_1-3\gamma y_1)\right].
\eeq As above condition ${\cal N}_6=0$ (respectively, ${\cal
N}_7=0$) is an affine invariant one and implies $d=e=c-f=0$
(respectively, $d=e=2c+f=0$). Moreover, in both  cases the
$GL$--comitant ${\cal N}_1$ will not depend on the coordinates of
the arbitrary point $(\gamma,\delta)$; i.e. it becomes a
$T$--comitant and the condition $N_1=0$ yields $a=0$. On the other
hand, for $f=c$ (respectively, for $f=-2c$ ) we obtain that ${\rm
sign} ({\cal N}_8)= -{\rm sign}(c)$ (respectively, ${\rm sign}
({\cal N}_8)= {\rm sign}(c)$). Thus, we obtain the next result.

\blm\label{lm:S6} System (\ref{SC6}) possesses invariant straight
lines with total multiplicity $9$ if and only if
$$
\hphantom{m}\hskip4cm {\cal N}_1= {\cal N}_2={\cal N}_3={\cal
N}_6{\cal N}_7= 0.
 \hskip5cm ({\cal C})
$$
Moreover, the configurations of the lines correspond to Figures
$14$ or $15$ for $({\cal C})$ with ${\cal N}_6\ne0$ and ${\cal
N}_8<0$ or ${\cal N}_8>0$, respectively; to Figures $16$, $17$ or
$18$ for $({\cal C})$ with ${\cal N}_6=0$ and ${\cal N}_8$
positive, negative or zero, respectively.
\elm

\noindent\textbf{ System (\ref{SC7})} It is clear that via a
translation of the origin of coordinates at the point
$(-g/3,-l/3)$, we can consider $g=l=0$ and, hence, we must examine
the  system \be\label{SC7:1} \dot x= a + c x + d y +
2hxy+ky^2+x^3,\quad
\dot y= b + e x + f y + 2mxy+ ny^2 +3x^2y.\\
\ee
For this system we have $C_3(x,y)=-2x^3y$ and therefore, there
exists two directions for the possible invariant straight lines:
$x=0$ and $y=0$.

\noindent\textbf{ Direction $x=0$.} In this case  $U=1,\ V=0$ and from
(\ref{eq_g}) we obtain \beq
&&A=1,\ B=C=0,\ D=-W,\ E=2h,\ F=W^2+c,  \\
&& Eq_7=k,\ Eq_9=-2hW+d,\ Eq_{10}=-W^3-cW+a.
\eeq
So, in order to have the maximum number of invariant straight
lines we obtain the conditions: $k=h=d=0$.

\noindent\textbf{ Direction $y=0$.} In this case  $U=0,\ V=1$ and from
(\ref{eq_g}) we obtain \beq
&&A=3,\ B=C=0,\ D=2m,\ E=n,\ F=-nW+f,\\
&& Eq_5=3W,\ Eq_8=-2mW+e,\ Eq_{10}=nW^2-fW+b. \eeq So, in this
direction can be only one simple invariant straight line with
$W_0=0$, and the necessary conditions are $e=b=0$.

In short, system (\ref{SC7:1}) becomes
$$
\dot x= a + c x +x^3,\quad
\dot y= f y + 2mxy+ ny^2 +3x^2y,\\
$$
for which we calculate ${\cal L}_1= 2^8 3^4 n x^3$, ${\cal L}_2=
2^8 3^4 n x (mx+9ny)$. Hence, in order to reach the total
multiplicity $9$, by Lemma \ref{lm:L12}, it is necessary ${\cal
L}_1= {\cal L}_2= 0$, and we obtain the additional condition:
$n=0$. This provides the system:
\be\label{SC7:2} \dot x= a + c x
+x^3,\quad
\dot y= f y + 2mxy +3x^2y.\\
\ee By Lemmas \ref{lm2} and \ref{lm3}, in order to determine the
possible invariant straight lines we shall use the affine
comitants ${\cal G}_i$ $(i=1,2,3)$. We consider the homogenized
system \be\label{SC7:2H} \dot X= aZ^3 + c XZ^2 + X^3,\quad
\dot Y= f YZ^2+2mXYZ+3X^2Y,\\
\ee corresponding to system (\ref{SC7:2}), and calculate the
following polynomial:
$$
  H(\ab,X,Y,Z)= \gcd({\cal G}_1,{\cal G}_2,{\cal G}_3)=
  2Y(aZ^3+cXZ^2+X^3),
$$
and
\be\label{RMS}
\ba{ll}
{\cal G}_1/H=& -6X^4-8mX^3Z- (3f+4m^2+3c)X^2Z^2\\
&-2(3a+2fm)XZ^3-(2am+f^2-cf)Z^4 =  T(X,Z),\\
{\cal G}_2/H=&
-3(3aZ^3+2mfZ^3+4Z^2m^2X+3Z^2Xf+3cXZ^2\\
  &+12mX^2Z+12X^3)(aZ^3+cXZ^2+X^3) =  S(X,Z).
\ea
\ee

By Lemma \ref{lm4}, in order to have the maximum number of
invariant straight lines it is necessary that $R(X,Z)\mid S(X,Z)$.
Hence, by Lemma \ref{Trudi:2}, the  conditions $ R^{(i)}_{X}\left(
T, S\right)=0\ (i=0,1,2,3), $ have to be satisfied, where
$R^{(0)}_{X}\left( T, S\right)= {\rm Res}_{X}\left( T, S\right)$.
So, we can calculate: $ R^{(3)}_{X}\left( T, S\right)=
72(27a+18cm+8m^3)Z^3= 0, $ and, hence, we obtain the condition:
$a= -(18cm+8m^3)/27$. Then, we have
\beq R^{(2)}_{X}\left( R,
S\right)&=&-36(9c-3f+4m^2)^2 (3c-3f+20m^2)^2 Z^8= 0, \eeq and
this implies the necessity to examine the two subcases:\
$9c-3f+4m^2=0$ and $9c-3f+4m^2\ne0$.

\noindent\textbf{ Case $9c-3f+4m^2=0$.} Then, $f=(9c+4m^2)/3$, and we
have $ R^{(1)}_{X}\left(T, S\right)=R^{(0)}_{X}\left(T,
S\right)=0, $ but at the same time this provides the degenerated
system:
$$
\dot x= -(2m-3x)(4m^2+6xm+9c+9x^2)/27,\quad
\dot y= y(4m^2+6xm+9c+9x^2)/3. \\
$$

\noindent\textbf{ Case $9c-3f+4m^2\ne0$.} Then, condition $
R^{(2)}_{X}\left(T, S\right)=0$ implies $f=(20m^2+3c)/3$, and we
can calculate \ $ R^{(1)}_{X}\left(T, S\right)= 2^{11} 3^2 m^3
(3c+28m^2)^3 (3c-8m^2)^3 Z^{15}= 0, $ and since condition
$9c-3f+4m^2\ne0$ yields $3c- 8m^2\ne 0$, we obtain condition
$m(3c+28m^2)=0$.

If $m\ne0$ then $c=-28m^2/3\ne0$, and this provides the
degenerated system:
$$
\dot x= (3x+10m)(-3x+8m)(-3x+2m)/27,\quad
\dot y= -y(4m+3x)(-3x+2m)/3. \\
$$
Thus, $m=0$ and we get the system :
\be\label{sys_conf:14-15}
\dot x=  c x  +x^3,\quad
\dot y=  c y+3x^2y, \\
\ee with $c\ne0$ (otherwise system becomes degenerated). We note
that by Remark \ref{rm:transf}  we can consider $c\in \{-1,1\}$. So, we obtain
Figure 19 for $c=-1$, and Figure 20 for $c=1$. It remains to
construct two perturbed systems.

For $c=-1$ the perturbed  system
\be\label{s:pert11} \dot x=  -x
+x^3,\quad
\dot y= -y+ 3x^2y-18\varepsilon xy^2+ 36\varepsilon^2 y^3,\\
\ee possesses the invariant straight lines $x=0$, $x=\pm1$, $y=0$,
$x- 3\varepsilon y= 0$, $x- 6\varepsilon y= 0$, $x- 6\varepsilon y
= \pm1$ in the configuration $(3,3,1,1)$; and for $c=1$ the system
\be\label{s:pert12} \dot x=x +x^3,\quad
\dot y= y+ 3x^2y-18\varepsilon xy^2+ 36\varepsilon^2 y^3,\\
\ee possesses the invariant straight lines $x=0$, $x=\pm i$,
$y=0$, $x- 3\varepsilon y= 0$, $x- 6\varepsilon y= 0$, $x- 6
\varepsilon y= \pm i$ in the same configuration $(3,3,1,1)$.

We construct necessary and sufficient affine invariant conditions
for the realization of the two possible configurations for system
(\ref{SC7:1}), see Figures $19$ and $20$.

As it was proved above system (\ref{SC7:1}) possesses $9$
invariant straight lines (counted with their multiplicities) if
and only if $h=k=m=n=0$, $d=e=c-f=0$, $a=b=0$. First, we shall
show that conditions $h=k=m=n=0$ are equivalent to conditions
${\cal N}_4= {\cal N}_5=0$. To get this goal  we shall consider
the system:
\be\label{SC7:1t1} \ba{ll}
\dot x_1=& p^\tau_0 + p^\tau_1 +3\gamma x_1^2+2hx_1y_1+ky_1^2+x_1^3,\\
\dot y_1=& q^\tau_0 + q^\tau_1 +3\delta  x_1^2+ 2(m+3\gamma )x_1y_1+ny_1^2 +3x_1^1y_1,\\
 \ea
\ee which is obtained from system (\ref{SC7:1}) via the
translation $x=x_1+\gamma$, $y=y_1+\delta$, where $(\gamma ,\delta
)$ is an arbitrary point of the phase plane of system
(\ref{SC7:1}) and $p^\tau_i(a(\gamma ,\delta  ),x_1,y_1)$,
$q^\tau_i(a(\gamma ,\delta  ),x_1,y_1)$ $(i=0,1)$ are the
corresponding homogeneous polynomials of degree $i$ in $x_1$ and
$y_1$.

For system (\ref{SC7:1t1}) we calculate the values of the
$GL$--comitants  ${\cal N}_4$ and ${\cal N}_5$ which depend only
on the coefficients of the quadratic and cubic parts of the
system:
$$
 {\cal N}_4 = 24x_1\left[3mx_1^2+5(h+n)x_1y_1+6ky_1^2\right], \quad
 {\cal N}_5 =-48 x_1y_1\left[(2h-n)x_1+3ky_1\right],
$$
The values of these polynomials do not depend on the coordinates
of the arbitrary point $(\gamma ,\delta  )$ and, consequently
conditions ${\cal N}_4=0$ and ${\cal N}_5=0$ are affine invariant.
These conditions yield $h=k=m=n=0$, i.e. all quadratic
coefficients vanish. Thus, system (\ref{SC7:1t1}) becomes
$$
\ba{ll}
\dot x_1=& (a+c\gamma +d\delta  +\gamma ^3) + (c+3\gamma ^2)x_1 +
dy_1 +3\gamma x_1^2+x_1^3,\\
\dot y_1=& (b +e\gamma +f\delta  +3\gamma ^2\delta  ) +(e +6\gamma
\delta) x_1 + (f+3\gamma ^2)y_1+ 3\delta  x_1^2+6\gamma x_1y_1+3x_1^2y_1,\\
 \ea
$$
for which we have:
\beq
 {\cal N}_6 &=& -144x_1^4\left[3ex_1^2+2(c-f)x_1y_1-17dy_1^2\right],\\
 {\cal N}_8 &=& -24(c+f)x_1^2,\\
 {\cal N}_1 &=& -216x_1^3\left[-18x_1(bx_1+ay_1)
 +(c-f)(23\delta  x_1+25\gamma y_1)x_1\right.\\
&&\left. -12e\gamma x_1^2+2dy_1(19\gamma y_1-13\delta x_1)\right].
\eeq
As above ${\cal N}_6=0$ is an affine invariant condition and
implies $d=e=c-f=0$. Moreover, for ${\cal N}_6=0$ we obtain that
the $GL$--comitant ${\cal N}_1$ does not depend on the coordinates
of the arbitrary point $(\gamma ,\delta  )$, i.e. it becomes a
$T$--comitant and condition $N_1=0$ yields $a=b=0$. On the other
hand for $f=c$  we obtain that ${\rm sign} ({\cal N}_8)=-{\rm
sign}(c)$. Thus, the next lemma follows.

\blm\label{lm:S7} Canonical system (\ref{SC7}) possesses invariant
straight lines with total multiplicity $9$ if and only if ${\cal
N}_1= {\cal N}_4= {\cal N}_5= {\cal N}_6= 0.$ Moreover, the potential
configuration of the lines corresponds with $(3,3,1,1)$ given in
Figure $19$ for ${\cal N}_8>0$ and in Figure $20$ for ${\cal
N}_8<0$.
\elm

\bigskip

\noindent\textbf{ 7.2  Systems  with  the potential  configuration
$(3,\,2,\,2,\,1)$}

By Lemma \ref{lm:4}
 if a cubic system possesses 8
invariant straight lines in the potential configuration $(3,\,2,\,2,\,1)$,
then  the conditions ${\cal V}_3= {\cal V}_4= {\cal U}_{\,2}=0$ hold.

We consider again (\ref{S6:1}). A straightforward computation of
the value of ${\cal V}_3$ yields:\linebreak $\displaystyle {\cal V}_3= 2^5
3^2 \sum_{j=0}^4 {\cal V}_{3j} x^{4-j} y^j$, where
$$\ba{ll}
{\cal V}_{30}=-u(u+3),& {\cal V}_{32}=-2ru+3r-v^2,\\
{\cal V}_{31}=-2uv, & {\cal V}_{33}=-2vr,\quad
{\cal V}_{34}=-r^2.
\ea
$$
Hence $r=v=0$ and $u(u+3)=0$. So, we have to examine two cases:
$u=-3$ and $u=0$.

For $u=-3$, by time rescaling $t\to -t/3$, we obtain the system
\be\label{SC8}
    x'=p_0 +p_1(x,y)+p_2(x,y) +2x^3, \quad
    y'=q_0 +q_1(x,y)+q_2(x,y) +3x^2y. \\
\ee Whereas for $u=0$, after the time rescaling $t\to t/3$, we get
system (\ref{SC6}) (see Remark \ref{rem:SC6}). It remains to note
that for system (\ref{SC8}) we have ${\cal V}_4= {\cal U}_{\,2}=
0$. So, for additional investigations we must consider only
system (\ref{SC8}).

It is clear that via the translation of the origin of coordinates
at the point $(-g/6,-l/3)$, we can consider $g=l=0$ and, hence, we
must examine the  system
\be\label{SC8:1}
\dot x= a + c x + d y + 2hxy+ky^2+2x^3,\quad
\dot y= b + e x + f y + 2mxy+ ny^2 +3x^2y.\\
\ee
For this system we have $C_3(x,y)=-x^3y$. Therefore, there exists
two directions for possible invariant straight lines: $x=0$ and
$y=0$.

\noindent\textbf{ Direction $x=0$.} In this case  $U=1,\ V=0$ and,
from  (\ref{eq_g}), we obtain \beq
&&A=2,\ B=C=0,\ D=-2W,\ E=2h,\ F=2W^2+c,   \\
&& Eq_7=k,\ Eq_9=-2hW+d,\ Eq_{10}=-2W^3-cW+a.
\eeq
So, in order to have the maximum number of invariant straight
lines we obtain the conditions: $k=h=d=0$.

\noindent\textbf{ Direction $y=0$.} In this case  $U=0,\ V=1$ and,
from  (\ref{eq_g}), we obtain \beq
&&A=3,\ B=C=0,\ D=2m,\ E=n,\ F=-nW+f,\\
&& Eq_5=-3W,\ Eq_8=-2mW+e,\ Eq_{10}=nW^2-fW+b.
\eeq
Hence, in this direction can be only one simple invariant straight
line with $W_0=0$, and the necessary conditions are $e=b=0$. Thus,
system (\ref{SC8:1}) becomes
$$
\dot x= a + c x +2x^3,\quad
\dot y= f y + 2mxy+ ny^2 +3x^2y,
$$
for which we calculate ${\cal L}_1= 2^93^4nx^3$, ${\cal L}_2= 2^9
3^4 n x (mx+3ny)$. Hence, in order to reach the total multiplicity
$9$, by Lemma \ref{lm:L12}, it is necessary  ${\cal L}_1={\cal
L}_2=0$. Then we obtain the additional condition: $n=0$. This
provides the system:
\be\label{SC8:2} \dot x= a + c x +2x^3,\quad
\dot y= f y + 2mxy +3x^2y. \ee We consider the homogenized system
$$
\dot X= aZ^3 + c XZ^2 + 2X^3,\quad
\dot Y= f YZ^2+2mXYZ+3X^2Y,\\
$$
corresponding to system (\ref{SC8:2}) and calculate the
polynomial:
$$
  H(\ab,X,Y,Z)= \gcd({\cal G}_1,{\cal G}_2,{\cal G}_3)=
  2Y(aZ^3+cXZ^2+2X^3),
$$
and
$$
\ba{ll}
{\cal G}_1/H=&
-3X^4-4mX^3Z-(3c+4m^2)X^2Z^2-2(3a+2fm)XZ^3\\
&-(2am-cf +f^2)Z^4 = \tilde T(X,Z),\\
{\cal G}_2/H=&
 -3(3aZ^3+2mfZ^3+3cXZ^2+4Z^2m^2X+6mX^2Z+6X^3)\cdot \\
 & \times (aZ^3+cXZ^2+2X^3) = \tilde S(X,Z).
\ea
$$
By Lemma \ref{lm4}, in order to have the maximum number of
invariant straight lines it is necessary that $\tilde T(X,Z)\mid
\tilde S(X,Z)$, i.e. the conditions $R^{(i)}_{X} \left(\tilde
T,\tilde S\right)= 0$, $(i=0,1,2,3), $ have to be satisfied. So,
we can calculate: \beq R^{(3)}_{X}\left(\tilde T,\tilde
S\right)&=&-18(54a-45cm+54fm-40m^3)Z^3=0, \eeq and obtain that $a=
(45cm-54fm+40m^3)/54$. Then, we have \beq R^{(2)}_{X}\left(\tilde
T,\tilde S\right)&=&-81(9c-6f+8m^2)^2 (c+2 f+3 m^2)^2 Z^8=0, \eeq
and this implies the necessity to examine two cases:\
$9c-6f+8m^2=0$ and $9c- 6f+ 8m^2\ne 0$.

\noindent\textbf{ Case $9c-6f+8m^2=0$.} Then $f=(9c+8m^2)/6$, and we
have $ R^{(1)}_{X}\left(\tilde T,\tilde S\right)= R^{(0)}_{X}
\left(\tilde T, \tilde S\right)= 0$, but at the same time this
provides the degenerated system:
$$
\dot x= -(2m-3x)(8m^2+12xm+9c+18x^2)/27,\quad
\dot y= y(8m^2+12xm+9c+18x^2)/6. \\
$$

\noindent\textbf{ Case $9c-6f+8m^2\ne0$.} Then condition $R^{(2)}_{X}
\left(\tilde T, \tilde S\right)= 0$ implies $f=-(c+3m^2)/2$, and
we have $ R^{(1)}_{X} \left(\tilde T,\tilde S\right)= -9m^3
(6c+13m^2)^3 (12c+17m^2)^3 Z^{15}=0, $ and since condition
$9c-6f+8m^2\ne0$ yields $12c+17m^2\ne0$, we obtain condition
$m(6c+13m^2)=0$.

If $m\ne0$ then $c=-13m^2/6\ne0$, and this provides the
degenerated system:
$$
\dot x= -(6x+5m)(-6x+7m)(3x+m)/54,\quad
\dot y= -y(6x+5m)(m-6x)/12. \\
$$
Thus, $m=0$ and we obtain the system (we set $c=2p$):
\be\label{sys_conf:22-23}
\dot x=  2p x  +2x^3,\quad
\dot y=  -p y+3x^2y, \\
\ee with $p\ne0$ (otherwise system becomes degenerated). We note
that by Remark \ref{rm:transf} we can consider $p=\in \{-1,1\}$. So, we obtain
Figure 21 for $p=-1$ and Figure 22 for $p=1$. It remains to
construct two perturbed systems.

For $p=-1$ the perturbed  system
\be\label{s:pert13} \dot x= -2x
+2x^3,\quad
\dot y= y+ 3x^2y- \varepsilon^2 y^3,\\
\ee possesses the invariant straight lines $x=0$, $y=0$, $x=\pm
1$, $x+ \varepsilon y= \pm 1$ and $x- \varepsilon y= \pm 1$ in the
configuration $(3,2,2,1)$; and for $p=1$ the system
\be\label{s:pert14} \dot x=2x +2x^3,\quad
\dot y= -y+ 3x^2y+ \varepsilon^2 y^3,\\
\ee possesses the invariant straight lines $x=0$, $y=0$, $x=\pm
i$, $\varepsilon y+ix= \pm 1$ and $\varepsilon y-ix= \pm 1$ in the
same configuration $(3,2,2,1)$.

We construct necessary and sufficient affine invariant conditions
for the realization of Figures $21$ and $22$ for system
(\ref{SC8:1}). As it was proved above, system (\ref{SC8:1})
possesses $9$ invariant straight lines (counted with their
multiplicities) if and only if $h=k=m=n=0$, $d=e=c+2f=0$, $a=b=0$.
First, we shall show that conditions $h=k=m=n=0$ are equivalent to
conditions  ${\cal N}_4= {\cal N}_5=0$. To get this goal we
consider the system: \be\label{SC8:1t1} \ba{ll}
\dot x_1=& p^\tau_0 + p^\tau_1 +6\gamma x_1^2+2hx_1y_1+ky_1^2+2x_1^3,\\
\dot y_1=& q^\tau_0 + q^\tau_1 +3\delta  x_1^2+ 2(m+3\gamma )x_1y_1+ny_1^2 +3x_1^1y_1,\\
 \ea
\ee which is obtained from system (\ref{SC8:1}) via the
translation $x=x_1+\gamma $, $y=y_1+\delta$, where $(\gamma
,\delta)$ is an arbitrary point of the phase plane of system
(\ref{SC8:1}), and $p^\tau_i(a(\gamma ,\delta  ),x_1,y_1)$,
$q^\tau_i(a(\gamma,\delta),x_1,y_1)$ $(i=0,1)$ are the
corresponding homogeneous polynomials of degree $i$ in $x_1$ and
$y_1$.

For system (\ref{SC8:1t1}) we  calculate the values of the
$GL$--comitants
$$
 {\cal N}_4 = 24x_1\left[2mx_1^2+3(h+n)x_1y_1+3ky_1^2\right], \quad
 {\cal N}_5 =24 x_1\left[2mx_1^2+3(n-2h)x_1y_1-6ky_1^2\right].
$$
As we can observe, the values of these polynomials do not depend
on the coordinates of the arbitrary point $(\gamma ,\delta  )$
and, consequently  ${\cal N}_4=0$ and ${\cal N}_5=0$ are affine
invariant conditions. These conditions yield $h=k=m=n=0$, i.e. all
quadratic coefficients vanish. Thus, system (\ref{SC8:1t1})
becomes
$$
\ba{ll}
\dot x_1=& (a+c\gamma +d\delta  +2\gamma ^3) + (c+6\gamma ^2)x_1 +
dy_1 +6\gamma x_1^2+2x_1^3,\\
\dot y_1=& (b +e\gamma +f\delta  +3\gamma ^2\delta  ) +(e +6\gamma
\delta  ) x_1 + (f+3\gamma ^2)y_1+ 3\delta  x_1^2+6\gamma x_1y_1+3x_1^2y_1,\\
 \ea
$$
for which we have:
\beq
 {\cal N}_7 &=& -216x_1^4\left[2ex_1^2-(c+2f)x_1y_1-5dy_1^2\right],\\
 {\cal N}_8 &=& -36(c+f)x_1^2,\\
 {\cal N}_1 &=& -648x_1^3\left[3x_1(4bx_1-7ay_1)
 +7(c+2f)(\delta  x_1-\gamma y_1)x_1\right.\\
&&\left. +6e\gamma x_1^2+22dy_1(\gamma y_1-\delta  x_1)\right].
\eeq As above condition ${\cal N}_7=0$ is affine invariant and
implies $d=e=c+2f=0$. Moreover, for ${\cal N}_7=0$ we obtain that
the $GL$--comitant ${\cal N}_1$ does not depend on the coordinates
of the arbitrary point $(\gamma ,\delta  )$, i.e. it becomes a
$T$--comitant and condition ${\cal N}_1=0$ yields $a=b=0$. On the
other hand, for $c=-2f$  we obtain that ${\rm sign} ({\cal
N}_8)={\rm sign}(f)$. Thus,  the following lemma holds.

\blm\label{lm:S8} Canonical system (\ref{SC8:1}) possesses
invariant straight lines with total multiplicity $9$ if and only
if ${\cal N}_1= {\cal N}_4= {\cal N}_5= {\cal N}_7= 0$. Moreover,
the potential configuration of the lines corresponds with $(3,2,2,1)$ given in
Figure $21$ for ${\cal N}_8>0$ and in Figure $22$ for ${\cal
N}_8<0$. \elm

By Lemmas \ref{lm:5}, \ref{lm:S6}, \ref{lm:S7} and \ref{lm:S8} the
next theorem follows.

\bth\label{th:14-22} We assume that for cubic system (\ref{s_1})
the conditions ${\cal D}_1= {\cal D}_3= {\cal D}_4=0$, ${\cal
D}_2\not=0$ hold (i.e. $C_3$ has $1$ triple and $1$ simple real
roots). Then, this system will possess the maximum number of
invariant straight lines (with total multiplicity $9$) if and only
if at least one of the following sets of conditions is fulfilled:
$$
\ba{ll}
\hskip2cm {\cal V}_{1}= {\cal N}_{1}= {\cal N}_{2}
 ={\cal N}_{3}= {\cal N}_{7}=0,\ {\cal L}_4<0; &\hskip3cm({\cal D}_1)\\
\hskip2cm {\cal V}_{1}= {\cal N}_{1}= {\cal N}_{2}
 ={\cal N}_{3}= {\cal N}_{6}=0,\ {\cal L}_4<0; &\hskip3cm({\cal D}_2)\\
\hskip2cm {\cal V}_{1}= {\cal N}_{1}= {\cal N}_{4}
 ={\cal N}_{5}= {\cal N}_{6}=0,\ {\cal L}_4>0; &\hskip3cm({\cal D}_3)\\
\hskip2cm {\cal V}_{3}= {\cal N}_{1}= {\cal N}_{4}
 ={\cal N}_{5}= {\cal N}_{7}=0. &\hskip3cm({\cal D}_4)\\
\ea
$$
Moreover, the configuration of the lines corresponds to Figures
$14$ or $15$ for $({\cal D}_1)$ if ${\cal N}_8$ negative or
positive, respectively; to Figures $16$, $17$ or $18$ for $({\cal
D}_2)$ if ${\cal N}_8$ positive, negative or zero, respectively;
to Figures $19$ or $20$ for $({\cal D}_3)$ if ${\cal N}_8$
positive or negative, respectively; to Figures $21$ or $22$ for
$({\cal D}_4)$ if ${\cal N}_8$ positive or negative, respectively.
\eth

\bigskip

\noindent\textbf{ 8. Cubic systems with $1$ real root of $C_3$ with
multiplicity $4$}

The objective of this section is to construct the cubic systems
with one real infinite singular point of multiplicity $4$ which
have invariant straight lines with total multiplicity $9$.

First we obtain  the   form for the homogeneous part of degree $3$
for these cubic systems.

\blm\label{lm:s9} Every cubic system with one real infinite
singular point of multiplicity $4$ which can admit invariant
straight lines with total multiplicity $9$ via a linear
transformation can be written as
\be\label{SC9}
    x'=p_0 +p_1(x,y)+p_2(x,y), \quad
    y'=q_0 +q_1(x,y)+q_2(x,y)-x^3. \\
\ee
\elm

\BProof As it  was shown in Section $3$ the cubic system having
one real infinite singular point of multiplicity $4$ via a linear
transformation becomes:
\be\label{S9} \ba{l}
    x'=p_0 +p_1(x,y)+p_2(x,y) +ux^3+vx^2y+rxy^2, \\
    y'=q_0 +q_1(x,y)+q_2(x,y) -x^3 +ux^2y+vxy^2+ry^3, \\
\ea \ee For system (\ref{S9}) we obtain $ C_3=x^4$ and hence, the
infinite singular point is situated on the ``end'' of the line
$x=0$.

By Lemmas \ref{lm:1} and \ref{lm:4}, we need to consider two
cases.

\noindent\textbf{ Case 1}: ${\cal V}_1= {\cal V}_2= {\cal U}_{\,1}=0$.
Then, a straightforward computation of the value of ${\cal
V}_1$ for system (\ref{S9}) yields: $\displaystyle {\cal
V}_1= 2^5 \sum_{j=0}^4{\cal V}_{1j}x^{4-j}y^j, $ where ${\cal
V}_{10}= u^2$, ${\cal V}_{11}= 2uv$, ${\cal V}_{12}= 2ur+ v^2$,
${\cal V}_{13}= 2rv$, ${\cal V}_{14}= r^2$. Consequently, relation
${\cal V}_1=0$ yield $v=r=u=0$ and hence, we obtain system
(\ref{SC9})  which must satisfy the relations: ${\cal V}_2={\cal
U}_{\,1}=0$.

\noindent\textbf{ Case 2}: ${\cal V}_3= {\cal V}_4= {\cal U}_{\,2}=0$.
 Now calculating
the value of ${\cal V}_3$ for system (\ref{S9}) we obtain:
$\displaystyle {\cal V}_3= 2^5 3^2 \sum_{j=0}^4 {\cal V}_{3j}
x^{4-j} y^j$, where ${\cal V}_{30}= -u^2+3v$, ${\cal V}_{31}=
-2uv+6r$, ${\cal V}_{32}= -2ur-v^2$, ${\cal V}_{33}= -2rv$, ${\cal
V}_{34}= -r^2$. Hence, condition ${\cal V}_3=0$ yields $r=v=u=0$
(then ${\cal V}_4= {\cal U}_2=0$) and, consequently we obtain
again system (\ref{SC9}).\EProof

We  examine system (\ref{SC9}). It is clear that via the
translation of the origin of coordinates at the point $(0,l/3)$,
we can consider parameter $l=0$ in the polynomial $q_2(x,y)$.
Thus, we must examine the  system
\be\label{sys_SC9:1} \dot x= a +
c x + d y +gx^2+ 2hxy+ky^2,\quad
\dot y= b + e x + f y + 2mxy+ny^2-x^3,\\
\ee for which $C_3(x,y)=x^4$. Therefore, there exists only one
direction for the possible invariant straight lines: $x=0$. In
this case  $U=1,\ V=0$ and from  (\ref{eq_g}) we obtain
\beq
&&A= B=C=0,\ D=g,\ E=2h,\ F=-gW+c,\\
&& Eq_7=k,\ Eq_9=-2hW+d,\ Eq_{10}=gW^2-cW+a. \eeq So, for having
the maximum number of invariant straight lines it is necessary
$k=h=d=0$ and then, we have ${\cal L}_1=0$, ${\cal L}_2= -2^{10}
3^4 n^2 x^2$. Hence, in order to reach the total multiplicity $9$,
by Lemma \ref{lm:L12}, it is necessary ${\cal L}_1= {\cal L}_2=
0$. Therefore, we obtain the additional condition: $n=0$.  This
provides the system:
\be\label{sys_SC9:2} \dot x= a + c x +
2gx^2,\quad
\dot y= b + e x + f y + 2mxy-x^3.\\
\ee We consider the homogenized system
\be\label{sys_SC9:2H} \dot
X= aZ^3 + c XZ^2 + 2gX^2Z,\quad \dot Y= bZ^3 + e XZ^2 + f YZ^2 +
2mXYZ-X^3, \ee
corresponding to system (\ref{sys_SC9:2}) and
calculate the following polynomial:
$$
  H(\ab,X,Y,Z)= \gcd({\cal G}_1,{\cal G}_2,{\cal G}_3)=
  2Z^2(aZ^2+cXZ+2gX^2).
$$
By Lemmas \ref{lm2} and \ref{lm3}, we obtain that system
(\ref{sys_SC9:2H}) has $4$ invariant straight lines (counted with
their multiplicities) $Z=0$ and $aZ^2+ cXZ+ 2gX^2= 0.$ So, for
having the total multiplicity $9$, by Lemmas \ref{lm2} and
\ref{lm3}, the polynomial $H(\ab,X,Y,Z)$ must be of the degree 8.
In order to find out the conditions to reach this situation we
shall calculate for  system (\ref{sys_SC9:2H}) the following
polynomials:
\be\label{S} \ba{ll}
{\cal G}_1/H=&  Y\tilde U(X,Z)+ \tilde V(X,Y)=   T(X,Y,Z),\\
{\cal G}_3/H=& 12(Z^2a+cZX+2gX^2)^3Z^2= 12 Z^2S^2(X,Z)^3, \ea \ee
where \beq
\tilde U(X,Z)&=&Z\left[4m(g-m)X^2+4f(g-m)XZ -(2am+f^2-fc)Z^2\right],\\
\tilde V(X,Z)&=&2(g+m)X^4 +(f+2c)X^3Z+ (3a+2ge-2me)X^2Z^2\\
   &&+(-2mb+4gb-fe)XZ^3-(ae-bc +fb)Z^4.
\eeq
By Lemma \ref{lm4}, in order to have the
maximum number of invariant straight lines it is necessary
$T(X,Y,Z)\mid Z^2S^2(X,Z)^3$. However, the second polynomial does
not depend on the variable $Y$, so we obtain  condition $\tilde
U(X,Z)= 0$, i.e. \be\label{cond:1a} m(g-m)=0,\quad f(g-m)=0,\quad
2am+f^2-fc=0. \ee Thus, we have to consider two cases: $g-m\ne0$
and $g-m=0$.

\noindent\textbf{ Case $g-m\ne0$.} Then, by (\ref{cond:1a}), we obtain
$m=f=0$, $g\ne0$, and then
\beq
T(X,Z)&=&2g X^4 +2cX^3Z+ (3a+2ge)X^2Z^2 +4gbXZ^3-(ae-bc)Z^4,\\
S^2(X,Z)&=&Z^2a+cZX+2gX^2. \eeq Since $g\ne0$, for obtaining the
maximum number of invariant straight lines it is necessary that
polynomial $T$ be proportional to polynomial $S^2$.

By Lemma \ref{Trudi:2}, we  have $ R^{(3)}_{X} \left(T,S^2\right)=
0,$ $R^{(2)}_{X} \left(T,\, S^2\right)= -4g^2 (2ag-c^2+4eg^2)^2
Z^4= 0$, and since $g\ne0$, we can consider $g=1$ (via the
transformation $x=x_1$, $y=y_1/g$ and $t=t_1/g$ to system
(\ref{sys_SC9:2})). Thus, we obtain $2a=c^2-4e$, and by setting
$c=2p,$ $ a=2(p^2-e)$, we get $ R^{(1)}_{X}\left(T,\, S^2\right)=
512(pe-p^3+b)^3Z^9=0. $ Therefore, the condition $b=p^3-pe$ holds.
Then we obtain: \beq
{\cal G}_1 &=& 8Z^2(Z^2p^2+XZp-eZ^2+X^2)^3,\\
{\cal G}_2&=& 48Z^3(2X+Zp)(Z^2p^2+XZp-eZ^2+X^2)^3,\\
{\cal G}_3&=& 384Z^4(Z^2p^2+XZp-eZ^2+X^2)^4.
\eeq Polynomials
${\cal G}_i$ $(i=1,2,3)$ have the common factor of degree 8, but
system (\ref{sys_SC9:2}) becomes degenerated:
$$
\dot x= 2(p^2-e+px+x^2),\quad
\dot y= (p-x)(p^2-e+px+x^2).
$$
So, in the case $g-m\ne0$ the system cannot possess invariant
straight lines with total multiplicity $9$.

\noindent\textbf{ Case $g-m=0$.} We take $m=g$. Then, conditions
(\ref{cond:1a}) yield $2ag+f(f-c)=0$, and we need to examine two
subcases: $g\ne0$ and $g=0$.

\noindent\textbf{ Subcase $g\ne0$.} Then, via the transformation
$x=x_1$, $y=y_1/g$ and $t=t_1/g$ to system (\ref{sys_SC9:2}), we
can consider $g=1$ (hence, $m=1$). Therefore, by setting $f=2p$,
we obtain the condition $a=p(c-2p)$. Thus, we can calculate:
$$
\ba{l}
R^{(3)}_{X}\left(T,\, S^2\right)=8(c-p)Z,\\
R^{(2)}_{X}\left(T,\, S^2\right)=64(p-c)(-ep+p^3+b)Z^4,\\
R^{(1)}_{X}\left(T,\,
S^2\right)=128(-ep+p^3+b)^2(4b+20p^3-c^3-24cp^2-4ep+9c^2p)Z^9.
\ea
$$
Hence, $c=p$ and then
$$
R^{(3)}_{X}\left(T,\, S^2\right)= R^{(2)}_{X}\left(T,\,
S^2\right)=0,\quad
 R^{(1)}_{X}\left(T,\, S^2\right)=512(-ep+p^3+b)^3Z^9=0,
$$
that implies $b=pe-p^3$. In this case we obtain: \beq
{\cal G}_1 &=& -2Z^2(X+Zp)^3(-2X+Zp)^3,\\
{\cal G}_2&=& -18Z^3X(X+Zp)^3(-2X+Zp)^3,\\
{\cal G}_3&=& 24Z^4(X+Zp)^4(-2X+Zp)^4. \eeq
i.e. ${\cal G}_1\mid
{\cal G}_2$ and ${\cal G}_1\mid {\cal G}_3$ but the system becomes
degenerated:
$$
\dot x= (p+x)(2x-p),\quad
\dot y= -(p+x)(p^2-px-e-2y+x^2).
$$

\noindent\textbf{ Subcase $g=0$.} Then, conditions (\ref{cond:1a})
yield $ m=0,\ f(f-c)=0, $ and we have to examine 2 subcases: $f=0$
and $f\ne0$.

For $f=0$  system (\ref{sys_SC9:2H}) becomes \be\label{SS0} \dot
X= aZ^3 + c XZ^2,\quad \dot Y= bZ^3 + e XZ^2 -X^3, \ee and we
obtain
$$
  H(\ab,X,Y,Z)= \gcd({\cal G}_1,{\cal G}_2,{\cal G}_3)=
  2(Za+cX)Z^4,
$$
and
$$
\ba{ll}
{\cal G}_1/H=&  2cX^3 +3aX^2Z+(bc-ae)Z^3 =  \hat T(X,Z),\\
{\cal G}_3/H=& 12(Za+cX)^3Z^4= 12 Z^4\hat S^3(X,Z).
\ea
$$
Thus, for  obtaining the  maximum number of invariant straight
lines it is necessary that the polynomial $\hat T$ divides
$Z^4\hat S^3$. This implies $c\ne0$, otherwise the second
polynomial does not depend on $X$ and, hence, it would be
necessary $a=0$, that yields the degenerated system (\ref{SS0})
with $c=a=0$. So, $c\ne0$ and then the polynomial $\hat T$ must be
proportional to $\hat S^3$. As above we calculate
$$
\ba{l}
R^{(2)}_{X}\left(\hat T,\,\hat S^3\right)=36ac^3Z=0,\\
R^{(1)}_{X}\left(\hat T,\,\hat S^3\right)=-432ac^4(a^3-c^2ae+c^3b)Z^4=0,\\
R^{(0)}_{X}\left(\hat T,\,\hat S^3\right)=1728c^3(a^3-c^2ae+c^3b)^3Z^9=0,\\
\ea
$$
and, hence, by $c\ne0$ we obtain $a=b=0$ that also provides a
degenerated system.

For $f\ne0$, from conditions (\ref{cond:1a}), we obtain $f=c\ne0$,
and we can consider $c=1$ via the transformation $x=x_1$, $y=
y_1/c$ and $t= t_1/c$. Then, system (\ref{sys_SC9:2H}) becomes
$$
\dot X= aZ^3 + XZ^2,\quad \dot Y= bZ^3 + e XZ^2+YZ^2 -X^3,
$$
for which we have
$$
  H(\ab,X,Y,Z)= \gcd({\cal G}_1,{\cal G}_2,{\cal G}_3)=
  2(Za+X)^2Z^4,
$$
and
$$
{\cal G}_1/H=  3X^2-eZ^2 =  \bar T(X,Z),\quad
{\cal G}_3/H= 12(Za+X)^2Z^4= 12 Z^4\bar S^2(X,Z).
$$
We calculate the resultant of the polynomials $\bar T$ and $\bar
S^2$ and its  first subresultant:
$$
\ba{l}
R^{(1)}_{X}\left(\bar T,\,\bar S^2\right)=6aZ=0,\\
R^{(0)}_{X}\left(\bar T,\,\bar S^2\right)=-(e-3a^2)^2Z^4=0.\\
\ea
$$
Consequently, we obtain $a=e=0$, and after translation of the
origin of coordinates at the point $(0,-b)$, we obtain the system:
\be\label{sys_conf:13} \dot x=  x,\quad
\dot y= y -x^3,\\
\ee which  is not degenerated. For the respective homogenized
system we have $H(\ab,X,Y,Z)= \gcd({\cal G}_1,{\cal G}_2,{\cal
G}_3)= X^4Z^4.$ Thus, we obtain Figure 23. Since for system
(\ref{sys_conf:13}) the conditions ${\cal V}_1= {\cal V}_2=0$ hold
as well as the conditions ${\cal V}_3={\cal V}_4=0$, we conclude
that this system has invariant straight lines with total
multiplicity 9 of both potential configurations: $(3,3,1,1)$ and
$(3,2,2,1)$. Indeed, the system
\be\label{s:pert15} \dot x=
x-4\varepsilon^2x^3,\quad
\dot y= y-x^3 -3\varepsilon^2x^2y+9\varepsilon^4xy^2-9\varepsilon^6y^3,\\
\ee with the invariant straight lines $x=0$, $2\varepsilon x= \pm
1$, $x- \varepsilon^2 y= 0$, $\varepsilon x- 3\varepsilon^3 y= \pm
1$, $3\varepsilon^2 y \pm x= 0$, has the configuration
$(3,3,1,1)$; and the system
\be\label{s:pert16} \dot x=
x+6\varepsilon x^2+8\varepsilon^2x^3,\quad
\dot y= y+6\varepsilon xy-x^3+9\varepsilon^2x^2y+9\varepsilon^4xy^2-9\varepsilon^6y^3,\\
\ee with the invariant straight lines
$$
   \ba{l} x=0,\ 2\varepsilon x=-1, \ \varepsilon^2y=x, \ \ 3\varepsilon^2 y \pm x=0,\\
        4\varepsilon x=-1,\ \        3\varepsilon(\varepsilon^2
        y-x)=1,\ \
        \varepsilon(3\varepsilon^2 y+x)=-1,
\ea
$$
has the configuration $(3,2,2,1)$. As above, system
(\ref{sys_SC9:1}) has the maximum number of invariant straight
lines with the configuration given by Figure 23 if and only if
\be\label{cond:tr} h=k=m=n=g=0,\quad d=f=c-f=0,\quad a=0. \ee We
give the respective affine invariant conditions in order to be
able to distinguish this class of cubic systems directly in the
space $\R^{20}$ of all cubic systems. To get this goal we consider
the system:
\be\label{sys_SC9:1t1} \dot x_1= p^\tau_0 + p^\tau_1
+gx_1^2+2hx_1y_1+ky_1^2,\quad
\dot y_1= q^\tau_0 + q^\tau_1 -3\gamma x_1^2+ 2mx_1y_1+ny_1^2-x_1^3,\\
\ee
which is obtained from system (\ref{sys_SC9:1}) via the
translation $x=x_1+\gamma $, $y=y_1+\delta  $, where $(\gamma
,\delta  )$ is an arbitrary point of the phase plane of system
(\ref{sys_SC9:1}) and $p^\tau_i(a(\gamma ,\delta  ),x_1,y_1)$,
$q^\tau_i(a(\gamma ,\delta  ),x_1,y_1)$ $(i=0,1)$ are the
corresponding homogeneous polynomials of degree $i$ in $x_1$ and
$y_1$.

For system (\ref{sys_SC9:1t1}) we calculate the $GL$-comitants
$$
 {\cal N}_2 = -6x_1^4\left[(g+m)x_1+(h+n)y_1\right],\quad
 {\cal N}_3 =
 -12x_1^3\left[(g-2m)x_1^2+2(2h-n)x_1y_1+3ky_1^2\right].
$$
The values of these polynomials do not depend on the coordinates
of the arbitrary point $(\gamma,\delta)$ and, consequently
conditions ${\cal N}_2=0$ and ${\cal N}_3=0$ are affine invariant
conditions. These conditions yield $g=h=k=m=n=0$, i.e. all
quadratic coefficients vanish. Thus, system (\ref{sys_SC9:1t1})
becomes the system
$$
\ba{ll}
\dot x_1=& (a+c\gamma +d\delta  ) + cx_1 + dy_1,\\
\dot y_1=& (b +e\gamma +f\delta  -\gamma ^3) +(e-3\gamma ^2)x_1 +
f y_1 -3\gamma x_1^2-x_1^3,
 \ea
$$
for which we have:
$$
 {\cal N}_9 = 3x_1^4\left[ex_1^2+(f-c)x_1y_1-dy_1^2\right],\quad
 {\cal N}_{10} = -12ax_1^3+6(f-c)\gamma x_1^3-12d\delta  x_1^3.
$$
As above ${\cal N}_9=0$ is affine invariant condition and implies
$d=e=c-f=0$. Moreover, in this case the $GL$--comitant ${\cal
N}_{10}$ is independent of the coordinates of the arbitrary point
$(\gamma ,\delta  )$, i.e. it becomes a $T$--comitant and
condition ${\cal N}_{10}=0$ yields $a=0$. Thus, it was proved,
that for system (\ref{sys_SC9:1}) conditions (\ref{cond:tr}) are
equivalent to ${\cal N}_2= {\cal N}_3= {\cal N}_9= {\cal
N}_{10}=0$. Hence, taking into account Lemma \ref{lm:5} we
obtain:

\bth\label{th:23} We assume that for a cubic system (\ref{s_1})
the conditions ${\cal D}_1= {\cal D}_2= {\cal D}_3= 0$, $C_3\not=
0$ hold (i.e. there exists only one (real) infinite singular
point). Then, this system will possess the maximum number of
invariant straight lines (with total multiplicity $9$) if and only
if $ {\cal V}_{1} ={\cal N}_{2} = {\cal N}_{3}= {\cal N}_{9}=
{\cal N}_{10}= 0. $ Moreover, there exists a unique configuration
given in Figure $23$. \eth

\bigskip

\noindent\textbf{ 9. Cubic systems whose infinite point configuration
do not allow to possess invariant straight lines with total
multiplicity $9$}

The goal of this section is to prove that  all others classes of
cubic systems enumerated in the statement of Lemma \ref{lm:5}
cannot have invariant straight lines with total multiplicity $9$.

\bigskip

\noindent\textbf{ 9.1 Systems with 4 imaginary  simple roots of $C_3$}

If a cubic system has $4$ imaginary  infinite singular points via
a linear transformation can written into the form (see, Section
3): \be\label{NS:1} \ba{l}
    x'=p_0 +p_1(x,y)+p_2(x,y) +ux^3+(p+q+v)x^2y+rxy^2+qy^3, \\
    y'=q_0 +q_1(x,y)+q_2(x,y) -px^3+ux^2y+vxy^2+ry^3, \\
\ea \ee for which $C_3(x,y)= (px^2+qy^2) (x^2+y^2)$, $pq>0$. By
Lemma \ref{lm:1}, in order to apply the conditions  ${\cal V}_1=
{\cal V}_2= {\cal U}_{\,1}= 0$, it is sufficient to calculate only
the T--comitant: $\displaystyle {\cal V}_1= 16\sum_{j=0}^4
{\cal V}_{1j}x^{4-j}y^j$, and to examine the relation ${\cal
V}_{10}= 2u^2+ 3p^2+ 3pq= 0. $ Indeed, since $pq>0$ the relation
${\cal V}_{10}=0$ cannot be satisfied and, therefore the
conditions of Lemma \ref{lm:1} cannot be satisfied for system
(\ref{NS:1}).

We examine the conditions of Lemma \ref{lm:4}: ${\cal V}_3= {\cal
V}_4= {\cal U}_2=0$. For system (\ref{NS:1}) we have
$\displaystyle {\cal V}_3= 32 \sum_{j=0}^4 {\cal V}_{3j} x^{4-j}
y^j$, and we examine the following two coefficients:  ${\cal
V}_{30}= -3p^2- 3pq+ 3pv- u^2$, ${\cal V}_{34}= -6q^2- 6pq- 3qv-
r^2$. Thus, conditions ${\cal V}_{30}= {\cal V}_{34}= 0$ imply
$3pv= 3p^2+ u^2+ 3pq>0$, $3qv= -(6q^2 +r^2 +6pq)< 0$, and hence,
$9pqv^2< 0$ in contradiction with $pq>0$. So, the conditions of
Lemma \ref{lm:4} cannot be verified for system (\ref{NS:1}).

\bigskip

\noindent\textbf{ 9.2 Systems with 1 double and 2 simple real roots of
$C_3$}

As it was shown in Section 3 in this case the cubic system via a
linear transformation can be written into the form:
\be\label{NS:2} \ba{l}
    x'=p_0 +p_1(x,y)+p_2(x,y) +(u+1)x^3+(v-1)x^2y+rxy^2, \\
    y'=q_0 +q_1(x,y)+q_2(x,y)+ ux^2y+vxy^2+ry^3, \\
\ea \ee for which $C_3(x,y)= x^2 y (x-y)$. By Lemma \ref{lm:1}, we
calculate:
$$
  {\cal V}_1=16\sum_{j=0}^4{\cal V}_{1j}x^{4-j}y^j,\qquad
  {\cal V}_2=8\sum_{j=0}^2{\cal V}_{2j}x^{2-j}y^j,\qquad
\hbox{where}
$$
$$
\ba{llll}
{\cal V}_{10}&=u(3+2u),     & {\cal V}_{14}&=2r^2,\\
{\cal V}_{11}&=-2u+4uv+3v,  & {\cal V}_{20}&=-3v-2u,\\
{\cal V}_{12}&=-1-2v+4ru+3r+2v^2,  &{\cal V}_{21}&=6r+4v-2,\\
{\cal V}_{13}&=2r(-1+2v), & {\cal V}_{22}&=-2r.\\
\ea
$$
Consequently, condition ${\cal V}_{14}= 2r^2= 0$ implies $r=0$,
and then we obtain the following contradictory relations: $ {\cal
V}_{10}=u(2u+3)=0$, $4{\cal V}_{20}+3{\cal V}_{21}= -2(4u+3)= 0$.
So, the conditions of Lemma \ref{lm:1} cannot be satisfied for
system (\ref{NS:2}).

We examine the conditions of Lemma \ref{lm:4}. For  system
(\ref{NS:2}) we have $\displaystyle {\cal V}_3= 32 \sum_{j=0}^4
{\cal V}_{3j} x^{4-j} y^j$, where ${\cal V}_{30}= -u(3+u)$, ${\cal
V}_{32}= 2+ v- 2ru+ 3r- v^2$, ${\cal V}_{31}= 2u(2-v)$, ${\cal
V}_{33}= -2r (1+v)$, ${\cal V}_{34}= -r^2$. Thus, condition ${\cal
V}_{34}= 0$ implies $r=0$, and then we again obtain the following
contradictory relations: ${\cal V}_{32}= -(v+1)(v-2)=0$, ${\cal
V}_4= 2^{10} 3^2 v (v-1) x^2 y (y-x)= 0$. So, the conditions of
Lemma \ref{lm:4} cannot be satisfied for system (\ref{NS:1}).

\bigskip

\noindent\textbf{ 9.3 Systems with 1 real double and 2 imaginary
simple roots of $C_3$}

In this case the cubic system  via a linear transformation can be
written into the form: \be\label{NS:3} \ba{l}
    x'=p_0 +p_1(x,y)+p_2(x,y) +u x^3+(v+1)x^2y+rxy^2, \\
    y'=q_0 +q_1(x,y)+q_2(x,y)-x^3+ ux^2y+vxy^2+ry^3, \\
\ea \ee for which $C_3(x,y)=x^2(x^2+y^2)$. By Lemma \ref{lm:1}, in
order to apply the conditions  ${\cal V}_1= {\cal V}_2= {\cal
U}_{\,1}=0$ it is sufficient to calculate only the T--comitant:
$\displaystyle {\cal V}_1=16\sum_{j=0}^4{\cal
V}_{1j}x^{4-j}y^j$, and to observe the contradiction ${\cal
V}_{10}= 2u^2+3= 0. $

We examine the conditions of Lemma \ref{lm:4}. For system
(\ref{NS:3}) we have $\displaystyle {\cal V}_3= 32 \sum_{j=0}^4
{\cal V}_{3j} x^{4-j} y^j$, where ${\cal V}_{30}= -u^2 -3+ 3v$,
${\cal V}_{31}= -2uv -4u +6r$, ${\cal V}_{32}= -2ur- v^2- v+2$,
${\cal V}_{33}= -2r(v-1)$, ${\cal V}_{34}= -r^2$. Thus, condition
${\cal V}_{34}= 0$ implies $r=0$, and then we again obtain the
following contradictory relations: ${\cal V}_{32}= -(v+2)(v-1)=
0$, ${\cal V}_4= 2^{10} 3^2 v(v+1) x^2 (x^2+y^2) =0.$ So, the
conditions of Lemma \ref{lm:4} cannot be satisfied for system
(\ref{NS:3}).

\bigskip

\noindent\textbf{ 9.4 Systems with $2$ double  real roots of $C_3$}

In this case the cubic system  via a linear transformation
becomes: \be\label{NS:4} \ba{l}
    x'=p_0 +p_1(x,y)+p_2(x,y) +u x^3+q x^2y+rxy^2, \\
    y'=q_0 +q_1(x,y)+q_2(x,y)+ ux^2y+vxy^2+ry^3, \\
\ea \ee for which $ C_3(x,y)=(q-v)x^2y^2,\quad q-v\ne0. $ By Lemma
\ref{lm:1}, in order to apply the conditions  ${\cal V}_1= {\cal
V}_2= {\cal U}_{\,1}= 0$ we calculate: $\displaystyle {\cal
V}_1= 16\sum_{j=0}^4 {\cal V}_{1j} x^{4-j} y^j,$
$\displaystyle {\cal V}_2= 16 \sum_{j=0}^2 {\cal
V}_{2j}x^{2-j}y^j, $\ where ${\cal V}_{10}=2u^2,$\  ${\cal
V}_{11}=2u(q+v),$\ $ {\cal V}_{12}=4ur-q^2+4qv-v^2$, $ {\cal
V}_{13}=2r(q+v)$,\ ${\cal V}_{14}=2r^2,$\   ${\cal
V}_{20}=u(q-v),$\ $ {\cal V}_{21}=(v-q)(v+q),$\  ${\cal
V}_{22}=r(q-v).$ Consequently, conditions ${\cal V}_{10}={\cal
V}_{14}=0$ imply $u=r=0$, and then from $q-v\ne0$, we obtain the
following contradictory relations: $ {\cal V}_{21}= (v-q) (v+q)=
0$, ${\cal V}_{12}- 3{\cal V}_{21}= 2v(2q-v)=0. $ So, the
conditions of Lemma \ref{lm:1} cannot be satisfied for system
(\ref{NS:4}).

We examine the conditions of Lemma \ref{lm:4}. For system
(\ref{NS:4}) we have $\displaystyle {\cal V}_3= 32 \sum_{j=0}^4
{\cal V}_{3j} x^{4-j} y^j$, where ${\cal V}_{30}=-u^2$, ${\cal
V}_{31}= 2u(v-2q),$\ ${\cal V}_{32}=-2ur+2q^2-5qv+2v^2,$\ ${\cal
V}_{33}=2r(q-2v),$\ ${\cal V}_{34}=-r^2.$\ Thus, conditions ${\cal
V}_{30}={\cal V}_{34}=0$ imply $u=r=0$ and then by condition
$q-v\ne0$ we again obtain the following contradictory relations:
$$
  {\cal V}_{32}= (2q-v)(q-2v)=0, \quad {\cal V}_4=2^{10}3^2
  qv(q-v)^2x^2y^2=0.
$$
So, the conditions of Lemma \ref{lm:4} also cannot be verified for
system (\ref{NS:4}).

\bigskip

\noindent\textbf{ 9.5 Systems with $2$ double imaginary roots of
$C_3$}

In this case the cubic system  via a linear transformation goes
over to \be\label{NS:5} \ba{l}
    x'=p_0 +p_1(x,y)+p_2(x,y) +u x^3+3(v+1) x^2y+rxy^2 +y^3, \\
    y'=q_0 +q_1(x,y)+q_2(x,y) -x^3 + ux^2y+(v-1)xy^2+ry^3, \\
\ea \ee for which $ C_3(x,y)=(x^2+y^2)^2. $ By Lemma \ref{lm:1},
in order to apply the conditions  ${\cal V}_1= {\cal V}_2= {\cal
U}_{\,1}= 0$ it is sufficient to calculate only the T--comitant:
$\displaystyle {\cal V}_1= 16 \sum_{j=0}^4 {\cal V}_{1j}
x^{4-j }y^j$, and to observe the contradiction ${\cal V}_{10}=
u^2+3= 0.$

We examine the conditions of Lemma \ref{lm:4}. For system
(\ref{NS:5}) we have $\displaystyle {\cal V}_3=32\sum_{j=0}^4{\cal
V}_{3j}x^{4-j}y^j, $ and we examine only the following two
coefficients: ${\cal V}_{30}= -u^2-9+3v,$\ ${\cal V}_{34}=
-r^2-9-3v$. Therefore, we obtain the following contradictory
condition: ${\cal V}_{30}+ {\cal V}_{34}= -(u^2+r^2+18)= 0,$ and,
hence, the conditions of Lemma \ref{lm:4} do not hold for system
(\ref{NS:5}).


\begin{thebibliography}{99}

\bibitem{AL} {\sc J. Artes and J.Llibre},
{\emph On the number of slopes of invariant straight lines for
polynomial differential systems}, {J. of Nanjing University \textbf{
13} (1996), 143--149.
}
\bibitem{Art_Llibre} {\sc J. Artes, B. Gr\"unbaum and J.Llibre},
\emph{On the number of invariant straight lines for polynomial
differential systems}, Pacific Journal of Mathematics \textbf{  184},
(1998), 317--327.

\bibitem{Calin} {\sc Iu. F. Calin},
\emph{Private communication}, 2001.

\bibitem{Llib_Per} {\sc C. Cristofer, J. Llibre and J. V. Pereira},
\emph{Multiplicity of invariant algebraic curves}, preprint 2002.

\bibitem{Sib_Dini} {\sc Dang Din' Bik and  K. S. Sibirskii},
\emph{Affine classification of a cubic differential system},
Studies in algebra, mathematical analysis, and their aplications,
Izdat. ``\c Stiin\c ta``, Kishinev, 1977, 43--52 (Russian).

\bibitem{Druzhkova} {\sc T.A. Druzhkova},
\emph{Quadratic differential systems with algebraic integrals},
Qualitative theory of differential equations, Gorky Universitet
\textbf{ 2} (1975), 34--42 (Russian).


\bibitem{Gr_Yng} {\sc J.H. Grace and A. Young},
\emph{The algebra of invariants}, New York: Stechert, 1941.

\bibitem{Trudi} {\sc A.S. Householder},
\emph{Bigradients and the problem of Routh and Hurwitz},
SIAM Review \textbf{ 10} (1968), 166--178.

\bibitem{Lyubim1} {\sc R.A. Lyubimova},
\emph{On some differential equation possesses invariant lines},
Differential and integral eequations, Gorky Universitet \textbf{ 1}
(1977), 19--22 (Russian).

\bibitem{Lyubim2} {\sc R.A. Lyubimova},
\emph{On some differential equation possesses invariant lines},
Differential and integral equations, Gorky Universitet \textbf{  8}
(1984), 66--69 (Russian).

\bibitem{Olver} {\sc P.J. Olver},
\emph{Classical Invariant Theory}, (London Mathematical Society
student texts: \textbf{ 44}), Cambridge University Press, 1999.

\bibitem{Popa1} {\sc M.N. Popa},
\emph{The number of comitants that are involved in determining the
number of integral lines of a cubic differential system}, Izv.
Akad. Nauk Moldav. SSSR Mat. \textbf{ 1} (1990), 67--69 (Russian).

\bibitem{Popa2} {\sc M.N. Popa},
\emph {Application of invariant processes to the study of
homogeneous linear particular integrals of a differential system},
Dokl. Akad. Nauk SSSR \textbf{ 317} (1991), 834--839 (Russian);
translation in Soviet Math. Dokl. \textbf{ 43} (1991), 550--555.

\bibitem{Popa3} {\sc M.N. Popa},
\emph {Conditions for the maximal multiplicity of an integral line
of a differential system with homogeneities of $m^{\hbox{th}}$
order}, Izv. Akad. Nauk Respub. Moldova, Mat. \textbf{ 1} (1992),
15--17 (Russian).

\bibitem{Popa_Sib1} {\sc M.N. Popa and K. S. Sibirskii},
\emph {Conditions for the existence of a homogeneous linear
partial integral of a differential system}, Differentsial'nye
Uravneniya \textbf{ 23} (1987), 1324--1331 (Russian).

\bibitem{Popa_Sib2} {\sc M.N. Popa and K. S. Sibirskii},
\emph {Conditions for the prezence of a  nonhomogeneous linear
partial integral in a quadratic differential system}, Izv. Akad.
Nauk Respub. Moldova, Mat. \textbf{ 3} (1991), 58--66 (Russian).

\bibitem{Popa_Sib3} {\sc M.N. Popa and K. S. Sibirskii},
\emph {Integral line of a general quadratic differential system},
Izv. Akad. Nauk Moldav.SSR, Mat. \textbf{ 1} (1991), 77--80 (Russian).

\bibitem{Dana_Vlp1} {\sc D. Schlomiuk and  N. Vulpe},
\emph{Geometry of quadratic differential systems in the
neighbourhood of the line at  infinity},  Report no. 2701, Centre
de recherches math\'ematiques et  D\'epartement de Math\'ematiques
et de Statistiques, Universit\'e de Montr\'eal, 2001.

\bibitem{Dana_Vlp2} {\sc D. Schlomiuk and  N. Vulpe},
\emph{Planar Quadratic Differential Systems with invariant lines},
preprint, 2001.

\bibitem{Sib1} {\sc K.S. Sibirskii},
\emph{Introduction to the algebraic theory of invariants of
differential equations}, Translated from the Russian, Nonlinear
Science: Theory and Applications, Manchester University Press,
Manchester, 1988.

\bibitem{Sib2} {\sc K.S. Sibirskii},
\emph{Conditions for the presence of a straight integral curve of
a quadratic differential system in the case of a center or a
focus}, Mat. Issled. No. 106, Differ. Uravneniya i Mat. Fizika,
1989, 114--118 (Russian).

\bibitem{Sokulski} {\sc J. Sokulski},
\emph{On the number of invariant lines for polynomial vector
fields}, Nonlinearity \textbf{  9} (1996), 479--485.

\bibitem{ZX} {\sc Zhang Xiang},
\emph{Number of integral lines of polynomial systems of degree
three and four}, J. of Nanjing University, Math. Biquartely \textbf{
10} (1993), 209--212.

\end{thebibliography}
\end{document}